\input amstex
\input epsf.sty
\input cyracc.def
\input amssym

\documentstyle{amsppt}
\document
\magnification=1200
\NoBlackBoxes
\TagsOnRight
\nologo

\newfam\cyifam
\font\tencyi=wncyi10
\font\sevencyi=wncyi7
\font\fivecyi=wncyi5
\def\cyi{\fam\cyifam\tencyi\cyracc}
\textfont\cyifam=\tencyi \scriptfont\cyifam=\sevencyi
  \scriptscriptfont\cyifam=\fivecyi

\define\lon{\longrightarrow}
\define\rar{\rightarrow}

\define\OM{\Omega^1 {\Cal M}}
\define\CP{\bold{P}}
\define\p{{\partial}}
\define\C{{\Bbb C}}
\define\ot{\otimes}
\define\mbox{\text}
\define\Id{\mbox{Id}}
\define\f{{\Cal O}}
\define\em{\it}
\define\cT{{\Cal T}}
\define\cM{{\Cal M}}

\define\cW{{\Cal W}}
\define\cV{{\Cal V}}
\define\cF{{\Cal F}}
\define\cE{{\Cal E}}
\define\cTM{{\Cal TM}}

\define\al{\alpha}
\define\be{\beta}
\define\ga{\gamma}
\define\var{\varepsilon}
\define\la{\lambda}
\define\om{\omega}
\define\dal{\dot{\alpha}}
\define\dbe{\dot{\beta}}
\define\dga{\dot{\gamma}}
\define\dmu{\dot{\mu}}
\define\dnu{\dot{\nu}}
\define\dde{\dot{\delta}}

\define\te{\tilde{e}}
\define\ts{\tilde{s}}
\define\tnabla{\tilde{\nabla}}

\define\tvar{\tilde{\varepsilon}}
\define\bvar{\bar{\varepsilon}}


\hyphenation{Schle-singer}


\centerline{\bf SEMISIMPLE FROBENIUS (SUPER)MANIFOLDS}

\medskip

\centerline{\bf AND QUANTUM COHOMOLOGY OF $\bold{P}^r$}

\medskip

\centerline{\bf Yu. I. Manin}

\medskip

\centerline{\it Max--Planck--Institut f\"ur Mathematik, Bonn, Germany}

\medskip

\centerline{\bf S. A. Merkulov}

\medskip

\centerline{\it Department of Mathematics, University of Glasgow, UK}

\bigskip

\line{\hfill {\cyi Posvyawaet\-sya O.~A.~Ladyzhensko\u\i.}}

\bigskip

{\bf Abstract.} We introduce and study a superversion
of Dubrovin's notion of semisimple Frobenius manifolds.
We establish a correspondence between semisimple
Frobenius (super)manifolds and special solutions to the (supersymmetric)
Schlesinger equations.
Finally, we calculate the Schlesinger initial conditions for solutions
describing quantum cohomology of projective spaces.

\bigskip

\centerline{\bf 0. Introduction}

\medskip

B.~Dubrovin introduced the notion of Frobenius manifold
and made an extensive study of it in [D]. Roughly
speaking, it is a triple $(M,g,\circ )$ where $M$
is a manifold, $g$ is a flat Riemannian metric on it,
and $\circ$ is a $\Cal{O}_M$--linear commutative and associative
multiplication on the tangent sheaf $\Cal{T}_M$,
with compatibility conditions (see 1.1.1 below).
An important class of examples is supplied by
{\it quantum cohomology} (see [KM].) Actually,
quantum cohomology furnishes versions of Dubrovin's
definition in which $M$ may be supermanifold,
or even a formal supermanifold.

\smallskip

An important subclass of Frobenius manifolds
consists of semisimple ones. This means that
tangent spaces with $\circ$--multiplication
are semisimple algebras. This is possible only
if $M$ has no odd coordinates, by purely formal
reasons.

\smallskip

In [Ma1], p.~41, one of the authors suggested that it would
be interesting to construct a natural superization of
the notion of semisimple Frobenius manifolds.
This is one of the goals of this paper.
To avoid any misunderstanding, 
we must stress that semisimple Frobenius supermanifolds
in the sense of this work {\it are not} Frobenius
manifolds in the category of supermanifolds  in the sense
of [KM] or [Ma1].

\smallskip

Our extension is based upon
Dubrovin's theory
of semisimple Frobenius manifolds reducing
their classification to that of isomonodromic deformations: see [D].
This reduction exists in two versions. The first version leads
to the deformation of a flat connection on a vector
bundle on $\bold{P}^1$
having two singular points, a regular and an irregular one
(see [D] and [S].) The second one deals with
connections having only regular singularities.
The two versions are related by the formal Laplace transform,
as was explained in [KM].

\smallskip

In the classical paper [Sch], L.~Schlesinger constructed
the universal deformation space of the connections with regular
singularities on $\bold{P}^1$ (see [Mal3] for a modern treatment.)
Schlesinger's equations govern the dependence
of the universal connection on the deformation parameters.
Semisimple Frobenius manifolds correspond to some
solutions to Schlesinger's equations, with the structure
group reduced to the orthogonal one, and  supplied
with an additional piece of data. These solutions
are called here ``strict special'' ones.

\smallskip

Our supeversion of the Dubrovin theory
includes a superization of Schlesinger's equations,
of the notion of strict special solutions,
and of the correspondence between them
and semisimple Frobenius manifolds
briefly described above.
This is one of the arguments for the naturality
of our definition. An additional detail
is that the (rather mysterious) structural action of the
Virasoro algebra on any semisimple Frobenius manifold
is now replaced by that of the Neveu--Schwarz superalgebra.
We hope that super--Schlesinger equations may present
an independent interest. 

\smallskip

Since the structure of the superversion is closely parallel
to that of the pure even one, we start with a report on the theory
of semisimple Frobenius manifolds
and with its application to the quantum cohomology
of projective spaces.

\smallskip

Quantum cohomology of a projective algebraic manifold $V$
is the pair $(H^*(V),\Phi )$
consisting of the usual cohomology space, say, with complex
coefficients, and  {\it the potential} $\Phi$, a formal function
on the cohomology space whose Taylor coefficients are numerical
invariants of $V$ counting the number of parametrized
rational curves subject to certain incidence conditions:
cf. e.~g. [KM], [M1] and [BM]. Its third derivatives 
form the structure constants
of {\it the quantum cohomology algebra of} $V.$ 

\smallskip 

For projective spaces, it turns out to be semisimple, and 
we characterize the relevant special solutions
by their initial conditions.

\smallskip

The paper is structured as follows. \S 1 contains an overview
of the basic facts about Frobenius manifolds. 
Omitted proofs can be found in [D] and [Ma1]. In \S 2, we discuss the
Schlesinger equations and interrelations between them
and Frobenius manifolds.
The version we explain here is taken from [Ma1]; for a closely
related treatment see [H]. An important amelioration
is the Theorem 1.14 and the related notion of strict speciality. 
\S 3 is devoted to the quantum
cohomology of projective spaces; see [Ma2] and [Ma1], \S III.5
for the special case of $\bold{P}^2.$
In \S 4, we supersymmetrize the notion of semisimple Frobenius
manifolds. Finally, in \S 5 we discuss a correspondence 
between semisimple Frobenius supermanifolds and special solutions of 
supersymmetric Schlesinger equations.

\newpage

\centerline{\bf \S 1. Frobenius manifolds}

\bigskip

{\bf 1.1. Frobenius manifolds.} Throughout this paper, we work
in the category of complex manifolds $M.$ A metric on $M$
is an even symmetric pairing 
$g: S^2(\Cal{T}_M)\to\Cal{O}_M$, inducing an isomorphism
$g^{\prime}:\Cal{T}_M\to\Cal{T}^*_M.$ Here $\Cal{O}_M$
is the structure sheaf, and $\Cal{T}_M$ is the tangent sheaf.

\smallskip

An affine flat structure on $M$ is a subsheaf $\Cal{T}_M^f\subset\Cal{T}_M$
of linear spaces of pairwise commuting vector fields,
such that $\Cal{T}_M=\Cal{O}_M\otimes_{\bold{C}}\Cal{T}_M^f.$
Sections of $\Cal{T}^f_M$ are called flat vector fields.
The metric $g$ is compatible with the structure
$\Cal{T}^f_M$, if $g(X,Y)$ is constant for flat $X,Y.$

\smallskip

An affine flat structure can be 
equivalently described by a complete atlas whose transition functions are affine linear, because for a maximal commuting set of linearly independent
vector fields $(X_a)$ one can find local coordinates such
that $X_a=\partial /\partial x^a$, and they are defined up to
a constant shift.
If a metric $g$ is compatible with an affine flat structure,
it is flat in the sense of the formalism
of Riemannian geometry. The parallel transport
endows $\Cal{T}^f_M$ with the structure of local system.

\medskip

\proclaim{\quad 1.1.1. Definition} Let $M$ be a manifold.
Consider a triple $(\Cal{T}^f_M, g, A)$ consisting of an
affine flat structure, a compatible metric, and an even symmetric tensor
$A: S^3(\Cal{T}_M)\to\Cal{O}_M.$

\smallskip

Define an $\Cal{O}_M$--bilinear symmetric multiplication $\circ =\circ_{A,g}$
on $\Cal{T}_M$:
$$
\Cal{T}_M\otimes\Cal{T}_M\to S^2(\Cal{T}_M)\,{\overset{A^{\prime}}\to\rightarrow}\,
\Cal{T}_M^*\,{\overset{g^{\prime}}\to\rightarrow}\,\Cal{T}_M:\ 
X\otimes Y\to X\circ Y
\eqno{(1.1)}
$$
where prime denotes a partial dualization, or equivalently,
$$
A(X,Y,Z)=g(X\circ Y,Z)=g(X,Y\circ Z).
\eqno{(1.2)}
$$
This means that the metric is invariant with respect to the
multiplication.

\smallskip

a). $M$ endowed with this structure is called a pre--Frobenius 
manifold.

\smallskip

b). A local potential $\Phi$ for $(\Cal{T}^f_M,A)$ is
a local even function such that for any flat local tangent fields
$X,Y,Z$
$$
A(X,Y,Z)=(XYZ)\Phi .
\eqno{(1.3)}
$$
A pre--Frobenius manifold is called potential one, if
$A$ everywhere locally admits a potential.

\smallskip

c). A pre--Frobenius manifold is called associative, if
the multiplication $\circ$ is associative.

\smallskip

d). A pre--Frobenius manifold is called Frobenius,
if it is simultaneously potential and associative.
\endproclaim

\medskip

If a potential $\Phi$
exists, it is unique up to adding a quadratic polynomial
in flat local coordinates.

\smallskip

In flat local coordinates $(x^a)$ (1.3) becomes
$A_{abc}=\partial_a\partial_b\partial_c\Phi$, and (1.2)
can be rewritten as
$$
\partial_a\circ\partial_b=\sum_cA_{ab}{}^c\partial_c,
\eqno{(1.4)}
$$
where 
$$
A_{ab}{}^c:=\sum_eA_{abe}g^{ec},\ (g^{ab}):=(g_{ab})^{-1}.
$$
Furthermore,
$$
(\partial_a\circ\partial_b)\circ\partial_c=
\left(\sum_e A_{bc}{}^e\partial_e\right)\circ\partial_c=
\sum_{ef}A_{ab}{}^eA_{ec}{}^f\partial_f,
$$
$$
\partial_a\circ (\partial_b\circ\partial_c)=
\partial_a\circ\sum_e A_{bc}{}^e\partial_e
=\sum_{ef}A_{bc}{}^eA_{ae}{}^f\partial_f=
$$
$$
=\sum_{ef}A_{bc}{}^eA_{ea}{}^f\partial_f.
\eqno{(1.5)}
$$

\smallskip

Comparing the coefficients of $\partial_f$ in (1.5),
lowering the superscripts and expressing $A_{abc}$
through a potential, we finally see that the notion
of the Frobenius manifold is a geometrization of the
following  highly non--linear and
overdetermined system of PDE:
$$
\forall a,b,c,d:\quad
\sum_{ef}\Phi_{abe}g^{ef}\Phi_{fcd}=
\sum_{ef}\Phi_{bce}g^{ef}\Phi_{fad}.
\eqno{(1.6)}
$$
They are called Associativity Equations, or
WDVV (Witten--Dijkgraaf--Verlinde--Verlinde) equations.

\smallskip

Following B.~Dubrovin, we will now express (1.6) as a flatness condition.

\smallskip

{\bf 1.2. The first structure connection.} Let $(M,g,A)$ be a pre--Frobenius
manifold (we omit $\Cal{T}^f_M$ in the notation, since it can
be reconstructed from $g$.)
Define the following objects:

\smallskip

a). The connection $\nabla_0: \Cal{T}_M\to\Omega_M^1\otimes\Cal{T}_M$
well determined by the condition that flat vector fields are
$\nabla_0$--horizontal.

\smallskip

Denote its covariant derivative along a vector field $X$ by
$$
\nabla_{0,X}(Y)=i_X(\nabla_0(Y)),\ i_X(df\otimes Z)=Xf\otimes Z.
$$

b). A pencil of connections depending on an even parameter 
$\lambda$:
$$
\nabla_{\lambda}: \Cal{T}_M\to\Omega_M^1\otimes\Cal{T}_M:\
\nabla_{\lambda ,X}(Y):=\nabla_{0,X}(Y)+\lambda X\circ Y.
\eqno{(1.7)}
$$ 
We will call $\nabla_{\lambda}$ the first structure connection of
$(M,g,A)$.

\smallskip

In flat coordinates (1.7) reads:
$$
\nabla_{\lambda ,\partial_a}(\partial_b )=
\lambda\sum_cA_{ab}{}^c\partial_c=
\lambda\,\partial_a\circ\partial_b=
(-1)^{ab}\lambda\partial_b\circ\partial_a =(-1)^{ab}
\nabla_{\lambda ,\partial_b}(\partial_a ).
$$
Therefore $\nabla_{\lambda}$
has vanishing torsion for any $\lambda$. In particular,
$\nabla_0$ is the Levi--Civita connection for $g$.

\smallskip

\proclaim{\quad 1.2.1. Theorem} Let $\nabla_{\lambda}$
be the structure connection of the pre--Frobenius
manifold $(M,g,A).$ Put $\nabla_{\lambda}^2=
\lambda^2R_2+\lambda\,R_1$ (there is no constant term
since $\nabla_0^2=0.)$ Then

\smallskip

a). $R_1=0\ \Longleftrightarrow\ (M,g,A)$ is potential.

\smallskip

b). $R_2=0\ \Longleftrightarrow\ (M,g,A)$  is associative.

\smallskip

Therefore $(M,g,A)$ is Frobenius, iff $\nabla_{\lambda}$
is flat.
\endproclaim

\smallskip

This can be proved by the direct computation.

\medskip

{\bf 1.3. Identity.} Let $(M,g,A)$ be a pre--Frobenius manifold.
A vector field $e$ on $M$ is called identity, if
$e\circ X=X$ for all $X.$

\smallskip

If $e$ exists at all, it is uniquely defined by $\circ$, hence by
$g$ and $A.$

\smallskip

Conversely, given $A$ and $e$, there can exist at most one metric
$g$ making $(M,g,A)$ a pre--Frobenius manifold with this identity:
$$
g(X,Y)=A(e,X,Y).
$$
This follows from (1.2). If $A$ has a potential $\Phi$,
this translates into a non--homogeneous linear differential
equation for $\Phi$ supplementing the Associativity Equations (1.6):
$$
\forall \ \roman{flat}\ X,Y,\quad eXY\Phi =g(X,Y).
\eqno{(1.8)}
$$
In most (although not all) important examples $e$ itself is {\it flat}.
If this is the case, one can everywhere locally find a flat
coordinate system $(x^0,\dots ,x^n)$ such that 
$e=\partial/\partial x^0=\partial_0$, and (1.8) becomes
$$
\forall\, a,b,\quad \Phi_{0ab}=g_{ab}.
\eqno{(1.9)}
$$
Since all $g_{ab}$ are constants, we get the following result.

\smallskip

On a potential pre--Frobenius
manifold with flat identity $e=\partial_0$ (in a flat coordinate
system) we have modulo terms of degree $\le 2$:
$$
\Phi (x^0,\dots ,x^n)=
\frac{1}{2}x^0\left(\sum_{a,b\ne 0}g_{ab}x^ax^b+\sum_{a\ne 0}g_{0a}x^0x^a+
\frac{1}{3}\,g_{00}(x^0)^2\right)+
\Psi (x^1,\dots ,x^n).
\eqno{(1.10)}
$$

\smallskip

The metric $g$ identifies $\Cal{T}_M$
and $\Cal{T}_M^*$. We will call {\it the co--identity} and denote
$\varepsilon$ the 1--form which is the image of $e$. 
$\varepsilon$ is defined by
$$
\forall X\in \Cal{T}_M,\quad i_X(\varepsilon)=g(X,e).
$$
If $(x^a)$ is a local coordinate system,then
$$
\varepsilon =\sum_a dx^ag(\partial_a,e).
$$
Finally, if $e$ and $(x^a)$ are flat, then $g(\partial_a,e)$
are constant, and
$$
\varepsilon =d\eta,\ \eta=\sum_ax^ag(\partial_a,e).
\eqno{(1.11)}
$$
\medskip

{\bf 1.4.  Euler field.} We will say that a vector field $E$
on a manifold with flat metric $(M,g)$ is {\it conformal},
if $\roman{Lie}_E(g)=Dg$ for some constant $D.$ In other words,
for all vector fields $X,Y$ we have
$$
E(g(X,Y))-g([E,X],Y)-g(X,[E,Y])=Dg(X,Y).
\eqno{(1.12)}
$$
It follows that in flat coordinates we have 
$E=\sum_aE^a(x)\partial_a$ where $E^a(x)$ are polynomials of degree
$\le 1.$ In fact, $E$ is a sum of infinitesimal rotation,
dilation and constant shift. Hence 
$[E,\Cal{T}^f_M]\subset \Cal{T}_M^f.$
Moreover, the operator
$$
\Cal{V}:\ \Cal{T}^f_M\to \Cal{T}^f_M,\ \Cal{V}(X):=
[X,E]-\frac{D}{2}\,X
$$
is skewsymmetric: 
$$
\forall\ \roman{flat\ }X,Y:\ g(\Cal{V}(X),Y)+g(X,\Cal{V}(Y))=0.
$$ 
\smallskip
\proclaim{\quad 1.4.1. Definition} Let $E$ be a vector field on
a pre--Frobenius manifold $(M,g,A).$ It is called
an Euler field, if it is conformal, and $\roman{Lie}_E(\circ )=d_0\circ$
for some constant $d_0$, that is, for all vector fields $X,Y$,
$$
[E,X\circ Y]-[E,X]\circ Y-X\circ [E,Y]=d_0X\circ Y.
\eqno{(1.13)}
$$

\endproclaim

\smallskip

Clearly, any scalar multiple
of an Euler field is also an Euler field.
One can use this remark in order to normalize $E$
by requiring that some non--vanishing eigenvalue becomes
one. A convenient choice is often $d_0=1$, if we
have reasons to restrict ourselves to the $d_0\ne 0$ case. 

\smallskip

\proclaim{\quad 1.4.2. Proposition} Let $E$ be a conformal vector field on
a Frobenius manifold $(M,g,\Phi ).$ Then $E$ is Euler, iff
$$
E\Phi=(d_0+D)\Phi+\roman{a\ quadratic\ polynomial\ in\ flat\ coordinates} .
$$
\endproclaim

\smallskip

\medskip

\medskip

{\bf 1.4.3. Case of semisimple $\roman{ad}\,E.$}  We will call the set of eigenvalues
of $-\roman{ad}\,E$ on $\Cal{T}^f_M$, together with $d_0$
and $D$, {\it the spectrum} of $E.$ We will say that $E$
is {\it semi--simple}, if $\roman{ad}\,E,$ acting on flat fields, is.
For semisimple $E$ we can construct many homogeneous elements
of $\Cal{O}_M(*)$ and $\Cal{T}_M(*)$ explicitly.

\smallskip

Let $(\partial_a)$ be a local basis of $\Cal{T}^f_M$
such that 
$$
[\partial_a,E]=d_a\partial_a
\eqno{(1.14)}
$$
where $(d_a)$ form a part of the spectrum of $E$.
Putting $E=\sum E^a(x)\partial_a,$ we find from
(1.14) that $\partial_aE^b=\delta_a^b\,d_a.$
Hence if $\partial_a=\partial /\partial x^a$,
we have
$$
E=\sum_{a:\,d_a\ne 0}(d_ax^a+r^a)\partial_a+\sum_{b:\,d_b=0}r^b\partial_b.
$$
By shifting $x^a$, we can make $r^a=0$ for $d_a\ne 0.$
Multiplying $x^b$ by a constant, we can make $r^b=0$ or $1$
for $d_b=0$. So finally we can choose local flat
coordinates in such a way that
$$
E=\sum_{a:\,d_a\ne 0}d_ax^a\partial_a+\sum_{\roman{some}\,b:\,d_b=0}
\partial_b. 
\eqno{(1.15)}
$$
Clearly, $E$ assigns definite degrees to the following local
functions:
$$
Ex^a=d_ax^a\ \roman{for}\ d_a\ne 0;\ E\,\roman{exp}\,x^b=\roman{exp}\,x^b\
\roman{or}\ 0\ \roman{for}\ d_b = 0.
\eqno{(1.16)}
$$
\smallskip

Assume now that $M$ has an identity $e$. From (1.13) we get
$$
[e,E]=d_0e.
\eqno{(1.17)}
$$
Hence our notation for the spectrum will be consistent, if in the case
of flat $e$ we put $e=\partial_0$, and otherwise do not
use $0$ as one of the subscripts in (1.14).

\smallskip

Formula (1.12) in the basis (1.13) becomes
$$
\forall \, a,b:\ g(d_a\partial_a,\partial_b)+
g(\partial_a,d_b\partial_b)=Dg_{ab}
$$
that is,
$$ 
(d_a+d_b-D)g_{ab}=0.
\eqno{(1.18)}
$$ 
In particular, $g(e,e)=0$ unless
$D=2d_0.$

\smallskip

{\bf 1.5.  Extended structure connection.} Let $M$ be
a pre--Frobenius manifold with a conformal vector field $E$.
Put $\widehat{M}:=M\times (\bold{P}^1_{\lambda}\setminus
\{0,\infty\} ),$ where $\bold{P}^1_{\lambda}$ is the completion
of $\roman{Spec}\,\bold{C}[\lambda ,\lambda^{-1}].$
Furthermore, put $\widehat{\Cal{T}}=\roman{pr}_M^*(\Cal{T}_M).$
If $X$ is a vector field on $M$, it may be lifted to
$\widehat{M}$ in two different guises: as a vector field annihilating
$\lambda$,
denoted again $X$, and as a section of $\widehat{\Cal{T}}$,
then denoted $\widehat{X}$. 

\smallskip

Choose a constant $d_0$ and put $\Cal{E}:=E-d_0\lambda\dfrac{\partial}{\partial\lambda}\in
\Cal{T}_{\widehat{M}}.$ 
Clearly, $\widehat{X}$ for flat $X$ span
 $\widehat{\Cal{T}}$, whereas flat $X$ and $\Cal{E}$ span
 $\Cal{T}_{\widehat{M}}$,
{\it provided $d_0\ne 0$}, which we will assume.

\smallskip

\proclaim{\quad 1.5.1. Definition} Let $M$ be a pre--Frobenius manifold
with a conformal field 
$E,$ and $d_0$ a non--zero constant. 
The extended
structure connection for $M$
is the connection $\widehat{\nabla}$ on the sheaf $\widehat{\Cal{T}}$
on $\widehat{M}$, defined by the following
formulas for its covariant derivatives: for any local vector fields 
$X\in \Cal{T}_M,\,Y\in\Cal{T}_M^f,$
$$
\widehat{\nabla}_X(\widehat{Y}):=\lambda\widehat{X\circ Y},
\eqno{(1.19)}
$$
$$
 \widehat{\nabla}_{\Cal{E}}(\widehat{Y}):=\widehat{[E,Y]}.
\eqno{(1.20)}
$$
\endproclaim

\smallskip

\proclaim{\quad 1.5.2. Theorem} The extended structure connection is flat
iff $M$ is Frobenius and $E$ is Euler with $\roman{Lie}_E\,(\circ )=
d_0\circ $.
\endproclaim

\medskip

From (1.19) and (1.20) one can derive a formula for the covariant
derivative in the $\lambda$--direction: if $Y$ is flat, we have
$$
\widehat{[E,Y]}=\widehat{\nabla}_{E-d_0\lambda\partial /\partial\lambda}
(\widehat{Y})=\widehat{\nabla}_{E}(\widehat{Y})
-d_0\lambda\widehat{\nabla}_{\partial /\partial\lambda}
(\widehat{Y})=\lambda\widehat{E\circ Y}-d_0\lambda
\widehat{\nabla}_{\partial /\partial\lambda}
(\widehat{Y})
$$ 
so that
$$
d_0\widehat{\nabla}_{\partial /\partial\lambda}
(\widehat{Y})=\widehat{E\circ Y}-\frac{1}{\lambda}
\widehat{[E,Y]}.
\eqno{(1.21)}
$$
\medskip

{\bf 1.6. Semisimple Frobenius manifolds.}
Let $(M,g,A)$ be an associative pre--Frobenius manifold
of dimension $n.$

\medskip

\proclaim{\quad 1.6.1. Definition} $M$ is called semisimple
(resp. split semisimple) if an isomorphism of the sheaves
of $\Cal{O}_M$--algebras
$$
(\Cal{T}_M,\circ )\,\widetilde{\to}\,(\Cal{O}^n_M,
\roman{componentwise\ multiplication})
\eqno{(1.22)}
$$
exists everywhere locally (resp. globally.)
\endproclaim

\smallskip

This means that in a local (resp. global) basis 
$(e_1,\dots ,e_n)$ of $\Cal{T}_M$ the multiplication takes form
$$
(\sum f_ie_i)\circ (\sum g_je_j)=\sum f_ig_ie_i,
$$
and in particular,
$$
e_i\circ e_j=\delta_{ij}e_j.
\eqno{(1.23)}
$$
Such a family of idempotents is well defined up to renumbering.
Another way of saying this is that a semisimple manifold
comes with the structure group of $\Cal{T}_M$ reduced to
$S_n.$ Notice that $e_i$ are generally not flat, so that this
reduction is not compatible with that induced by
$\Cal{T}^f_M,$ with the structure group $GL(n).$

\smallskip

Hence if $M$ is semisimple, there exists an unramified covering
of degree $\le n!$, upon which the induced pre--Frobenius structure
is split.

\smallskip

Denote by $(\nu^i)$ the basis of 1--forms dual to $(e_i).$
>From (1.2) and (1.23) we find
$$
g(e_i,e_k)=g(e_i\circ e_i,e_k)=g(e_i,e_i\circ e_k)=\delta_{ik}g_{ii}.
$$
We will denote $g_{ii}$ by $\eta_i$. We see that 
the symmetric 2-form representing $g$ is diagonal
in the basis $(\nu^i)$:
$$
g=\sum_i\eta_i(\nu^i)^2.
\eqno{(1.24)}
$$
Moreover, according to (1.2), $A(e_i,e_j,e_k)=\delta_{ij}\delta_{ik}\eta_i$,
so that the symmetric 3-form representing $A$, is diagonal
with the same coefficients:
$$
A=\sum_i\eta_i(\nu^i)^3.
\eqno{(1.25)}
$$
Finally, $e:=\sum_ie_i$ is the identity in $(\Cal{T}_M,\circ ),$
and the co--identity, defined in 2.1.2, nicely complements
(1.24) and (1.25):
$$
\varepsilon =\sum_i\eta_i\nu^i.
\eqno{(1.26)}
$$
Thus the Definition 1.61 can be restated as follows:
\smallskip

\proclaim{\quad 1.6.2. Definition} The structure of the semisimple
pre--Frobenius manifold on $M$ is determined by the following data:

\smallskip

a). A reduction of the structure group of $\Cal{T}_M$ to $S_n$,
specified by a choice of local bases $(e_i)$ and dual bases
$(\nu^i).$

\smallskip

b). A flat metric $g$, diagonal in $(e_i),\,(\nu^i).$

\smallskip

c). A diagonal cubic tensor $A$ with the same coefficients as $g.$
\endproclaim

Associativity of $(\Cal{T}_M,\circ )$ is automatic in
both descriptions. However, potentiality (and the flatness of $g$
which we postulated) are non--trivial conditions.

\smallskip

\proclaim{\quad 1.7. Theorem} The structure described in the
Definition 3.2 is Frobenius iff the following
conditions are satisfied:

\smallskip

a). $[e_i,e_j]=0,$ or equivalently, $e_i=\partial /\partial u^i,\,
\nu^i=du^i$ for a local coordinate system $(u^i)$ called canonical
one.

\smallskip

b). $\eta_i=e_i\eta$ for a local function $\eta$ defined up to
addition of a constant. Equivalently, $\varepsilon$ is closed.
\endproclaim

\smallskip

We will call $\eta$ {\it the metric potential} of this structure.
(Sometimes this term refers to $h$ such that $g_{ab}=\partial_a\partial_b h;$
our meaning is different.)

\smallskip

Canonical coordinates are defined up to renumbering and constant
shifts.

\medskip 

{\bf 1.8. The Darboux--Egoroff equations.} The Theorem 1.7
establishes a (not very explicit) equivalence between
the following functional spaces on $M$ (modulo self--evident
equivalence):

\smallskip

a). Flat coordinates $(x^1,\dots ,x^n),$ flat metric $g_{ab}$,
function $\Phi (x)$ satisfying the Associativity Equations 
(1.6) and semisimplicity.

\smallskip

b). Canonical coordinates $(u^1,\dots ,u^n),$ function
$\eta (u)$ such that the metric $g=\sum e_i\eta (du^i)^2$
is flat, where $e_i=\partial /\partial u^i.$

\smallskip

The constraints on $\eta$, implicit in b), are called the
Darboux--Egoroff equations. In order to write them down explicitly,
let us introduce the rotation coefficients of the potential
metric:
$$
\gamma_{ij}:=\frac{1}{2}\frac{\eta_{ij}}{\sqrt{\eta_i\eta_j}}
\eqno{(1.27)}
$$
where as before, $\eta_i=e_i\eta ,$ $\eta_{ij}=e_ie_j\eta .$

\smallskip

\proclaim{\quad 1.8.1. Proposition} The diagonal potential metric
$g=\sum e_i\eta (du^i)^2$ is flat iff 
$\forall k\ne i\ne j\ne k:$
$$
\quad e_k\gamma_{ij}=\gamma_{ik}\gamma_{kj}
\eqno{(1.28)}
$$
and
$$
e\gamma_{ij}=0.
\eqno{(1.29)}
$$
\endproclaim

\medskip
\proclaim{\quad 1.8.2. Proposition} Let $e$ be the identity,
and $\varepsilon$ the co--identity of the semisimple
Frobenius manifold. Then

\smallskip

a). $\varepsilon =d\eta ,$ where $\eta$ is the metric potential.

\smallskip

b). $e$ is flat iff for all $i$, $e\eta_i=0$, or equivalently,
$e\eta = g(e,e) =\roman{const}.$ This condition is satisfied
in the presence of an Euler field with $D\ne 2d_0$
(see (1.12), (1.13), (1.14).)

\smallskip

c). If $e$ is flat, and $(x^a)$ is a flat coordinate system, then
$$
\eta = \sum_a x^ag(\partial_a ,e)+\roman{const} .
\eqno{(1.30)}
$$
\endproclaim
\smallskip

The formula (1.30) shows that in the passage from 
the $(x^a,\Phi )$--description to the $(u^i,\eta )$--description
the main information is encoded in the transition formulas
$u^i=u^i(x)$, at least in the presence of flat identity.

\smallskip

Like the identity, the Euler field
is almost uniquely defined by the canonical coordinates,
if it exists at all.

\smallskip

\proclaim{\quad 1.9. Theorem} Let $E$ be a vector field
on the semisimple Frobenius manifold $M$, $d_0$ a constant.

\smallskip

a). We have
$\roman{Lie}\,(\circ )=d_0(\circ )$, iff 
$$
E=d_0\sum_i (u^i+c^i)e_i,
\eqno{(1.31)}
$$
where $c^i$ are some constants.

\smallskip

b). For the field of the form (3.17) and a constant $D$, we have $\roman{Lie}_E(g)=Dg$ iff for all $i$, $E\eta_i=(D-2d_0)\eta_i$,
or equivalently
$$
E\eta=(D-d_0)\eta + \roman{const}.
\eqno{(1.32)}
$$
\endproclaim

Thus in the presence of a non--vanishing
Euler field we may and will normalize the canonical coordinates
so that $E=d_0\sum u^ie_i.$

\medskip

{\bf 1.10. A pencil of flat metrics.} Equations (1.28) are stable with
respect to a semigroup of coordinate changes. Namely, let $f_i$ be
arbitrary
functions of one variable such that $\check {u}^i:=f_i(u^i)$
form a local coordinate system, $\check {e}_i=\partial /\partial
\check {u}^i, \check {\eta}_i=\check {e}_i\eta$ etc.
\smallskip
\proclaim{\quad 1.10.1. Proposition} If $(e_i,\gamma_{ij})$
satisfy (1.28), then $(\check {e}_i,\check {\gamma}_{ij})$
satisfy (1.28) as well.
\endproclaim

\medskip

In order to satisfy (1.29) as well, we will have to restrict ourselves
to the one--parameter family of local coordinate changes
$$
\check {u}^i=\roman{log}\,(u^i-\lambda ),\
\check {e}_i= (u^i-\lambda)e_i,\
\check {g}_{\lambda}=\sum_i(u^i-\lambda )^{-1}e_i\eta (du^i)^2
\eqno{(1.33)}
$$
which make sense on $M_{\lambda}:=\{x\in M\,|\,\forall i,\,
u^i\ne\lambda\}.$

\smallskip

\proclaim{\quad 1.10.2. Theorem} Let $M$ be a semisimple Frobenius manifold
with canonical coordinates $(u^i)$ and metric potential $\eta .$ 
Then the following statements are equivalent.

\smallskip

a). For all $\lambda ,$ the structure (1.33) is semisimple Frobenius on 
$M_{\lambda}.$

\smallskip

b). The same for a particular value of $\lambda .$

\smallskip

c). For all $i\ne j,$
$$
\sum_ku^ke_k\gamma_{ij}=-\gamma_{ij}.
\eqno{(1.34)}
$$
Moreover, (1.34) is satisfied if $E=\sum_ku^ke_k$ is the Euler field
on $M$ with $d_0=1.$
\endproclaim
\smallskip

Notice that generally $\check {e}=\sum\check {e}_k$ is not
flat for $\check {g}_{\lambda}$ and $\check {E}=\sum
\check {u}^k\check {e}_k$ is not an Euler field.

\smallskip

If $E$ is Euler, the metric $\check {g}_{\lambda}$
in (1.33) can be written in coordinate free form:
$$
\check {g}_{\lambda}(X,Y)=g((E-\lambda )^{-1}\circ X,Y).
\eqno{(1.35)}
$$
In fact (1.35) is flat on any Frobenius manifold with
semisimple Euler field on it: cf. [D].

\medskip

{\bf 1.11. The second structure connection.} From now on, 
we will restrict ourselves to the
case of semisimple complex Frobenius manifolds carrying an Euler
field with $d_0=1$ and admitting a global system of canonical
coordinates $(u^i).$ We will call {\it the second structure
connection} $\check \nabla_{\lambda}$ the Levi--Civita
connection of the flat metric (1.35),
depending on a parameter $\lambda$ and defined on the open
subset $M_{\lambda}\subset M$ where $u^i\ne\lambda$
for all $i.$ Put $\check M:=\cup_{\lambda}(M_{\lambda}\times\{\lambda\})
\subset M\times\bold{P}^1_{\lambda}$ and denote by
$\check \Cal{T}$ the restriction of $\roman{pr}^*_M(\Cal{T}_M)$
to $\check M.$

\smallskip

We will construct a flat extension of
$\check \nabla$ of $\check \nabla_{\lambda}$ to $\check \Cal{T}$
which will also be referred to as the second structure connection.
Both extensions $\widehat{\nabla}$ and $\check \nabla$ will be further
studied as isomonodromic deformations of their restrictions
to the $\lambda$--direction parametrized by $M.$

\smallskip
 
More precisely, assume that $\Cal{T}^f_M$ is a trivial local system
(for instance, because $M$ is simply connected.) Put $T:=\Gamma (M,\Cal{T}^f_M).$
Then $\widehat{\nabla}$ (resp. $\check \nabla$) induces an integrable family
of connections with singularities on the trivial bundle on $\bold{P}^1_{\lambda}$ with the fiber $T.$ The first connection 
$\widehat{\nabla}$ is singular only at $\lambda =0$
and $\lambda =\infty$ but whereas $0$ is a regular (Fuchsian)
singularity, $\infty$ is irregular one, so that
$\widehat{\nabla}$ cannot be an algebraic geometric Gauss--Manin 
connection, and its monodromy involves the Stokes phenomenon.
To the contrary, the second connection $\check \nabla$
generally has only regular singularities at infinity and at
 $\lambda =u^i$
whose positions thus depend on the parameters. It is determined 
by the conventional monodromy representation and has a chance
to define a variation of Hodge structure. For more details, see the
next section.

\smallskip

It turns out that both deformations have a common moduli space
and deserve to be studied together. In fact, fiberwise they are more or less
formal Laplace transforms of each other. More to the point,
they form a dual pair in the sense of J.~Harnad.

\smallskip

In our calculations the key role will be played by the $\Cal{O}_M$--linear
skew symmetric operator $\Cal{V}:\,\Cal{T}_M\to\Cal{T}_M$
which is the unique extension of the operator defined in 1.4
on flat vector fields by the formula
$$
\Cal{V}(X)=[X,E]-\frac{D}{2}X\ \roman{for}\ X\in\Cal{T}^f_M.
\eqno{(1.36)}
$$
\smallskip

\proclaim{\quad 1.11.1. Proposition} a). We have for arbitrary 
$X\in\Cal{T}_M:$
$$
\Cal{V}(X)=\nabla_{0,X}(E)-\frac{D}{2}\,X.
\eqno{(1.37)}
$$
\smallskip

b). Let $e_j=\partial /\partial u^j,\,f_j=e_j/\sqrt{\eta_j}.$
Then
$$
\Cal{V}(f_i)=\sum_{j\ne i}(u^j-u^i)\gamma_{ij}f_j.
\eqno{(1.38)}
$$
\endproclaim

Formula (1.37) defines an $\Cal{O}_M$--linear endomorphism of
$\Cal{T}_M$ which coincides with (1.36) on the flat fields,
as a calculation in flat coordinates shows.
To check (1.38), we can use (1.37) and the classical
explicit expressions for the Levi--Civita connection,
cf. below.

\smallskip

We can now state the main result of this section. In addition
to (1.36), define the operator $\Cal{U}:\,\Cal{T}_M\to\Cal{T}_M:$
$$
\Cal{U}(X):=E\circ X,
\eqno{(1.39)}
$$
so that $\Cal{U}(f_i)=u^if_i.$
\smallskip

\proclaim{\quad 1.12. Theorem} For $X,Y\in\roman{pr}^{-1}_M(\Cal{T}_M)\subset
\Cal{T}_{\check M}$ (meromorphic vector fields on 
$\Cal{T}_{M\times \bold{P}^1_{\lambda}}$ independent on $\lambda$) put
$$
\check \nabla_X(Y)=\nabla_{0,X}(Y)-(\Cal{V}+\frac{1}{2}\,\roman{Id})\,(\Cal{U}
-\lambda )^{-1}(X\circ Y),
\eqno{(1.40)}
$$
$$
\check \nabla_{\partial /\partial \lambda}(Y)=(\Cal{V}+\frac{1}{2}
\,\roman{Id})\,(\Cal{U}
-\lambda )^{-1}(Y).
\eqno{(1.41)}
$$
Then $\check \nabla$ is a flat connection on $\check \Cal{T}$ whose
restriction on $M\times\{\lambda\}$ defined by (1.40) is the Levi--Civita
connection for $\check g_{\lambda}.$
\endproclaim

\smallskip

{\bf Remark.} Rewriting $\widehat{\nabla}$ in the same notation, we get
$$
\widehat{\nabla}_X(Y)=\nabla_{0,X}(Y)+
\lambda X\circ Y,
\eqno{(1.42)}
$$
$$
\widehat{\nabla}_{\partial /\partial \lambda}(Y)=
\left[\Cal{U}+\frac{1}{\lambda}(\Cal{V}+\frac{D}{2}\,\roman{Id})\right](Y).
\eqno{(1.43)}
$$

\smallskip

{\bf Proof.} We will first  calculate
the Levi--Civita connection for $\check g_{\lambda}$ in coordinates
$\check u^i=\roman{log}\,(u^i-\lambda ):$
$$
\check e_i=\frac{\partial}{\partial \check u^i}=(u^i-\lambda )e_i,\
\check \eta_i=(u^i-\lambda )\eta_i,\
\check \eta_{ij}=(u^i-\lambda )(u^j-\lambda )\eta_{ij}+
\delta_{ij}(u^i-\lambda )\eta_i,
$$
$$
\check \gamma_{ij}=\gamma_{ij}(u^i-\lambda )^{1/2}(u^j-\lambda )^{1/2}.
$$
Then for $i\ne j$
$$
\check \nabla_{\check e_i}(\check e_j)=
\frac{1}{2}\frac{\check \eta_{ij}}{\check \eta_{i}}\check e_i+
\frac{1}{2}\frac{\check \eta_{ij}}{\check \eta_{j}}\check e_j=
\frac{1}{2}(u^i-\lambda )(u^j-\lambda )
\left( \frac{\eta_{ij}}{\eta_{i}}e_i+
\frac{\eta_{ij}}{\eta_{j}}e_j\right)
$$
so that
$$
\check \nabla_{e_i}(e_j)=
\frac{1}{2}\frac{\eta_{ij}}{\eta_{i}}e_i+
\frac{1}{2}\frac{\eta_{ij}}{\eta_{j}}e_j
=\nabla_{0,e_i}(e_j).
\eqno{(1.44)}
$$
Similarly,
$$
\check \nabla_{\check e_i}(\check e_i)=
\frac{1}{2}\frac{\check \eta_{ii}}{\check \eta_{i}}\check e_i-
\frac{1}{2}\sum_{j\ne i}\frac{\check \eta_{ij}}{\check \eta_{j}}\check e_j=
$$
$$
\frac{1}{2}(u^i-\lambda )^2\left[\frac{\eta_{ii}}{\eta_i}-
\frac{1}{u^i-\lambda}\right]e_i-
\frac{1}{2}\sum_{j\ne i}({u^i-\lambda})\,({u^j-\lambda})\,
\frac{\eta_{ij}}{\eta_j}\,e_j
$$
so that
$$
\check \nabla_{e_i}(e_i)=
\frac{1}{2}\left[\frac{\eta_{ii}}{\eta_i}-
\frac{1}{u^i-\lambda}\right]e_i-
\frac{1}{2}\sum_{j\ne i}\frac{u^j-\lambda}{u^i-\lambda}\,
\frac{\eta_{ij}}{\eta_j}\,e_j.
\eqno{(1.45)}
$$
Subtracting from this the Levi--Civita covariant derivative, we get
$$
(\check \nabla_{e_i}-\nabla_{0,e_i})(e_i)=
-\frac{1}{2}\,\frac{1}{u^i-\lambda}\,e_i-
\frac{1}{2}\sum_{j:\,j\ne i}\frac{u^j-u^i}{u^i-\lambda}\,
\frac{\eta_{ij}}{\eta_j}\,e_j
\eqno{(1.46)}
$$
and
$$
(\check \nabla_{e_i}-\nabla_{0,e_i})(f_i)=
-\frac{1}{2}\,\frac{1}{u^i-\lambda}\,f_i-
\sum_{j\ne i}\frac{u^j-u^i}{u^i-\lambda}\,
\gamma_{ij}\,f_j.
\eqno{(1.47)}
$$
In view of (1.38), we can write (1.44) and (1.45) together as
$$
(\check \nabla_{e_i}-\nabla_{0,e_i})(f_j)=
-\left(\Cal{V}+\frac{1}{2}\,\roman{Id}\right)(\Cal{U}-\lambda )^{-1}
(e_i\circ f_j )
\eqno{(1.48)}
$$
because $e_i\circ f_j=\delta_{ij}f_j.$ This family of formulas
is equivalent to (1.40) so that (1.40) is the Levi--Civita
connection for $\check g_{\lambda}.$ In particular, it is flat
for each fixed $\lambda .$

\smallskip

Since $[X,{\partial}/{\partial \lambda}]=0$ for
$X\in\roman{pr}^{-1}_M(\Cal{T}_M),$ it remains to show that
the covariant derivatives (1.40) and (1.41) commute on $\check M$
i.~e.~, that for all $i,j$
$$
\check \nabla_{e_i}\check \nabla_{\partial /\partial \lambda}(e_j)=
\check \nabla_{\partial /\partial \lambda}\check \nabla_{e_i}(e_j).
\eqno{(1.49)}
$$
First of all, from (1.41) and (1.42) we find
$$
\check \nabla_{\partial /\partial\lambda}(e_j)=
\frac{1}{2}\,\frac{1}{u^j-\lambda}\,e_j
+\frac{1}{2}\sum_{k\ne j}\frac{u^k-u^j}{u^j-\lambda}\,
\frac{\eta_{jk}}{\eta_k}\,e_k.
\eqno{(1.50)}
$$
Together with (1.44) and (1.45) this gives for $i\ne j$:
$$
\check \nabla_{\partial /\partial \lambda}\check \nabla_{e_i}(e_j)=
\frac{1}{2}\,\frac{\eta_{ij}}{\eta_j}
\left[\frac{1}{2}\,\frac{1}{u^j-\lambda}\,e_i
+\frac{1}{2}\sum_{k\ne i}\frac{u^k-u^i}{u^i-\lambda}\,
\frac{\eta_{ik}}{\eta_k}\,e_k\right]+(i\,\leftrightarrow\,j),
\eqno{(1.51)}
$$
$$
\check \nabla_{e_i}\check \nabla_{\partial /\partial \lambda}(e_j)=
\frac{1}{2}\,\frac{1}{u^j-\lambda}
\left( \frac{1}{2}\frac{\eta_{ij}}{\eta_{i}}e_i+
\frac{1}{2}\frac{\eta_{ij}}{\eta_{j}}e_j\right)+
$$
$$
+\frac{1}{2}\sum_{k\ne j}e_i\left(\frac{u^k-u^j}{u^j-\lambda}\,
\frac{\eta_{jk}}{\eta_k}\right)e_k+
\frac{1}{2}\sum_{k\ne j,i}\frac{u^k-u^j}{u^j-\lambda}\,
\frac{\eta_{jk}}{\eta_k}\,
\left( \frac{1}{2}\,\frac{\eta_{ik}}{\eta_{i}}\,e_i+
\frac{1}{2}\,\frac{\eta_{ik}}{\eta_{k}}\,e_k\right)+
$$
$$
\frac{1}{2}\,\frac{u^i-u^j}{u^j-\lambda}\,
\frac{\eta_{ik}}{\eta_k}
\left[ \frac{1}{2}\,\left(\frac{\eta_{ii}}{\eta_i}
-\frac{1}{u^i-\lambda}\right)\,e_i
-\frac{1}{2}\sum_{j\ne i}\frac{u^j-\lambda}{u^i-\lambda}\,
\frac{\eta_{ij}}{\eta_j}\,e_j
\right].
\eqno{(1.52)}
$$
The coincidence of coefficients of $e_k$ in (1.51) and (1.52)
for $i\ne j\ne k\ne i$ can be checked with the help of the
following identities which are equivalent to the
Darboux--Egoroff equations:
$$
\eta_{ijk}=\frac{1}{2}\left( \frac{\eta_{ik}\eta_{jk}}{\eta_k}+
\frac{\eta_{ij}\eta_{ik}}{\eta_i}+\frac{\eta_{ij}\eta_{jk}}{\eta_j}
\right).
$$
The coincidence of the coefficients of $e_i$ requires
a little more work, and we will give some details,
again for the case $i\ne j.$

\smallskip

In (1.51) the coefficient of $e_i$ is
$$
\frac{1}{4}\,\frac{1}{u^i-\lambda}\,\frac{\eta_{ij}}{\eta_i}+
\frac{1}{4}\,\frac{u^k-u^j}{u^j-\lambda}\,\frac{\eta_{ij}^2}{\eta_i\eta_j},
\eqno{(1.53)}
$$
whereas in (1.52) we get
$$
\frac{1}{4}\,\frac{1}{u^j-\lambda}\,\frac{\eta_{ij}}{\eta_i}+
\frac{1}{2}\,e_i\,\left(\frac{u^i-u^j}{u^j-\lambda}\,
\frac{\eta_{ij}}{\eta_i}\right)+
$$
$$
+\frac{1}{4}\sum_{k\ne i,j}\frac{u^k-u^j}{u^j-\lambda}\,
\frac{\eta_{ik}\eta_{jk}}{\eta_i\eta_k}+
\frac{1}{2}\left(\frac{\eta_{ii}}{\eta_i}-\frac{1}{u^i-\lambda}\right).
\eqno{(1.54)}
$$
To identify (1.53) and (1.54) we have to get rid of the
sum $\sum_k$ in (1.54). This can be done with the help
of (1.28), (1.29)  and (1.34):
$$
\frac{1}{4}\sum_{k\ne i,j}\frac{u^k-u^j}{u^j-\lambda}\,
\frac{\eta_{ik}\eta_{jk}}{\eta_i\eta_k}=
\frac{1}{u^j-\lambda}\,\frac{\eta_j^{1/2}}{\eta_i^{1/2}}\,
\left[\sum_{k\ne i,j}u^k\gamma_{ik}\gamma_{kj}-
u^j\sum_{k\ne i,j}\gamma_{ik}\gamma_{kj}\right]=
$$
$$
=\frac{1}{u^j-\lambda}\,\frac{\eta_j^{1/2}}{\eta_i^{1/2}}\,
\left[-\gamma_{ij}-u^ie_i\gamma_{ij}-u^je_j\gamma_{ij}
+u^j(e_i+e_j)\gamma_{ij}\right]=
$$
$$
=\frac{1}{2}\,\frac{1}{u^j-\lambda}\,\left[
-\frac{\eta_{ij}}{\eta_i}+(u^j-u^i)\left(
\frac{\eta_{iij}}{\eta_{i}}-\frac{1}{2}\,
\frac{\eta_{ij}\eta_{ii}}{\eta_{i}^2}
-\frac{1}{2}\,
\frac{\eta_{ij}^2}{\eta_{i}\eta_j}\right)\right].
$$
The remaining part of the calculation is straightforward, and we leave it
to the reader, as well as the case $i=j$ which is treated
similarly.

\medskip

{\bf 1.13. Formal Laplace transform.} Assume now that $\Cal{T}^f_M$
is a trivial local system. This means that if we put
$T:=\Gamma (M,\Cal{T}^f_M),$ there is a natural isomorphism
$\Cal{O}_M\otimes T\to\Cal{T}_M.$

\smallskip

Formulas (1.41) (resp. (1.43)) define two families of
connections with singularities on the trivial vector
bundle on $\bold{P}^1_{\lambda}$ with fiber $T$, parametrized
by $M.$ Namely, denote by $\partial_{\lambda}$ the
covariant derivative along
$\partial /\partial \lambda$ on this
bundle for which the constant sections are horizontal.
Then the two connections are
$$
\check \nabla_{\partial /\partial \lambda} =\partial_{\lambda}+\left(\Cal{V}+\frac{1}{2}\,\roman{Id}\right)
\,(\Cal{U}-\lambda )^{-1},
\eqno{(1.55)}
$$
$$
\widehat{\nabla}_{\partial /\partial \lambda} =\partial_{\lambda}+\Cal{U}+
\frac{1}{\lambda}\left(\Cal{V}+\frac{D}{2}\,\roman{Id}\right).
\eqno{(1.56)}
$$
\smallskip

Let $M,N$ be two $\bold{C}[\lambda ,\partial_{\lambda}]$--modules.
{\it A formal Laplace transform} $M\to N:\ Y\mapsto Y^t$
is a $\bold{C}$--linear map for which
$$
(-\lambda Y)^t=\partial_{\lambda}(Y^t),\ (\partial_{\lambda}Y)^t=\lambda Y^t.
\eqno{(1.57)}
$$
The archetypal Laplace transform is the Laplace integral
$$
Y^t(\mu )=\int e^{-\lambda\mu}Y(\lambda )d\lambda
\eqno{(1.58)}
$$
taken along a contour (not necessarily closed) in $\bold{P}^1(\bold{C}).$
In an analytical setting we have to secure the convergence
of (1.58), the possibility to derivate under the integral sign
and the identity 
$$
\int \partial_{\lambda}(e^{-\lambda\mu}Y(\lambda ))d\lambda =0.
$$
However, (1.58) may admit other interpretations,
for instance, in terms of asymptotic series.

\smallskip

Let now $M$ (resp. $N$) be two  
$\bold{C}[\lambda ,\partial_{\lambda}]$--modules of local
(or formal, or distribution) sections of $\bold{P}^1_{\lambda}\times T$
so that the operators $\check \nabla \cdot (\Cal{U}-\lambda )$
(resp. $\lambda\widehat{\nabla}$) make sense in $M$ (resp. $N$)
(cf. (1.55), resp. (1.56)), and assume that we are given a formal Laplace
transform $M\to N.$

\smallskip

\proclaim{\quad 1.13.1. Proposition} We have:
$$
[\check \nabla_{\partial /\partial \lambda} ((\Cal{U}-\lambda )Y)]^t=
(\lambda\widehat{\nabla}_{\partial /\partial \lambda}+\frac{1-D}{2})\,Y^t=
\lambda^{\frac{D+1}{2}}\widehat{\nabla}_{\partial /\partial
\lambda}(\lambda^{\frac{1-D}{2}}Y^t). 
$$
In particular, $\lambda^{\frac{1-D}{2}}Y^t$ is $\widehat{\nabla}$--horizontal,
if $(\Cal{U}-\lambda )Y$ is $\check \nabla$--horizontal.
\endproclaim

\smallskip

{\bf Proof.} Using (1.55)--(1.57), we find:
$$
[\check \nabla_{\partial /\partial \lambda}((\Cal{U}-\lambda )Y)]^t=
\left[(\partial_{\lambda}\cdot (\Cal{U}-\lambda )+
\Cal{V}+\frac{1}{2}\,\roman{Id})Y\right]^t=
$$
$$
=\left[\lambda \,(\Cal{U}+\partial_{\lambda})+
\Cal{V}+\frac{1}{2}\,\roman{Id}\right]Y^t=
$$
$$
=\left[\lambda\widehat{\nabla}_{\partial /\partial \lambda}+
\frac{1-D}{2}\,\roman{Id}\right]Y^t=
\lambda^{\frac{D+1}{2}}\widehat{\nabla}_{\partial /\partial \lambda}(
\lambda^{\frac{1-D}{2}}Y^t).
$$
For a more detailed discussion of the formal Laplace transform,
see  [S], 1.6.

\smallskip

For later use we note that the connection $\check{\nabla}$ defined by
(1.40), (1.41) can be further deformed. Namely, for any
constant $s$ put
$$
\check{\nabla}_X^{(s)}(Y)=\check{\nabla}_X(Y)-
s(\Cal{U}-\lambda )^{-1}(X\circ Y),
\eqno(1.59)
$$
$$
\check{\nabla}_{\partial /\partial \lambda}^{(s)}(Y)=
\check{\nabla}_{\partial /\partial \lambda}(Y)+
s(\Cal{U}-\lambda )^{-1}(Y).
\eqno(1.60)
$$
\smallskip

\proclaim{\quad 1.14. Theorem} $\check{\nabla}_X^{(s)}$
is a flat connection on $\check{\Cal{T}}.$
\endproclaim

\smallskip

This can be checked by a direct calculation similar to that
in the proof of the Theorem 1.12. The formal Laplace transform
of $\check{\nabla}_{\partial /\partial \lambda}^{(s)} $ is given by
$$
[\check \nabla^{(s)}_{\partial /\partial \lambda} ((\Cal{U}-\lambda )Y)]^t=
(\lambda\widehat{\nabla}_{\partial /\partial \lambda}+\frac{1-D+2s}{2})\,Y^t=
\lambda^{\frac{D+1-2s}{2}}\widehat{\nabla}_{\partial /\partial
\lambda}(\lambda^{\frac{1-D+2s}{2}}Y^t). 
$$
In particular, $\lambda^{\frac{1-D+2s}{2}}Y^t$ is $\widehat{\nabla}$--horizontal,
if $(\Cal{U}-\lambda )Y$ is $\check \nabla^{(s)}$--horizontal.

\newpage

\centerline{\bf \S 2. Schlesinger equations}

\bigskip

{\bf 2.1. Singularities of meromorphic connections.} Let $N$
be a complex manifold, $D\subset N$ a closed complex submanifold
of codimension one, $\Cal{F}$ a locally free sheaf of finite rank
on $N.$ A meromorphic connection with singularities on $D$
is given by a covariant differential
$\nabla :\Cal{F}\to\Cal{F}\otimes\Omega_N^1((r+1)D)$
for some $r\ge 0.$ It is called flat (or integrable) if it is
flat outside $D.$ We start with a list of elementary notions and
constructions that will be needed later. They depend only on the local
behaviour of $\Cal{F}$ and $\nabla$ in a neighborhood
of $D$, so we will assume $D$ irreducible.

\medskip

{\it i) Order of singularity.} We will say that $\nabla$ as above
is of order $\le r+1$ on $D$ if $\nabla_X(\Cal{F})\subset
\Cal{F}(rD)$ for any vector field $X$ tangent to $D$ (i.~e. satisfying
$XJ_D\subset J_D$ where $J_D$ is the ideal of $D$), and
$\nabla_X(\Cal{F})\subset\Cal{F}((r+1)D)$ in general.
Locally, if $(t^0,t^1,\dots ,t^n)$ is a coordinate system
on $N$ such that $t^0=0$ is the equation of $D,$
the connection matrix of $\nabla$ in a basis of $\Cal{F}$
can be written as
$$
G_0\frac{dt^0}{(t^0)^{r+1}}+
\sum_{i=1}^nG_i\frac{dt^i}{(t^0)^{r}}
\eqno{(2.1)}
$$
where $G_i=G_i(t^0,t^1,\dots ,t^n)$ are holomorphic matrix functions.

\smallskip

Note that $G_0(0,t^1, \dots, t^n)\in H^0(D, \text{End} {\Cal F})$ 
is well-defined, i.e.\ it  does not depend on the choice of local
coordinates. It is called the residue of $\nabla$ at $D$ and is denoted
by $\text{res}_D(\nabla)$.

\medskip

{\it ii) Restriction to a transversal submanifold.}
Let $i:\,N^{\prime}\to N$ be a closed embedding of a submanifold
transversal to $D,\, D^{\prime}=N^{\prime}\cap D,\,
\Cal{F}^{\prime}=i^*(\Cal{F}).$ Then the induced connection
$\nabla^{\prime}=i^*(\nabla )$ on $\Cal{F}^{\prime}$ is flat
and of order $\le r+1$ on $D^{\prime}$ if $\nabla$ has these properties.

\medskip

{\it iii) Residual connection.} 
 Assume that $\nabla$ is of order $\le 1$
on $D.$ For any given local trivialization $f$ of $[D]$, one can define a
connection without
singularities $\nabla^{D,f}$ on $j^*(\Cal{F})$ where $j$ is the embedding of
$D$ in $N.$ Namely, to define $\nabla_{X^{\prime}}^{D,f}(s^{\prime})$ where
$s^{\prime}\in j^*(\Cal{F}),\,X^{\prime}\in \Cal{T}_D$, we extend
locally $s^{\prime}$ to a section $s$ of $\Cal{F}$, $X^{\prime}$
to a vector field $X$ on $N$, $\text{res}_D(\nabla)$ to a section 
$\text{res}(\nabla)$ of ${\Cal F}\otimes {\Cal F}^*$ on $N$,  
 calculate $(\nabla_X - \frac{Xf}{f}
\text{res}(\nabla))(s)$  and restrict it
to $D.$ One checks that the result does not depend on the choices made.
In the notation of (2.1), the matrix of the residual connection
can be written as ($r=0$):
$$
\sum_{i=1}^nG_i(0,t^1,\dots ,t^n)dt^i.
\eqno{(2.2)}
$$
If $\nabla$ is flat, 
 $\nabla^{D,f}$ is flat for any local trivialization $f$ of $[D]$.


\medskip

{\it iv) Principal part of order $r+1.$} Similarly to (2.2),
we can consider the matrix function on $D$
$$
G_0(0,t^1,\dots ,t^n)
\eqno{(2.3)}
$$
which we will call the principal part of order $r+1$ of $\nabla .$
In more invariant terms, it is the $\Cal{O}_D$--linear
map $j^*(\Cal{F})\to j^*(\Cal{F})$ induced by
$\Cal{F}\to j^*(\Cal{F}):\, s\mapsto
(t^0)^{r+1}\nabla_{\partial /\partial t^0}(s)\,|_D$.
For $r\ge 1$ it depends on the choice
of local coordinates, and is multiplied by an invertible local function
on $D$ when this choice is changed. Hence its spectrum 
is well defined globally on $D$ for $r=0$, and the simplicity of the spectrum
makes sense for any $r$.

\medskip

{\it v) Tameness and resonance.} Two general position conditions are
important in the study of meromorphic singularities of order $\le r+1.$

\smallskip

If $r\ge 1$ (irregular case), the singularity is called {\it tame},
if the spectrum of its principal part at any point of $D$
is simple.

\smallskip

If $r=0$ (regular case), the singularity is called {\it non--resonant},
if it is tame and moreover, the difference of any two eigenvalues
never takes an integer value on $D.$ 

\medskip

{\bf 2.1.1. Example: the structure connections of Frobenius
manifolds.} As above, we will assume
that $\Cal{T}^f_M$ is trivial, and its fibers are identified
with the space $T$ of global flat vector fields.

\smallskip

Put  $N=M\times\bold{P}^1_{\lambda},\,\Cal{F}=\Cal{O}_N\otimes T.$
We can apply the previous considerations to $\widehat{\nabla}$
and $\check \nabla.$

\smallskip

{\it Analysis of  $\widehat{\nabla}$.} Clearly, $\widehat{\nabla}$
has singularity of order 1 at $\lambda =0$
(i.~e. on $D_0=M\times\{0\}$) and of order 2 at
$\lambda =\infty$ (i.~e. on $D_{\infty}=M\times\{\infty\}$):
cf. (1.42) and (1.43). Restricting $\widehat{\nabla}$ to
$\{y\}\times\bold{P}^1_{\lambda}$ for various $y\in M$ we get a family
of meromorphic connections on $\bold{P}^1_{\lambda}$ parametrized
by $M.$ 

\smallskip

The residual connection is defined on $D_0=M$ and it coincides
with the Levi--Civita connection of $g.$ The principal part of order 1
on $D_0$ is $\Cal{V}+\dfrac{D}{2}\,\roman{Id}.$ The eigenvalues of this
operator do not depend on $y\in D_0$: in 1.5 they were denoted
$(d_a).$ In the case of quantum cohomology
the principal part
{\it is always resonant.}

\smallskip

The principal part of order 2 on $D_{\infty}=M$ is (proportional to)
$\Cal{U}$ (cf. (1.43), use the local equation $\mu =\lambda^{-1}=0$
for $D_{\infty }.$) Its eigenvalues now depend on $y\in M:$ they are just
the canonical coordinates $u^i(y).$ We will call the point $y$ 
{\it tame} if $u^i(y)\ne u^j(y)$ for $i\ne j.$ We will call
$M$ tame, if all its points are tame. Every $M$ contains 
the maximum tame subset which is open and dense.

\medskip

{\it Analysis of $\check \nabla .$} According to (1.40), (1.41),
$\check \nabla$ has singularities of order 1 at the divisors
$\lambda =u^i$ and $\lambda =\infty .$ These divisors
do not intersect pairwise iff $M$ is tame. The principal part
of order 1 at $\lambda =u^i$ is  
$-(\Cal{V}+\dfrac{1}{2}\,\roman{Id})\cdot (e_i\circ ).$

\smallskip

The residual connection of $\check \nabla$ on $\lambda =\infty$
is again the Levi--Civita connection $\nabla_0$ of $g.$
In fact, using (1.48) we find
$$
\check \nabla=d\lambda\check \nabla_{\partial /\partial\lambda}
+\sum_idu^i\check \nabla_{e_i}=
$$
$$
=d\lambda\check \nabla_{\partial /\partial\lambda}
+\sum_idu^i[\nabla_{0,e_i}-
(\Cal{V}+\frac{1}{2}\,\roman{Id})\,(\Cal{U}-\lambda )^{-1}(e_i\circ )].
$$
Replacing $\lambda$ by the local parameter $\mu=\lambda^{-1}$ at infinity,
we have
$$
\check \nabla =d\mu\check \nabla_{\partial /\partial\mu}
+\sum_idu^i[\nabla_{0,e_i}-\mu\,
(\Cal{V}+\frac{1}{2}\,\roman{Id})\,(\mu\,\Cal{U}-\roman{Id} )^{-1}(e_i\circ )]
$$
so that the expression (2.2) (with $(\mu,u^1,\dots ,u^m)$
in lieu of $(t^0,t^1,\dots ,t^n)$) becomes
$\sum_idu^i\nabla_{0,e_i}=\nabla_0.$

\medskip

{\bf 2.2. Versal deformation.} We will now review the basic results
on the deformation of meromorphic connections on
$\bold{P}^1_{\lambda}$, restricting ourselves to the case
of singularities of order $\le 2.$ This suffices for applications
to both structure connections, on the other hand, this is precisely
the case treated in full detail by B.~Malgrange in [Mal4],
Theorem 3.1. It says that the positions of finite poles and 
the spectra of the principal parts of order 2 form coordinates on the
coarse moduli space with tame singularities. To be more precise,
one has to rigidify the data slightly.

\smallskip

Let $\nabla^0$ be a meromorphic connection on a locally free sheaf
$\Cal{F}^0$ on $\bold{P}^1_{\lambda}$ of rank $p$, with
$m+1\ge 2$ tame singularities (including $\lambda =\infty$)
of order $\le 2.$ Call {\it the rigidity} for $\nabla^0$
the following data:

\smallskip

a). A numbering of singular points: $a_0^1,\dots ,a_0^m,\,a^{m+1}=\infty.$

\smallskip

b). The subset $I\subset \{1,\dots ,m+1\}$ such that $a_0^j$
is of order 2 exactly when $j\in I.$

\smallskip

c). For each $j\in I,$ a numbering $(b_0^{j1},\dots , b_0^{jp})$
of the eigenvalues of the principal part at $a_0^j.$

\smallskip

Construct the space $B=B(m,p,S)$ as the universal covering
of
$$
(\bold{C}^m\setminus\roman{diagonals})\times\prod_{j\in I}
(\bold{C}^p\setminus\roman{diagonals})
$$
with the base point $(a_0^i;\,b_0^{jk})$, let $b_0\in B$ be its
lift. We denote by $a^i,\,b^{jk}$ the coordinate functions lifted to
$B.$ Let $i:\, \bold{P}^1_{\lambda}\to B\times \bold{P}^1_{\lambda}$
be the embedding $\lambda\mapsto (b_0,\lambda ),$ and $D_j$ the  divisor
$\lambda =a^j$ in $B\times \bold{P}^1_{\lambda}.$

\smallskip

\proclaim{\quad 2.2.1. Theorem ([Mal4], Th. 3.1)} For a given
$(\nabla^0,\Cal{F}^0)$ with rigidity, there exists
a locally free sheaf $\Cal{F}$ of rank $p$ on $ \bold{P}^1_{\lambda}\times B,$
a flat meromorphic connection $\nabla$ on it, and an isomorphism
$i^0:\,i^*(\Cal{F},\nabla )\to (\Cal{F}^0,\nabla^0)$ with the following
properties:

\smallskip

$D_j,\,j=1,\dots ,m+1,$ are all the poles of $\nabla ,$ of order
1 (resp. 2) if $j\ne I$ (resp. $j\in I.$) If $j\in I,$ then
$(b^{j1},\dots ,b^{jp})$ (as functions on $D_j$) form the spectrum
of the principal part of order 2 of $\nabla$ at $D_j.$

\smallskip

It follows that the restrictions of $\nabla$ to the fibers
$\{b\}\times \bold{P}^1_{\lambda}$ are endowed with the induced
rigidity, and $i^0$ is compatible with it.

\smallskip

The data $(\Cal{F},\nabla ,i^0)$ are unique up to unique
isomorphism.
\endproclaim

\smallskip

{\bf 2.2.2. Comments on the proof.} a). The case when all singularities
are of order 1 is easier. It is treated separately in [Mal3], Th. 2.1.
Since the second structure
connection satisfies this condition, we sketch Malgrange's
argument in this case.

\smallskip

Choose  base points $a\in U:=\bold{P}^1_{\lambda}\setminus\cup_{j=1}^{m+1}\{a^j_0\}$
and $(b_0,a)\in B\times\bold{P}^1_{\lambda}.$ Notice that $(b_0,a)$ belongs to
$V:=B\times\bold{P}^1_{\lambda}\setminus\cup_{j=1}^{m}D_j.$

\smallskip

The restriction of $(\Cal{F}^0,\nabla^0)$ to $U$ is determined uniquely
up to unique isomorphism by the monodromy action
of $\pi_1(U,a)$ on the space $F,$ the geometric fiber  $\Cal{F}^0(a)$
at $a,$ which can be arbitrary. Similarly, there is a bijection
between flat connections $(\Cal{F},\nabla )$ on $V$ with fixed
identification $\Cal{F}^0(a)\to\Cal{F}(a)=F$ and actions of
$\pi_1(V,(a,b))$ on $F.$ Hence to construct an extension  
$(\Cal{F},\nabla )$ to $V$ together with an isomorphism  
 of its restriction to $U$ with $(\Cal{F}^0,\nabla^0),$ 
it suffices to check that $i$ induces an isomorphism $\pi_1(U,a)\to
\pi_1(V,(a,b)),$ which follows from the homotopy exact sequence and the fact
that $B$ is contractible.

\smallskip

This argument explains the term ``isomonodromic deformation.''

\smallskip

Next, we must extend $(\Cal{F},\nabla )$ to $B\times\bold{P}^1_{\lambda}.$
It suffices to do this separately in a tubular neighborhood of each $D_j$
disjoint from other $D_k.$ The coordinate change
$\lambda \mapsto \lambda -a^j$ (or $\lambda\mapsto\lambda^{-1}$)
allows us to assume that the equation of $D_j$ is $\lambda =0.$
Take a neighborhood $W$ of $0$ in which $\Cal{F}^0$
can be trivialized, describe $\nabla^0$ by its connection
matrix, lift $(\Cal{F}^0,\nabla^0)$ to $B\times W$ and restrict to a tubular neighborhood of $D_j.$ On the complement to $D_j,$ this lifting
can be canonically identified with $(\Cal{F},\nabla )$ through their
horizontal sections. Clearly, it is of order $\le 1$ at $D_j.$

\smallskip

It remains to establish that any two extensions are canonically isomorphic.
Outside singularities, an isomorphism exists and is unique.
An additional argument which we omit shows that it extends
holomorphically to $B\times\bold{P}^1_{\lambda}.$

\smallskip

b). When $\nabla$ admits singularity of order 2, this argument must be
completed. The extension of $(\Cal{F}^0,\nabla^0)$ first to
$V$ and then to the singular divisors of order $\le 1$ can be done exactly
as before. But both the existence and the uniqueness of the extension
to the irregular singularities requires an additional
local analysis in order to show that the simple spectrum of the principal polar
part determines the singularity. When formulated in terms of
the asymptotic behaviour of horizontal sections, this analysis
introduces the Stokes data as a version of irregular monodromy,
which also proves to be deformation invariant.

\medskip

{\bf 2.3. The theta divisor and Schlesinger's equations.}
In this subsection we will assume that 
$\Cal{F}^0=T\otimes\Cal{O}_{\bold{P}^1_{\lambda}}$ where $T$
is a finite dimensional vector space which can be identified with the space
of global sections of $\Cal{F}^0.$ This is the case of
the two structure connections, when the local system $\Cal{T}^f_M$
is trivial.

\smallskip

Then there exists a divisor $\Theta ,$ eventually empty, such that
the restriction of $\Cal{F}$ to all fibers $\{b\}\times\bold{P}^1_{\lambda},
\,b\notin \Theta ,$ is free. This can be proved using the fact that
a locally free sheaf $\Cal{E}$ on $\bold{P}^1$ is free iff
$H^0(\bold{P}^1,\Cal{E}(-1))=H^1(\bold{P}^1,\Cal{E}(-1))=0,$
and that the cohomology of fibers is semi--continuous.
For an analytic treatment, see [Mal4], sec.~4 and 5.

\smallskip

Moreover, assume that $\lambda =\infty$ is a singularity of order 1
(to achieve this for the first structure connection, we must replace $\lambda$
by $\lambda^{-1}.$) Then we can identify the inverse image of
$\Cal{F}$ on $B\setminus\Theta\times\bold{P}^1_{\lambda}$
with $T\otimes\Cal{O}_{B\setminus\Theta\times\bold{P}^1_{\lambda}}$
compatibly with the respective trivialization of $\Cal{F}^0.$
To this end trivialize $\Cal{F}$ along $\lambda =\infty$
using the residual connection (see 2.1 iii))and then take the constant
extension of each residually horizontal section along
$\bold{P}^1_{\lambda}.$ (If there are no poles of order 1,
one can extend this argument using a different version of the residual
connection, see [Mal4], p.430, Remarque 1.4.)

\smallskip

Using this trivialization, we can define a meromorphic integrable
connection $\partial$ on $\Cal{F}$ with the space of horizontal sections $T$
on $B\setminus\Theta\times\bold{P}^1_{\lambda}.$
As sections of $\Cal{F}$, they develop a singularity at $\Theta .$
Therefore, the respective connection form $\nabla -\partial$
is a meromorphic matrix one--form with eventual pole at $\Theta .$

\smallskip

The following classical result clarifies the structure of this
form in the case {\it when all poles of $\nabla$ are of order 1.}

\smallskip

\proclaim{\quad 2.3.1. Theorem} a). Let $(a^1,\dots ,a^m)$ be the functions
on $B$ describing the $\lambda$--coordinates of finite poles of $\nabla$
(with given rigidity.) Then
$$
\nabla =\partial +\sum_{i=1}^mA_i(a^1,\dots ,a^m)\,
\frac{d(\lambda -a^i)}{\lambda -a^i}
\eqno{(2.4)}
$$
where $A_i$ are meromorphic functions $B\to \roman{End}\,(T)$
which can be considered as multivalued meromorphic functions
of $a_i.$

\smallskip

b). The connection (2.4) is flat iff $A_i$ satisfy the Schlesinger 
equations
$$
\forall j,\qquad dA_j=\sum_{i\ne j}[A_i,A_j]\,
\frac{d(a^i-a^j)}{a^i -a^j}.
\eqno{(2.5)}
$$
\smallskip

c). Fix a tame point $a_0=(a^1_0,\dots ,a^m_0).$ Then arbitrary
initial conditions $A_i^0=A_i(a_0)$ define a solution of
(2.5) holomorphic on $B\setminus\Theta ,$ with eventual pole at 
$\Theta$ of order 1.

\smallskip

d). For any such solution $\nabla$ of (2.5), define the meromorphic
1-form on $B$:
$$
\omega_{\nabla}:=\sum_{i<j}\roman{Tr}\,(A_iA_j)\,\frac{d(a^i-a^j)}{a^i -a^j}.
\eqno{(2.6)}
$$
This form is closed, and for any local equation $t=0$ of $\Theta$
the form $\omega_{\nabla}-\dfrac{dt}{t}$ is locally holomorphic.
\endproclaim

\smallskip

\proclaim{\quad 2.3.2. Corollary} For any solution $\nabla$ of
(2.5), there exists a holomorphic function $\tau_{\nabla}$
on $B$ such that $\omega_{\nabla}=d\,\roman{log}\,\tau_{\nabla}.$
It is defined uniquely up to a multiplication by a constant.
\endproclaim

\smallskip

In fact, $B$ is simply connected.

\medskip

For a proof of Theorem 2.3.1, we refer to [Mal4]:
a), b), and c) are proved on pp. 406--410,
d) on pp. 420--425.

\medskip

{\bf 2.4. Special solutions.} Slightly generalizing (2.5),
we will call {\it a solution to Schlesinger's equations}
any data $(M,(u^i),T,(A_i))$ where $M$ is a complex
manifold of dimension $m\ge 2$; $(u^1,\dots ,u^m)$ a system
of holomorphic functions on $M$ such that $du^i$ freely generate
$\Omega_M^1$ and for any $i\ne j,\,x\in M$, we have 
$u^i(x)\ne u^j(x)$; $T$ a finite dimensional complex
vector space; $A_j:\,M\to \roman{End}\,T,\,j=1,\dots ,m,$
a family of holomorphic matrix functions such that
$$
\forall j:\qquad dA_j=\sum_{i:\,i\ne j} [A_i,A_j]\,
\frac{d(u^i-u^j)}{u^i-u^j}.
\eqno{(2.7)}
$$

Let such a solution be given. Summing (3.1) over all $j,$
we find $d(\sum_jA_j)=0.$ Hence $\sum_jA_j$ is a constant
matrix function; denote its value by $\Cal{W}.$

\smallskip

\proclaim{\quad 2.4.1. Definition} A  solution to Schlesinger's
equations as above is called special, if 
$\roman{dim}\,T=m=\roman{dim}\,M$; $T$ is endowed with
a complex nondegenerate quadratic form $g$;\
$\Cal{W}=-\Cal{V}-\dfrac{1}{2}\,\roman{Id},$ where
$\Cal{V}\in\roman{End}\,T$ is a skew symmetric operator
with respect to $g$, and finally
$$
\forall j:\qquad A_j=-(\Cal{V}+\dfrac{1}{2}\,\roman{Id})P_j
\eqno{(2.8)}
$$
where $P_j:\,M\to \roman{End}\,T$ is a family of holomorphic
matrix functions whose values at any point of $M$
constitute a complete system of orthogonal projectors
of rank one with respect to $g$:
$$
P_iP_k=\delta_{ik}P_i,\quad
\sum_{i=1}^mP_i=\roman{Id}_T,\quad g(\roman{Im}\,P_i,\roman{Im}\,P_j)={0}
\eqno{(2.9)}
$$
if $i\ne j.$
Moreover, we require that $A_j$ do not vanish at any point of $M$.
\endproclaim

\smallskip

{\bf 2.4.2. Comment.} We committed a slight abuse of language:
the notion of special solution involves a choice of additional
data, the metric $g.$ However, when it is chosen, the
rest of the data is defined unambiguously if it  exists
at all.

\smallskip

In fact, assume that $A_j=\Cal{W}P_j$ as above do not vanish
anywhere. Then they have constant rank one. Hence at any
point of $M$ we have
$$
\roman{Ker}\,A_j=\roman{Ker}\,\Cal{W}P_j=\roman{Ker}\,P_j=
\oplus_{i:\,i\ne j}\roman{Im}\,P_i,
$$
so that
$$
\roman{Im}\,P_i=\cap_{j:\,j\ne i}\oplus_{k:\,k\ne j}\roman{Im}\,P_k=
\cap_{j:\,j\ne i}\roman{Ker}\,A_j.
$$
This means that $P_j$ can exist for given $A_j$ only if
the spaces $\Cal{T}_j=\cap_{j:\,j\ne i}\roman{Ker}\,A_j$
are one--dimensional and pairwise orthogonal at any
point of $M.$

\smallskip

Conversely, assume that this condition is satisfied.
Define $P_j$ as the orthogonal projector onto $\Cal{T}_j.$
Then $A_iP_j=0$ for $i\ne j$ because
$\Cal{T}_j=\roman{Im}\,P_j\subset\roman{Ker}\,A_i.$
Hence
$$
A_j=A_j(\sum_{i=1}^mP_i)=A_jP_j=(\sum_{i=1}^mA_i)P_j=
\Cal{W}P_j.
\eqno{(2.10)}
$$

Notice that all $A_j$ are conjugate to $\roman{diag}\,
(-\dfrac{1}{2},0,\dots ,0)$ and satisfy
$A_j^2+\dfrac{1}{2}\,A_j=0.$ These conditions, as well
as $\sum_jA_j=-(\Cal{V}+\dfrac{1}{2}\,\roman{Id})$, are
compatible with the equations (2.8) and so must be checked
at one point only.

\medskip

{\bf 2.4.3. Strictly special solutions.} A special solution
to Schlesinger's equations as above is called
{\it strictly special} if the operators
$$
A_j^{(t)}:=A_j+tP_j
$$
also satisfy Schlesinger's equations for any $t\in \bold{C}.$

\smallskip

\proclaim{\quad 2.4.4. Lemma} If $\Cal{W}$ is invertible,
then any special solution with given $\Cal{W}$ is
strictly special.
\endproclaim

\smallskip

{\bf Proof.} Inserting $A_j^{(t)}$ into (2.7) one sees
that the solution is strictly special iff
$$
\forall j:\qquad dP_j=\sum_{i:\,i\ne j}
(P_i\Cal{W}P_j-P_j\Cal{W}P_i)\,\frac{d(u^i-u^j)}{u^i-u^j}.
$$
On the other hand, replacing $A_k$ by $\Cal{W}P_k$ in (2.7),
one sees that after left multiplication by $\Cal{W}$
this becomes a consequence of (2.7).

\medskip

{\bf 2.5. From Frobenius manifolds to special solutions.}
Given a semisimple Frobenius manifold with flat
identity and an Euler field $E$ with $d_0=1$, we can produce
a special solution to the Schlesinger equations
rephrasing the results of the previous two sections.

\smallskip

Namely, we first pass to a covering $M$ of the subspace of tame
points of the initial manifold such that $\Cal{T}^f_M$
is trivial and a global splitting can be chosen,
represented by the canonical coordinates $(u^i).$
Then we put $T=\Gamma(M,\Cal{T}^f_M)$ and
$A_i=$ the coefficients of the second structure connection
written as in (2.4).

\smallskip

Since this connection is flat, $(M,(u^i),T,(A_i))$
form a solution of (2.7).

\smallskip

Moreover, this solution is special. In fact, $T$ comes equipped
with the metric $g.$ The operator $A_i$ is the principal
part of order 1 of $\check \nabla$ at $\lambda =u^i$
which is of the form (2.8), with $P_j=e_j\circ$.

\smallskip

Finally, this special solution comes with one more piece of data,
the identity $e\in T$. We will axiomatize its properties in the
following definition.

\smallskip

\proclaim{\quad 2.5.1. Definition} Consider a special solution to 
Schlesinger's equations as in the Definition 2.4.1. A vector
$e\in T$ is called an identity of weight $D$ for this solution,
if

\smallskip

a). $\Cal{V}(e)=(1-\dfrac{D}{2})\,e.$

\smallskip

b). $e_j:=P_j(e)$ do not vanish at any point of $M$.

\endproclaim

\smallskip

For Frobenius manifolds with $d_0=1,$ a) is satisfied
in view of (1.17) and (1.36).

\medskip

Theorem 1.14 shows moreover
that in this way we always obtain strictly special solutions, although the
operator $\Cal{W}$ need not be
invertible. For example, 
for quantum cohomology of $\bold{P}^r$
(which is semisimple, cf. below) the spectrum of $\Cal{W}$
is $\{a-\dfrac{r+1}{2}\,|\,a=0,\dots ,r\}.$ It contains $0$
if $r$ is odd.

\medskip

{\bf 2.6. From special solutions to Frobenius manifolds.}
Let $(M,(u^i),T,g,(A_i))$ be a strictly special solution,
and $e\in T$ an identity of weight $D$ for it.

\smallskip

\proclaim{\quad 2.6.1. Theorem} These data
come from the unique structure of semisimple split
Frobenius manifold on $M$, with flat identity and
Euler field, as it was described in 2.5.
\endproclaim

\smallskip

{\bf Proof.} Proceeding as in 2.5, but in the reverse direction,
we are bound to make the following choices.

\smallskip

Put $e_j=P_j(e)\subset\Cal{O}_M\otimes T,\,j=1,\dots ,m$.
Identify $\Cal{O}_M\otimes T$ with $\Cal{T}_M$ by
setting $e_j=\partial /\partial u^j.$ Transfer the metric
$g$ from $T$ to $\Cal{T}_M.$ Because of 2.5.1 c), it will be
flat, with $T$ as the space of flat vector fields.
Define the multiplication on
$\Cal{T}_M$ for which $e_i\circ e_j=\delta_{ij}e_j.$
Put $\eta_i:=g(e_i,e_i).$

\smallskip

Let $\Cal{T}^f_M$ be the image of $T$ under this identification. 
We will first check that it is an abelian Lie subalgebra of $\Cal{T}_M.$
It will then follow that $g$ is flat, so that we get a structure of
semisimple pre--Frobenius manifold in the sense
of the Definition 1.6.2.

\smallskip

Choose $t\in\bold{C}$ in such a way that
$$
\Cal{W}^{(t)}:=\sum_j\,A_j^{(t)}=\Cal{W}+t\,\roman{Id}\in\roman{End}\,T
$$
is invertible. The section $X=\sum_jf_je_j$ 
of $\Cal{O}_M\otimes T$ lands in $\Cal{T}^f_M$
iff
$$
\Cal{W}^{(t)}X=\sum_jf_j\Cal{W}^{(t)}P_j(e)=
\left(\sum_jf_jA_j^{(t)}\right)(e)\in T.
$$
Let $\nabla$ be the connection on $\Cal{T}_M$ for which
$\Cal{T}^f_M$ is horizontal. Applying it to
$\left(\sum_jf_jA_j^{(t)}\right)(e)$ we see that the last condition is in turn
equivalent to
$$
\forall k:\qquad \sum_j\frac{\partial f_j}{\partial u^k}\,
A_j^{(t)}(e)=
-\sum_jf_j\,\frac{\partial A^{(t)}_j}{\partial u^k}\,(e).
\eqno{(*)}
$$
We can similarly rewrite the condition $Y:=\sum_jg_je_j\in
\Cal{T}^f_M.$

\smallskip

Commutator of vector fields induces on $\Cal{O}_M\otimes T$
the bracket
$$
[X,Y]=\sum_{j,k}\left(f_j\,\frac{\partial g_k}{\partial u^j}- 
g_j\,\frac{\partial f_k}{\partial u^j} \right)\,e_k.
$$
From (*) for $Y$ and $X$ we find:
$$
\sum_{j,k}f_j\,\frac{\partial g_k}{\partial u^j}\,A_k^{(t)}(e)=
-\sum_jf_j\sum_kg_k\,\frac{\partial A_k^{(t)}}{\partial u^j}\,(e),
$$
$$
\sum_{j,k}g_j\,\frac{\partial f_k}{\partial u^j}\,A_k^{(t)}(e)=
-\sum_jg_j\sum_kf_k\,\frac{\partial A_k^{(t)}}{\partial u^j}\,(e).
$$
The terms $j=k$ in the right hand sides are the same.
For $j\ne k,$ using the strict speciality of our solution, we find
$$
f_jg_k\,\frac{\partial A_k^{(t)}}{\partial u^j}=
f_jg_k\,\frac{[A_j^{(t)},A_k^{(t)}]}{u^j-u^k}
$$
so that the  $(j,k)$--term of the first identity
cancels with the $(k,j)$--term of the second one.

\smallskip

To establish that this structure is Frobenius, it suffices to prove that
$e_i\eta_j=e_j\eta_i$ for all $i,j$: see Theorem 1.7.

\smallskip

We have $\eta_j=g(e,e_j).$ Therefore
$$
g(e,A_j(e))=-g(e,(\Cal{V}+\frac{1}{2}\,\roman{Id})\,P_je)=
g(\Cal{V}e,e_j)-\frac{1}{2}\,g(e,e_j)=\frac{1-D}{2}\,\eta_j
\eqno{(2.11)}
$$
since $\Cal{V}$ is skewsymmetric, and $e$ is an eigenvector
of $\Cal{V}.$ Furthermore, let $\nabla$ be the Levi--Civita
connection of the flat metric $g$. Then derivating
(2.11) we find for every $i,j$:
$$
\frac{1-D}{2}\,\frac{\partial}{\partial u^i}\,\eta_j=
g(\nabla_{e_i}(e),A_j(e))+g(e,\nabla_{e_i}(A_j(e)))=
$$
$$
=g(e,\frac{\partial A_j}{\partial u^i}\,(e)),
\eqno{(2.12)}
$$
because $e\in T$ so that $\nabla (e)=0.$ If $i\ne j,$
we find from (2.5)
$$
\frac{\partial A_j}{\partial u^i}=
\frac{[A_i,A_j]}{u^i-u^j}=
\frac{\partial A_i}{\partial u^j}.
\eqno{(2.13)}
$$
This shows that if $D\ne 1,\, e_i\eta_j=e_j\eta_i.$

\smallskip

To see that $D=1$ is not exceptional, one can replace in this argument
$A_j$ by $A_j^{(t)}$ for any $t\ne 0$, so that
$\dfrac{1-D}{2}$ in (2.11) will become 
$\dfrac{1-D}{2}+t.$

\smallskip

It remains to check that $E=\sum_iu^ie_i$ is the Euler field.
According to the Theorem 1.9, we must prove that
$E\eta_j=(D-2)\eta_j$ for all $j$. Insert (2.13) into (2.12)
and sum over $i\ne j$. We obtain:
$$
\frac{1-D}{2}\,E\eta_j=\frac{1-D}{2}\,\sum_{i:\,i\ne j} 
u^i\,\frac{\partial \eta_j}{\partial u^i}+
\frac{1-D}{2}\,u^j\,\frac{\partial \eta_j}{\partial u^j}=
$$
$$
=\sum_{i:\,i\ne j}g\left( e,u^i\frac{[A_i,A_j]}{u^i-u^j}\,(e)\right) +
u^jg(e,\frac{\partial A_j}{\partial u^j}\,(e)).
\eqno{(2.14)}
$$
>From (2.8) it follows that
$$
\frac{\partial A_j}{\partial u^j}=-
\sum_{i:\,i\ne j}\frac{[A_i,A_j]}{u^i-u^j}.
\eqno{(2.15)}
$$
On the other hand,
$$
u^i\,\frac{[A_i,A_j]}{u^i-u^j}=[A_i,A_j]+
u^j\,\frac{[A_i,A_j]}{u^i-u^j}.
\eqno{(2.16)}
$$
Inserting (2.15) and (2.16) into (2.14), we find
$$
\frac{1-D}{2}\,E\eta_j=
\sum_{i:\,i\ne j}g( e,[A_i,A_j]\,(e))+
u^j\sum_{i:\,i\ne j}g\left( e,\frac{[A_i,A_j]}{u^i-u^j}\,(e)\right)+
$$
$$
+u^jg(e,\frac{\partial A_j}{\partial u^j}\,(e))=
g( e,[\sum_{i:\,i\ne j}A_i,A_j]\,(e))=
$$
$$
=-g( e,[\Cal{V}+\frac{1}{2}\,\roman{Id},(\Cal{V}+\frac{1}{2}\,\roman{Id})
P_j]\,(e)).
\eqno{(2.17)}
$$
Using the skew symmetry of $\Cal{V}$,  we see that the last expression
in (2.17) equals $\dfrac{1-D}{2}\,(D-2)\eta_j.$
Hence $E\eta_j=(D-2)\eta_j$ if $D\ne 1.$

\smallskip

Again, replacing in this argument $A_j$ by $A_j^{(t)}$
we see that the restriction $D\ne 1$ is irrelevant.

\medskip

{\bf 2.7. Special initial conditions.} Theorem 2.3.1  shows that
arbitrary initial conditions for the Schlesinger equations
determine a global meromorphic solution on the universal covering
$B(m)$ of $\bold{C}^m\setminus\{\roman{diagonals}\}$,\, $m\ge 2.$

\smallskip

Fix a base point $b_0\in B(m).$ Studying the special solutions,
we may and will identify $T$ with the tangent space at
$b_0$ thus eliminating the gauge freedom. This tangent space
is already coordinatized: we have $e_i$ and $e.$

\smallskip

We will call a family of matrices $A_1^0,\dots ,A_m^0\in
\roman{End}\, T$ {\it special initial conditions}
if we can find a diagonal metric $g$ and a skew symmetric
operator $\Cal{V}$ such that $A_j^0=-(\Cal{V}+\frac{1}{2}\,\roman{Id})\,P_j$
where $P_j$ is the projector onto $\bold{C}\,e_j.$

\smallskip

We will describe explicitly the space $I(m)$ of the special
initial conditions.

\medskip

{\bf 2.7.1. Notation.} Let $R$ be any equivalence relation 
on $\{ 1,\dots ,m\},$ $|R|$ the number of its classes.
Put $F(m)=(\roman{End}\,\bold{C}^m)^m$, Furthermore, denote
$F_R(m)$ the subset of families $(A_1,\dots ,A_m)$ in $F(m)$
such that $R$ coincides with the minimal equivalence relation
for which $iRj$ if $\roman{Tr}\,A_iA_j\ne 0$, and put
$I_R(m)=F_R(m)\cap I(m).$

\medskip

{\bf 2.7.2. Construction.} Denote by $\overline{I}(m)
\subset \bold{C}^m\times\bold{C}^{m(m-1)/2}$
the locally closed subset defined by the equations:
$$
\sum_{i=1}^m\eta_i=0,\quad \eta_i\ne 0\quad \roman{for\ all\ } i;
\eqno{(2.18)}
$$
$$
v_{ij}\eta_j=-v_{ji}\eta_i\quad \roman{for\ all\ } i,j;
\eqno{(2.19)}
$$
$$
\sum_{i=1}^mv_{ij}:=1-\frac{D}{2}\quad  \roman{does\ not\ depend\ on\ } j.
\eqno{(2.20)}
$$
Each point of $\overline{I}(m)$ determines the diagonal metric $g(e_i,e_i)=
\eta_i$ and the operator
$\Cal{V}:\ e_i\mapsto \sum_iv_{ij}e_j$ which is skew symmetric
with respect to $g$ and for which $e$ is an eigenvector.
Setting $A_i=-(\Cal{V}+\dfrac{1}{2}\,\roman{Id})\,P_i$ we get
a point in $I(m).$

\smallskip

This amounts to forgetting $(\eta_i)$ which furnishes the surjective map
$\overline{I}(m)\to I(m)$ because 
$$
A_i(e_j)=0\ \roman{for\ }i\ne j,\quad
A_i(e_i)=-\frac{1}{2}\,e_i-\sum_{j=1}^mv_{ij}e_j.
$$

\smallskip

\proclaim{\quad 2.7.3. Theorem} a). The space $\overline{I}(m)$
can be realized as a Zariski open dense subset
in $\bold{C}^{m+(m-1)(m-2)/2}$.

\smallskip

b). Inverse image in $\overline{I}(m)$ of any point in $I_R(m)$
is a manifold of dimension $1$ for $|R|=1$, $|R|-1$ for $|R|\ge 2.$
\endproclaim  
\smallskip

{\bf Proof.} Fixing $\eta_i$, we can solve (2.19) and (2.20)
explicitly. Put $w_{ij}=v_{ij}\eta_j$ so that
$w_{ij}=-w_{ji}$ and (2.20) becomes
$$
\forall\,j:\qquad \sum_{i=1}^mw_{ij}=\eta_j(1-\frac{D}{2}).
\eqno{(2.21)}
$$
If we choose arbitrarily the values $(w_{ij})$ for all
$1\le i<j\le m-1,$ we can find $w_{mj}$ from the
first $m-1$ equations (2.21), and then the last equations
will hold automatically:
$$
w_{mk}=\eta_k(1-\frac{D}{2})-\sum_{i=1}^{m-1}w_{ik},
$$
$$
\sum_{i=1}^{m}w_{im}=-\sum_{k=1}^{m}w_{mk}=
-\sum_{k=1}^{m-1}\eta_k\,(1-\frac{D}{2})+
\sum_{i,k=1}^{m-1}w_{ik}=
\eta_m\,(1-\frac{D}{2})
$$
because of (2.18).

\smallskip

It remains to determine the fiber of the projection onto $I(m).$

\smallskip

We have for $i\ne j$: $\roman{Tr}\,A_iAj=v_{ij}v_{ji}.$
Hence in the generic case when all these traces
do not vanish, we can reconstruct $\eta_i$
compatible with given $v_{ij}$ from (2.19) uniquely
up to a common factor. Generally, for $i,j$
in the same $R$--equivalence class, (3.12) allows us
to determine the value $\eta_i /\eta_j$  so that
we have $|R|$ overall arbitrary factors constrained by
(3.11).

\medskip

{\bf 2.7.4. Question.} If we choose a special initial condition
for the Schlesinger equation, does the solution remain special
at every point?

\smallskip

Generically, the answer is positive.
If this is the case, we obtain the action of the braid group
$\roman{Bd}_m$ as the group of deck transformations
on the space $I(m).$

\medskip 

{\bf 2.8. Analytic continuation of the potential.}
The picture described in this section gives  a good grip
on the analytic continuation of a germ of
semisimple Frobenius manifold $(M_0,m_0)$ in terms of its
canonical coordinates. Namely, construct the universal
covering $M$ of the subset of the tame points
of $M_0,$ then fix at the point $b_0 = (u^i(m_0))\in B(m)$
the initial conditions of $M$ at $m_0$. This provides
an open embedding $(M,m_0)\subset (B(m),b_0)$.
Loosely speaking, we find in this way a maximal tame
analytic continuation of the initial germ.

\smallskip

Now construct some global flat coordinates $(x^a)$ on $B(m)$
corresponding to a given Frobenius structure.
They map $B(m)$ to a subdomain in $\bold{C}^m$.
This is the natural domain of the analytic continuation of
the potential $\Phi$ of this Frobenius structure,
which is the most important object for Quantum Cohomology.
Unfortunately, its properties are not clear from
this description.  

\newpage

\centerline{\bf \S 3. Quantum cohomology of projective spaces}

\bigskip

In this section we will apply the developed
formalism to the study of the quantum cohomology
of projective spaces $\bold{P}^r,\,r\ge 2.$
 Our main goal is the
calculation of the initial conditions of the relevant solutions
to the Schlesinger's equations. 

\smallskip

{\bf 3.1. Notation.} We start with introducing
the basic notation.
Put $H=H^*(\bold{P}^r,\bold{C})=\sum_{a=0}^r\bold{C}\Delta_a,$
$\Delta_a=$ the dual class of $\bold{P}^{r-a}\subset \bold{P}^r.$
Denote the dual coordinates on $H$ by $x_0,\dots ,x_r$
(lowering indices for visual convenience), $\partial_a=\partial /
\partial x_a .$ The Poincar\'e form is $(g_{ab})=(g^{ab})=(\delta_{a+b,r}).$
The classical (cubic) part of the Frobenius
potential is
$$
\Phi_{\roman{cl}}(x):=\frac{1}{6}\sum_{a_1+a_2+a_3=r}x_{a_1}x_{a_2}x_{a_3}.
\eqno{(3.1)}
$$
The remaining part of the potential is the sum of physicists'
instanton corrections to the self--intersection form:
$$
\Phi_{\roman{inst}}(x):=
\sum_{d=1}^{\infty}\Phi_d(x_2,\dots ,x_r)e^{dx_1},
\eqno{(3.2)}
$$
where we will write $\Phi_d$ as
$$
\Phi_d(x_2,\dots ,x_r)=\sum_{n=2}^{\infty}
\sum_{{a_1+\dots +a_n=}\atop{r(d+1)+d-3+n}}I(d;a_1,\dots ,a_n)\,
\frac{x_{a_1}\dots x_{a_n}}{n!}.
\eqno{(3.3)}
$$ 
This means that if we assign the weight $a-1$ to $x_a,\,a=2,\dots ,n,$ 
$\Phi_d$ becomes the weighted homogeneous polynomial of weight $(r+1)d+r-3.$
Moreover, if we assign to $e^{dx_1}$ the weight $-(r+1),$
$\Phi_{\roman{cl}}$ and $\Phi$
become weighted homogeneous formal series of weight $r-3.$
(Notice that $e$ in the expressions $e^{dx_1}$ and alike is
$2,71828 \dots$, whereas in other contexts $e$ means the identity
vector field. This cannot lead to confusion.)

\smallskip

The starting point of our study in this section will be the following result.

\smallskip

\proclaim{\quad 3.2. Theorem} a). For each $r\ge 2$, there exists a unique
formal solution of the Associativity Equations (1.6) of the form
$$
\Phi (x)=\Phi_{\roman{cl}}(x)+\Phi_{\roman{inst}}(x)
\eqno{(3.4)}
$$
for which $I(1;r,r)=1.$

\smallskip

b). This solution has a non--empty convergence domain in $H$
on which it defines the structure of semisimple
Frobenius manifold $H_{\roman{quant}}(\bold{P}^r)$
with flat identity $e=\partial_0$ and Euler field
$$
E=\sum_{a=0}^{r}(1-a)x^a\partial_a+(r+1)\partial_1
\eqno{(3.5)}
$$
with $d_0=1, D=2-r.$

\smallskip

c). The coefficient $I(d;a_1,\dots ,a_n)$ is the number
of rational curves of degree $d$ in $\bold{P}^r$ intersecting
$n$ projective subspaces of codimensions $a_1,\dots a_n\ge 2$
in general position.
\endproclaim

\smallskip

Uniqueness of the formal solution can be established by showing that
the Associativity Equations imply recursive relations
for the coefficients of $\Phi$ which allow one to express
all of them through $I(1;r,r).$
This is an elementary exercise for $r=2.$
 A more general result (stated
in the language of Gromov--Witten invariants but of essentially
combinatorial nature) is proved in [KM], Theorem 3.1,
and applied to the projective spaces in [KM], Claim 5.2.2.

\smallskip

Existence is a subtler fact. The algebraic geometric
(or symplectic) theory of the Gromov--Witten invariants
provides the numbers $I(d;a_1,\dots a_n)$ satisfying
the necessary relations, together with their
numerical interpretation: see [KM], [BM].
Another approach consists in calculating {\it ad hoc} the
``special initial conditions'' for the semisimple
Frobenius manifold $H_{\roman{quant}}(\bold{P}^r)$
in the sense of the previous section and identifying the
appropriate special solution to the Schlesinger
equations with this manifold. 
For $r=2,$ direct estimates of the coefficients
showing convergence can be found in [D], p.~185.
Probably, they can be generalized to all $r.$

\smallskip

Our approach in this section consists in taking 
Theorem 3.2 for granted and investigating the passage to
the Darboux--Egoroff picture as a concrete
illustration of the general theory. The net outcome
are formulas (3.18) and (3.19) for the special initial
conditions.

\smallskip

Conversely, starting with them, we can construct the Frobenius
structure on the space $B(r+1)$ as was explained in 2.8 above.
Expressing the $E$--homogeneous flat coordinates $(x_0,\dots ,x_r)$
on this space satisfying (3.17)
in terms of the canonical coordinates and then calculating the
multiplication table of the flat vector fields, we
can reconstruct the potential which now will be
a germ of holomorphic function of $(x_a).$ Because of the
unicity, it must have the Taylor series (3.4). So the Theorem 3.2 a),b)
can be proved essentially by reading this section in the
reverse order. Of course, the last statement is of different nature.

\medskip

{\bf 3.3. Tensor of the third derivatives.} Most of our
calculations in $(\Cal{T},\circ )$ will be restricted
to the first infinitesimal neighborhood of the plane
$x_2=\dots =x_r=0$ in $H.$ This just suffices for the
calculation of the Schlesinger initial conditions.
We denote by $J$ the ideal $(x_2,\dots ,x_r).$

\smallskip

Multiplication by the identity $e=\partial_0$ is described
by the components $\Phi_{0a}{}^b=\delta_{ab}$ of the structure
tensor. Of the remaining components, we will need only
$\Phi_{1a}{}^b$ which allow us to calculate multiplication by
$\partial_1,$ and proceed inductively. This is where
the Associativity Equations are implicitly used. 

\smallskip

Obviously, $\Phi_{10}{}^b=\delta_{1b}.$
  
\smallskip
\proclaim{\quad 3.3.1. Claim} We have
$$
\roman{for\ }1\le a\le r-1:\qquad
\Phi_{1a}{}^b=\delta_{a+1,b}+x_{r+1-a+b}e^{x_1}+O(J^2),
\eqno{(3.6)}
$$
$$
\Phi_{1r}{}^b=\delta_{b0}e^{x_1}+x_{b+1}e^{x_1}+O(J^2).
\eqno{(3.7)}
$$
(Here and below we agree that $x_c=0$ for $c>r.$
\endproclaim

\smallskip

{\bf Proof.} The term $\delta_{a+1,b}$ in (3.6) comes from
$\Phi_{\roman{cl}}.$ The remaining terms are provided by
the summands of total degree $\le 3$ in $x_2,\dots ,x_r$
in 
$$
\partial_1\Phi_{\roman{inst}}=
\sum_{d\ge 1}de^{dx_1}\left(
\sum I(d;a_1,a_2)\frac{x_{a_1}x_{a_2}}{2}+
\sum I(d;a_1,a_2,a_3)\frac{x_{a_1}x_{a_2}x_{a_3}}{6}\right)
+O(J^4).
$$
For $n=2,$ the grading condition means that
$d=1,\,a_1=a_2=r.$ For $n=3,$ it means that
$d=1,\,a_1+a_2+a_3=2r+1.$  We know that $I(1;r,r)=1.$
Similarly, $I(1;a_1,a_2,a_3)=1$ in this range.
This can be deduced formally from the Associativity Equations.
A nice exercise is to check that this agrees also
with the geometric description (for instance, only one line
intersects two given generic lines and passes through a 
given point in the three space.) So finally
$$
\partial_1\Phi_{\roman{inst}}=
\left(\frac{x_r^2}{2}+\frac{1}{6}\sum_{a_1+a_2+a_3=2r+1}
x_{a_1}x_{a_2}x_{a_3}\right)e^{x_1}+O(J^4).
$$
The term $\delta_{b0}e^{x_1}$ in (3.7) comes from
$\dfrac{x_r^2}{2}.$ Furthermore,
$$
\Phi_{\roman{inst};1ab} =x_{2r+1-a-b}e^{x_1}+O(J^2)
$$
and
$$
\Phi_{\roman{inst};1a}{}^b =
\Phi_{\roman{inst};1,a,r-b} =x_{r+1-a+b}e^{x_1}+O(J^2).
$$

\medskip

{\bf 3.4. Multiplication table.} The main formula of this subsection is
$$
\partial_1^{\circ (r+1)}=e^{x_1}\left(\partial_0+
\sum_{b=1}^{r-1}(b+1)x_{b+1}\partial_{b}\right) +O(J^2).
\eqno{(3.8)}
$$
We will prove it by consecutively calculating the
powers $\partial_1^{\circ a}$. The intermediate
results will also be used later. (Notice that $O(J^2)$
in (3.8) now means $O(\sum_iJ^2\partial_i).$)

\smallskip

First, we find from (3.6) and (3.7) for $1\le a\le r-1$:
$$
\partial_1\circ\partial_a=\sum_{b=0}^r\Phi_{1a}{}^b\partial_b=
\partial_{a+1}+e^{x_1}\sum_{b=0}^{a-1}
x_{r+1-a+b}\partial_b+O(J^2),
\eqno{(3.9)}
$$
$$
\partial_1\circ\partial_r=\sum_{b=0}^r\Phi_{1r}{}^b\partial_b=
e^{x_1}\left(\partial_{0}+\sum_{b=1}^{r-1}
x_{b+1}\partial_b\right) +O(J^2).
\eqno{(3.10)}
$$
Then using (3.9) and induction, we obtain
$$
\roman{for\ }1\le a\le r:\qquad
\partial_1^{\circ a}=
\partial_{a}+e^{x_1}\sum_{b=0}^{a-2}(b+1)
x_{r+2-a+b}\partial_b+O(J^2).
\eqno{(3.11)}
$$
Multiplying this formula for $a=r$ by $\partial_1$ and using (3.10),
we finally find (3.8).

\smallskip

>From (3.11) it follows that $\partial_1^{\circ a}$ for
$0\le a\le r$ freely span the tangent sheaf.

\medskip

{\bf 3.5. Idempotents.} Formula (3.8) allows us to calculate
all $e_i$ $\roman{mod}\,J^2$ thus demonstrating semisimplicity.
Namely, denote by $q$ the $(r+1)-$th root of the right hand side
of (3.8) congruent to $e^{\frac{x_1}{r+1}}\,\roman{mod}\,J$ and put
$\zeta =\roman{exp}\,(\dfrac{2\pi i}{r+1}).$ Then
$$
e_i=\frac{1}{r+1}\sum_{j=0}^r\zeta^{-ij}(\partial_1\circ q^{-1})^{\circ j}
\eqno{(3.12)}
$$
satisfy
$$
e_i\circ e_j=\delta_{ij}e_i,\ \sum_ie_i=\partial_0
$$
for all $i= 0,\dots r.$ A straightforward check shows this. 

\smallskip

\proclaim{\quad 3.5.1. Proposition} We have
$$
e_i=\frac{1}{r+1}\sum_{j=0}^r\zeta^{-ij}e^{-x_1\frac{j}{r+1}}
\left(e^{x_1}\sum_{b=0}^{j-2}
\frac{(b+1-j)(r+1-j)}{r+1}\,x_{r+b+2-j}\partial_b\right.+
$$
$$
+\left.\partial_j-\sum_{b=j+1}^r\frac{(b+1-j)j}{r+1}\,x_{b+1-j}\partial_b
\right) +O(J^2).
\eqno{(3.13)}
$$
\endproclaim 

\smallskip

{\bf Proof.} We have
$$
q^{-1}=e^{-\frac{x_1}{r+1}}\left(\partial_0-
\sum_{b=1}^{r-1}\frac{b+1}{r+1}x_{b+1}\partial_b\right) +O(J^2).
$$
Together with (3.9) this gives
$$
\partial_1\circ q^{-1}=e^{-\frac{x_1}{r+1}}\left(\partial_1-
\sum_{b=1}^{r-1}\frac{b+1}{r+1}x_{b+1}\partial_{b+1}\right) +O(J^2).
$$
Hence
$$
(\partial_1\circ q^{-1})^j=e^{-\frac{jx_1}{r+1}}\left(\partial_1^{\circ j}-
j\partial_1^{\circ (j-1)}\circ \sum_{b=1}^{r-1}\frac{b+1}{r+1}x_{b+1}\partial_{b+1}\right) +O(J^2).
$$
Inserting this into (3.12) and using (3.9)--(3.11) once again,
we finally obtain (3.13).

\medskip

{\bf 3.6. Metric coefficients in canonical coordinates.} The metric potential
$\eta$ is simply $x_r$ (see (1.11).) Hence we can easily calculate
$\eta_i=e_ix_r.$ The answer is
$$
\eta_i=\frac{\zeta^i}{r+1}e^{-x_1\frac{r}{r+1}}-
\sum_{b=2}^r\frac{\zeta^{ib}}{(r+1)^2}\,b(r+1-b)\,e^{-x_1\frac{r+1-b}{r+1}}x_b
+O(J^2).
\eqno{(3.14)}
$$
As an exercise, the reader can check that the same answer results
from the (longer) calculation of $\eta_i=g(e_i,e_i).$

\medskip

{\bf 3.7. Derivatives of the metric coefficients.} We now see that
the chosen precision   just suffices to calculate
the restriction of $\eta_{ij}$, $\gamma_{ij}$ and the matrix
elements of $A_j$ to the plane $x_2=\dots =x_r=0$
any point of which can be taken as initial one.

\smallskip

\proclaim{\quad 3.7.1. Claim} We have
$$
\eta_{ki}=e_k\eta_i=-2\,\frac{\zeta^{i-k}}{(\zeta^{i-k}-1)^2}\,
\frac{e^{-x_1}}{(r+1)^2}+O(J).
\eqno{(3.15)}
$$
\endproclaim

\smallskip

Notice that (3.15) is symmetric in $i,k$ as it should be.

\smallskip

This is obtained by a straightforward calculation from
(3.13) and (3.14). The numerical coefficient in (3.16)
comes as a combination of $\sum_{j=1}^rj\zeta^j$ and
$\sum_{j=1}^rj^2\zeta^j$ which are then summed by standard tricks.

\medskip

{\bf 3.8. Canonical coordinates.} We find $u^i$ from the formula
$E\circ e_i=u^ie_i.$ To calculate $E\circ e_i,$ use (3.5), (3.13)
and (3.9)--(3.11). We omit the details. The result is:

\smallskip

\proclaim{\quad 3.8.1. Claim} We have
$$
u^i=x_0+\zeta^i(r+1)e^{\frac{x_1}{r+1}}+
\sum_{a=2}^r\zeta^{ai}e^{\frac{ax_1}{r+1}}x_a+O(J^2).
\eqno{(3.16)}
$$
\endproclaim

\smallskip

The reader can check that $e_iu^j=\delta_{ij}+O(j).$

\medskip

{\bf 3.9. Schlesinger's initial conditions.} Recall that the matrix
residues $A_i$ of Schlesinger's equations for Frobenius manifolds
are
$$
A_j(e_i)=0\roman{\ for\ }i\ne j,
$$
$$
A_j(e_j)=-\frac{1}{2}e_j-\frac{1}{2}\sum_k(u^k-u^j)
\frac{\eta_{jk}}{\eta_k}e_k
\eqno{(3.17)}
$$ 
(cf (1.46).) Substituting here (3.14), (3.15) and (3.16), we finally get
the main result of this section.

\smallskip

\proclaim{\quad 3.9.1. Theorem} The point $(x_0,x_1,0,\dots ,0)$
has canonical coordinates
$u^i=x_0+\zeta^i(r+1)e^{\frac{x_1}{r+1}}.$

\smallskip

The special initial conditions at this point (in the sense of 2.7)
corresponding to $H_{\roman{quant}}(\bold{P}^r)$ are given by
$$
v_{jk}=-\frac{\zeta^{j-k}}{1-\zeta^{j-k}}
\eqno{(3.18)}
$$
and
$$
\eta_i=\frac{\zeta^i}{r+1}\,e^{-x_1\frac{r}{r+1}}.
\eqno{(3.19)}
$$
\endproclaim

\smallskip

As an exercise, the reader can check that
$$
-\sum_{k:\,k\ne j}\frac{\zeta^{j-k}}{1-\zeta^{j-k}}=1-\frac{D}{2}
=\frac{r}{2}.
$$

\newpage

\centerline{\bf \S 4. Semisimple  Frobenius supermanifolds}

\bigskip

{\bf 4.1. Supermanifolds and SUSY-structures.} 
A (smooth, analytic, etc.) {\em supermanifold}\ of dimension $(m|n)$ 
is a locally ringed space $(\cM,\f_{\cM})$ with the following
properties [Ma3]: (i) the structure sheaf $\f_{\cM} = \f_{\cM,0} \oplus 
\f_{\cM,1}$ is a sheaf of ${\Bbb Z}_2$-graded supercommutative rings; (ii) 
$\cM_{\text{red}} = (\cM, \f_{\cM, \text{red}}:= 
\f_{\cM}/[\f_{\cM,1} + \f_{\cM,1}^2])$ is a (smooth, analytic, etc.)
classical manifold of dimension $m$;  
(iii)~$\f_{\cM}$ is locally isomorphic to the exterior algebra $\Lambda(E)$
of a locally free $\f_{\cM, \text{red}}$-module $E$ of rank~ $n$.
If $\phi: \Lambda(E)\rar \f_{\cM}$ is any such local isomorphism,
$\bar{x}^1, \ldots, \bar{x}^m$ are local coordinates on  
$\cM_{\text{red}}$ and $\bar{\theta}^{1}, \ldots, \bar{\theta}^{n}$
are free local generators of $E$, then the set of $m+n$ sections
$$
x^1=\phi(\bar{x}^1), \ldots, x^m=\phi(\bar{x}^m), \ 
\theta^1= \phi(\bar{\theta}^1), \ldots, \theta^n=\phi(\bar{\theta}^n)
$$
of the structure sheaf $\f_{\cM}$ form a local coordinate system on $\cM$. 
Any local function $f$ on $\cM$ can be expressed as a polynomial
in anticommuting odd coordinates $\theta^{\al}$, $\al=1,\ldots,n$,
$$
f(x,\theta) = \sum_{k=0}^n \sum_{\al_1, \ldots, \al_k=1}^n f_{\al_1\ldots \al_k}
(x) \theta^{\al_1}\cdots \theta^{\al_k}
$$
whose coefficients $f_{\al_1\ldots \al_k}(x)$ are classical (smooth, analytic,
etc.) functions of the commuting variables $x^a$, $a=1,\ldots,m$.

\smallskip

When a need arises to use odd constants in the structure sheaf of a 
supermanifold $M$, one simply replaces $\cM$ by its relative version, i.e.\ by a 
submersion of supermanifolds $\pi: M \rar S$ whose 
typical fibre is $\cM$. Then (odd) constants
are just (odd) elements of $\pi^{-1}(\f_{S})$. The necessary changes
are routine, see [Ma3] and [Ma4].

\medskip

\proclaim{\quad 4.1.1. Definition} Let $\cM$ be an $(m|n)$-dimensional
supermanifold. A SUSY-structure on $\cM$ is a rank  $0|n$
locally split subsheaf  $\cT_1\subset
\cT\cM$ such that the associated Frobenius form
$$
\matrix
       \Phi: & \Lambda^2 \cT_1 & \lon & \cT_0 := \cT\cM / \cT_1 \\ 
             &  X\ot Y & \lon & \frac{1}{2} [X,Y]\bmod \cT_1
\endmatrix
$$
is surjective. 
\endproclaim

With any SUSY-structure on $\cM$ there is canonically associated
an extension
$$
0 \lon \cT_1 \overset{i}\to\lon \cTM \overset{p}\to\lon \cT_0 \lon 0,
\eqno{(4.1)}
$$
i.e.\ an element ${\frak t}\in \text{Ext}^1_{\f_{\cM}}(\cT_0, \cT_1)
\simeq H^1(\cM, \cT_1\ot \cT_0^*)$.

\smallskip

{\bf 4.1.2. Examples.} 1).  A $(1|1)$-dimensional supermanifold with a
SUSY-structure is called a SUSY$_1$-curve [Ma4]. 2). A SUSY-structure
on a $(3|2)$-dimensional is equivalent to a simple conformal
supergravity in 3 dimensions [Ma3]. 3). A SUSY-structure on a $(4|4)$-dimensional
supermanifold with $\cT_1$ being a direct sum of two integrable
rank $(0|2)$ distributions $\cT_l$ and $\cT_r$ is the same as a
simple superconformal supergravity in 4 dimensions [Ma3]. 

\medskip

{\bf 4.2. Pre-Frobenius supermanifolds.} 
Let $S$ be a module over a
supercommutative ring $R$. A {\it  left odd 
involution}\, on $S$ is by definition a map 
$$
\matrix 
\Pi_l : & S & \lon &  S \\
        & X & \lon & \Pi_lX
\endmatrix
$$
such that $\Pi_l^2=\text{Id}$ and  $\Pi_l(aX)= (-1)^a a \Pi_l(X)$, 
$\Pi_l(Xa) = \Pi_l(X)a$
for any $X\in S$, $a\in R$.

A {\it  right odd involution}\,
$$
\matrix 
\Pi_r :& S & \lon &  S \\
        & X & \lon & X\Pi_r
\endmatrix
$$
also satisfies, by definition, $\Pi_r^2=\text{Id}$ but has
different linearity properties: $(aX)\Pi_r= a
(X\Pi_r)$,  $(Xa)\Pi_r= (-1)^a (X\Pi_r)a$.

\proclaim{\quad 4.2.1. Definition}  
A\, {\it  pre-Frobenius structure} on an 
$(n|n)$-dimensional supermanifold $\cM$ is a quaternary 
$(\cT_1,s, \Pi_l, \Pi_r)$ consisting of a SUSY-structure (4.1),
an even isomorphism $s: \cT\cM \lon \cT_1 \oplus \cT_0$ 
and a pair of left and right odd involutions  
$$
\Pi_{l,r}: \cT_1 \oplus \cT_0 \lon \cT_1 \oplus \cT_0
$$
such that 
\medskip

(i) $\Pi_l(\cT_0)= (\cT_0)\Pi_r= \cT_1$, $\Pi_l(\cT_1)=(\cT_1) \Pi_r = \cT_0$;

\medskip

(ii) a product defined by
$$
  X\circ Y  = \left\{\matrix
                   \Phi(X,Y), &\text{for any }\  X\in \cT_1, Y\in \cT_1, \\ 
                   \Pi_l \Phi(\Pi_l X, Y), &\text{for any }\  X\in \cT_0, Y\in
                   \cT_1,\\ 
                   \Phi(X,Y\Pi_r)\Pi_r, &\text{for any }\  X\in \cT_1, Y\in
                   \cT_0,\\ 
                   \Phi(X\Pi_r, \Pi_l Y), &\text{for any }\  X\in \cT_0, Y\in \cT_0
              \endmatrix\right.
$$
makes $\cT:= \cT_1 \oplus
\cT_0$  a sheaf of associative algebras;

\medskip

(iii) $s\mid_{\cT_1}=i$ and $\left.s\right|_{\cT_0}$ is a
splitting of the extension (4.1), i.e.  $p\circ
\left.s\right|_{\cT_0}=\text{\rm Id}_{\cT_0}$. 

\endproclaim

\smallskip

\proclaim{\quad 4.2.2. Lemma }  There is an even morphism of $\f_{\cM}$-modules,
 $c: \cT \ot \cT \rar \cT$, such that the product $\circ: \cT \times \cT \rar \cT$
factors through the map
$$
\cT \times \cT \lon \cT\ot \cT \overset{c}\to\lon \cT.
$$
\endproclaim

{\bf Proof.} The product $\circ$ obviously satisfies
$(aX)\circ Y = a(X\circ Y)$, $X\circ (a Y)= (Xa)\circ Y$ and
$X\circ(Ya)=(X\circ Y)a$. 

\medskip


{\bf 4.3. Semisimple pre-Frobenius supermanifolds.} 
Let $\cM$ be a pre-Frobenius supermanifold.
\smallskip

\proclaim{\quad 4.3.1. Definition}   A pre-Frobenius structure is called  
 (split)\, {\it almost semisimple}\,  if there exists a local
(global)  basis 
$\{e_{\dal}\}$, $\dal=\dot{1},\ldots,\dot{n}$, of $\cT_1$ called
canonical, such that 
$\Phi(e_{\dal}, e_{\dbe})=\delta_{\dal\dbe} \Pi_l(e_{\dal})
= \delta_{\dal\dbe}(e_{\dal})\Pi_r$.
\endproclaim
\smallskip 

Since in the almost semisimple case $\Pi_l$ completely determines $\Pi_r$ and vice
versa, we can and will omit the subscripts $l,r$. Note that
the definition 4.3.1 makes sense, because the extra condition 
$\Phi(e_{\dal}, e_{\dbe})=\delta_{\dal\dbe} \Pi(e_{\dal})$ is consistent
with 4.2.1(ii) (and in fact implies the latter).
Indeed, in the basis $\{e_{\dal},e_{\al}:=\Pi(e_{\dal})\}$ of 
$\cT_1\oplus \cT_0$ the multiplication $\circ$ takes the form
$$
\align
e_{\dal}\circ e_{\dbe} & =  \delta_{\dal\dbe} e_{\al} ,\\
e_{\dal}\circ e_{\be} & =  e_{\al}\circ e_{\dbe} = \delta_{\al\be} e_{\dal} ,\\
e_{\al}\circ e_{\be} & =  \delta_{\al\be} e_{\al}. 
\endalign
$$
This algebra structure is evidently associative. This table together
with Lemma 4.2.2 immediately imply that
$$
e:= \sum_{\al=1}^n e_{\al}
$$
is the identity, i.e. $e\circ X= X\circ e= X$ for any $X\in \cT$. It
also follows that  
$$
\var := \sum_{\dal=1}^n e_{\dal} = \Pi e
$$
satisfies $\var\circ X = \Pi X$, $X\circ \var = X\Pi$ for any $X\in \cT$.
We call $\var$ the $\Pi$-{\it identity}. 
\medskip

{\bf 4.3.2. Odd multiplication.} If $\cM$ is an almost semisimple pre-Frobenius
supermanifold, then the formulae
$$
  X\bullet Y := X\circ (\Pi Y)=  \left\{\matrix
              \Phi(X,\Pi Y\Pi)\Pi, & \text{for any }\  X\in \cT_1, Y\in \cT_1, \\ 
               \Phi(X\Pi, Y), &\text{for any }\  X\in \cT_0, Y\in
                   \cT_1,\\ 
               \Phi(X,\Pi Y), &\text{for any }\  X\in \cT_1, Y\in
                   \cT_0,\\ 
              \Pi\Phi(\Pi X, \Pi Y), & \text{for any }\  X\in \cT_0, Y\in \cT_0
              \endmatrix\right.
$$
define an odd associative multiplication in $\cT$. Indeed, in the basis
$\{e_{\dal}, e_{\al}\}$ one has 
$$
\align
e_{\dal}\bullet e_{\dbe} & =  \delta_{\dal\dbe} e_{\dal} ,\\
e_{\dal}\bullet e_{\be} & =  e_{\al}\bullet e_{\dbe} = \delta_{\al\be} e_{\al} ,\\
e_{\al}\bullet e_{\be} & =  \delta_{\al\be} e_{\dal}. 
\endalign
$$
The roles of $e$ and $\var$ get interchanged: $\var$ is  the identity,
that is $\var \bullet X= X\bullet \var = X$,
while $e$ is the $\Pi$-identity, that is 
 $e\bullet X = \Pi X$, $X\bullet e= X\Pi$ for any $X\in \cT$. 
\medskip


\proclaim{\quad 4.3.3. Definition} 
 A pre-Frobenius structure on $\cM$ is called  
(split)\, {\it semisimple}\, if 
\medskip

(i) it is (split) almost semisimple;
\medskip

(ii) there is a local (global) coordinate system $\{u^{\al}, \theta^{\dal}\}$ 
called canonical, such that the isomorphism 
$s: \cT_1 \oplus \cT_0 \lon \cT\cM$ is given by
$$
s(e_{\dal})= \p_{\dal} + \theta^{\dal} \p_{\al}, \ \ \ s(e_{\al})= \p_{\al},
$$
where $(e_{\al}, e_{\dal})$ is the canonical basis,
 $\p_{\dal}= \p/\p \theta^{\dal}$ and $\p_{\al}=\p/\p u^{\al}$.

\medskip

(iii) there is an odd metric $g$ on $\cT\cM$ such that
$\cT_1\subset \cTM$ is isotropic and
$$
g(\p_{\al}, \p_{\be})= -\delta_{\al\be}\eta_{\be}, \ \ \ \  
g(e_{\al}, s(e_{\dbe}) )= \delta_{\al\dbe}\eta_{\dal},  
\eqno{(4.2)}
$$
where $\eta_{\al}= \p_{\al} \Psi$,  $\eta_{\dal}= (\p_{\dal} + \theta^{\dal}
\p_{\al})\Psi$
and $\Psi$ is an odd function. Such a metric is called an\, {\it
Egoroff} metric.
\endproclaim

\smallskip

Since the restriction of $s: \cT_1 \oplus \cT_0 \lon \cT\cM$ to $\cT_1$
coincides with $i$ and hence is rigidly fixed by the choice of a SUSY structure
on $\cM$, we identify from now on
$X$ and $s(X)$ for any $X\in \cT_1$. In particular, whenever
the pre-Frobenius structure is semisimple, we always assume that
$e_{\dal}= \p_{\dal} + \theta^{\dal}\p_{\al}$. 

\smallskip

Note that canonical coordinates are defined up to a transformation
$$
\align
\theta^{\dal} & \lon  \hat{\theta}^{\dal} = \theta^{\dal} + c^{\dal} \\
u^{\al} & \lon  \hat{u}^{\al}= u^{\al}+  c^{\al} + \theta^{\dal}c^{\dal} 
\endalign
$$
which satisfy $\hat{\p}_{\al}= \p_{\al}$, $\hat{e}_{\dal}= e_{\dal}$
and hence leave all the defining relations invariant.

\medskip

{\bf 4.3.4. Example.} If $\cM$ is a SUSY$_1$-curve, then any 
odd isomorphism $\cT_1 \rar \cT_0$ together with a
global nowhere vanishing section of $\cT_1$ (if any) and an odd 
function $\Psi$ equips $\cM$ with a semisimple pre-Frobenius structure. 

\medskip

{\bf 4.3.5. Algebra structure in $\cT\cM$.}  The isomorphism
$s: \cT \rar \cT\cM$ translates the product $\circ$ from
$\cT$ to $\cT\cM$ which we denote by the same symbol. The element
$\bar{e}:= s(e) =\sum_{\al} \p_{\al}$ is the identity in $(\cT\cM,\circ)$.

\smallskip

It is easy to check that in the basis
$\{\p_{\dal}, \p_{\al}\}$ the induced multiplication $\circ$ takes the form
$$
\align
\p_{\dal}\circ \p_{\dbe} & =  \delta_{\dal\dbe}\p_{\al} ,\\
\p_{\dal}\circ \p_{\be} & =  \p_{\al}\circ \p_{\dbe} = \delta_{\al\be} \p_{\dal} ,\\
\p_{\al}\circ \p_{\be} & =  \delta_{\al\be} \p_{\al}, 
\endalign
$$
implying that the element $\bar{\var} = \sum_{\dal=1}^n \p_{\dal}$ is the
$\bar{\Pi}$-identity in $(\cT\cM, \circ)$, where the odd automorphism 
$\bar{\Pi}: \cT\cM \rar \cT\cM$ is defined by
$$
\bar{\Pi}(\p_{\dal}) = \p_{\al}, \ \ \  \bar{\Pi}(\p_{\al}) = \p_{\dal}.
$$

\medskip

{\bf 4.3.6. Representation of the Neveu-Schwarz Lie superalgebra in $\cTM$.}
Let $\cM$ be a semisimple pre-Frobenius supermanifold.
Consider vector fields
$$
\align
E & = \sum_{\al=1}^n (u^{\al}\p_{\al} + \frac{1}{2} \theta^{\dal} e_{\dal}) \\
F & = \sum_{\al=1}^n u^{\al} e_{\dal} \\
\endalign
$$
on $\cM$. A direct calculation shows that the vector fields
$$
\align
{\frak e}_a & :=  E^{\circ(a+1)}  =  \sum_{\al=1}^n \left[
(u^{\al})^{a+1} \p_{\al} + \frac{a+1}{2}(u^{\al})^{a}
\theta^{\dal} e_{\dal}\right], \  a=0,1,2,\ldots \\	
{\frak f}_{i+1/2} & := \Pi^{i}F^{\circ(i+1)} =  \sum_{\al=1}^n (u^{\al})^{i+1}
e_{\dal}, \ i=0,1,2,\ldots
\endalign
$$
satisfy the following commutation relations
$$
\align
 \left[{\frak e}_a, {\frak e}_b\right] & =  (b-a) {\frak e}_{a+b}, \\
 \left[{\frak e}_a, {\frak f}_{i}\right] & =  (i - \frac{a}{2}) {\frak f}_{i+a}, \\
 \left[{\frak f}_{i}, {\frak f}_{j}\right] & =  2 {\frak e}_{i+j}, 
\endalign
$$
where $a,b=0,1,2,\ldots$ and $i,j=\frac{1}{2}, \frac{3}{2}, \frac{5}{2},\ldots$. 

\medskip


{\bf 4.4. Flat connections on pre-Frobenius supermanifolds}.
On a supermanifold $\cM$ with a SUSY-structure 
one may define the differential operator $\delta: \f_{\cM}\rar \cT_1^*$ as the
composition
$$
\delta: \f_{\cM} \lon \OM \overset{i^*}\to\lon \cT_1^*.
$$

Let $\cM$ be a semisimple pre-Frobenius
supermanifold with the associated commutative  diagram
$$
\CD
@.@. 0 @. \\
@.@. @VVV @. \\
@.@. \cT_0^* @. \\
@.@. @Vp^*VV @. \\
0 @>>> \f_{\cM} @>d>> \OM @>d>> \Omega^2\cM   \\
@.@. @AsAA  @VV{\ga}V \\
@.@. \cT_1^* @>\delta'>> \Lambda^2 \cT_1^*\\
@.@. @AAA  @VVV\\
@.@. 0 @. 0 \\
\endCD
$$
and let $\Phi^*: \cT_0^* @>{\ga\circ d \circ p^*}>>\Lambda^2 \cT_1^*$
be the (dual) Frobenius form of the SUSY-structure. Then one has

\smallskip

\proclaim{\quad 4.4.1. Proposition}  There is a one-to-one correspondence 
between locally flat connections in a locally free sheaf $\cF$
on $\cM$ and covariant differentials
$$
\nabla: \cF\lon  \cT_1^*\ot \cF
$$
such that 
$$
\nabla(fv)= \delta(f)v + f\nabla(v), \ \ \ \ \ \forall f\in \f_{\cM}, v\in
\cF, 
$$
and the composition 
$$
C(\nabla^2): \cF\overset{\nabla}\to\lon  \cT_1^*\ot \cF \overset{\nabla}\to\lon
\Lambda^2\cT_1^* \ot \cF \lon
\Lambda^2\cT_1^*\ot \cF / \Phi^*(\cT_0)\ot \cF 
$$
(which is $\f_M$-linear) is zero, where the action of $\nabla$ 
on $\cT_1^*\ot \cF$ is defined as follows
$$
\nabla(t\ot v) = (-1)^{\tilde{t}} \ga\ot\text{\rm Id}_{\cF}(t\ot \nabla(v)) +
\delta'(t)\ot v, \ \ \ \ \ \ t\in \cT_1^*, v\in \cF.
$$
\endproclaim
\smallskip

{\bf Proof}. Given a linear connection $D: \cF \rar \OM\ot \cF$,
one defines $\nabla$ as the composition
$$
\nabla: \cF \overset{D}\to\lon \OM \ot \cF @>{i^*\ot \text{Id}}>> \cT_1^*\ot \cF.
$$

Since the curvature tensor, $F\in \Omega^2\cM\ot \cF\ot \cF^*$, of $D$ satisfies
$$
F(X,Y)v = [D_X,D_Y]v - D_{[X,Y]}v
$$
for any $X,Y\in \cT\cM$ and $v\in \cF$, 
the composition $\nabla^2:=\nabla\circ \nabla$,
viewed as a morphism $\Lambda^2 \cT_1 \ot \cF \rar  \cF$, can be written 
explicitly as
$$
\matrix
\nabla^2:& \Lambda^2 \cT_1 \ot \cF & \lon & \cF \\
&                         X\ot Y \ot v      & \lon &
 D_{2\Phi(X,Y)} v + F(X,Y) v
\endmatrix
$$
implying that $C(\nabla^2)$ is essentially $\ga\ot \text{Id}_{\cF\ot \cF^*}(F)$.
Hence $C(\nabla^2)$
is always $\f_{\cM}$-linear and vanishes when $D$ is flat.

\smallskip

In the other direction, let $\nabla: \cF\rar \cT_1^*\ot \cF$
 be a covariant differential such that $C(\nabla^2)=0$. Then $\nabla^2$
factors through the composition
$$
\nabla^2: \Lambda^2 \cT_1 \ot \cF @>{\Phi\ot \text{Id}}>> 
\cT_0\ot \cF \overset{\nabla^0}\to\lon \cF,
$$
for some covariant differential operator $\nabla^0: \cF\rar \cT_0^*\ot \cF$.
Define 
$$
D: \cF\rar \OM\ot \cF \simeq \cT_1^*\ot \cF \ \oplus \ \cT_0^*\ot \cF
$$
as $\nabla\oplus \nabla^0$. A simple calculation in the canonical
coordinates (which we omit) shows that $D$ is flat. This completes the
proof. 
\medskip

Note that in the presence of a SUSY-structure with invertible
Frobenius form {\em any}\ covariant differential $\nabla: \cF\rar \cT_1^*\ot
\cF$ can be canonically extended to a  linear connection $D: \cF\rar \OM \ot
\cF$ as $\nabla\oplus \nabla^0$, where 
$$
\nabla^0: \cF @>{\nabla^2}>>\Lambda^2 \cT_1^* \ot \cF
 @>{{\Phi^*}^{-1}}>> \cT_0^*\ot \cF.
$$
More generally, such an extension is possible whenever $\cM$ comes
equipped with a monomorphism $\Theta: \cT_0 \rar \Lambda^2 \cT_1$
satisfying $\Theta \circ \Phi = \text{Id}$. For example, 
if $\cM$ is almost semisimple pre-Frobenius, then 
$$
\matrix
\Theta: & \cT_0 & \lon & \Lambda^2 \cT_1 \\
& e_{\al} & \lon & e_{\dal}\ot e_{\dal}
\endmatrix
$$
does have this property. When we call a covariant differential 
$\nabla: \cF\rar \cT_1^*\ot \cF$ a connection on $\cF\rar \cM$, we mean precisely this 
SUSY extension. 

\medskip

{\bf 4.5. Semisimple Frobenius structures.}  Let $(\cM, 
\cT_1,s,\Pi,g)$ be a semisimple pre-Frobenius
supermanifold.

\smallskip

\proclaim{\quad 4.5.1. Definition} $\cM$ is called \,  
{\it semisimple Frobenius} if $g$ is flat.

\endproclaim

\medskip

{\bf 4.5.2. Levi-Civita connection.} 
If $g$ is an odd metric on a supermanifold $\cM$, 
$g_{AB} := (-1)^{\tilde{A}}g(e_A,e_B)$ are the components of $g$
 in a basis $e_A$ of $\cT\cM$, then 
the Christoffel symbols, $\nabla_{e_A} e_B = \sum_{C}
\Gamma_{AB}^C e_C$, of the associated Levi-Civita connection $\nabla$ 
are given by
$$
\align
\Gamma_{AB}^{C} & =  \frac{1}{2} \sum_{D} 
\left[ e_A g_{BD} +
 (-1)^{AB} e_{B} g_{AD} - (-1)^{D(A+B+1)+B} e_D g_{AB} \right. \\
& + \sum_M(\left. C_{AB}^{M}g_{MD} - (-1)^{BD +B +D} C_{AD}^M g_{MB} - 
(-1)^{A(B+D+1)+D} C_{BD}^M g_{MA})\right] g^{DC}, 
\endalign
$$
where $C_{AB}^M$ are defined by 
$$
[e_A, e_B]= \sum_M C_{AB}^M e_M.
$$

\smallskip

Consider now  a special case when $g$ is an Egoroff metric on a
semisimple pre-Frobenius supermanifold $\cM$.

\smallskip

\proclaim{\quad 4.5.3. Proposition-definition (Darboux-Egoroff equations)}  
The metric  $g$ is flat if and only if $\Psi$
satisfies the equations
$$
e_{\dmu}\gamma_{\dal\dbe} =   \gamma_{\dmu\dal}\gamma_{\dmu\dbe} \ \
\text{for all}\ \dmu\neq \dal\neq \dbe\neq \dmu, 
\eqno{(4.3)}
$$
$$
\bar{e}\gamma_{\dal\dbe}  =  0 \ \ \ \  \text{for all}\ 
\dal\neq \dbe ,
\eqno{(4.4)}
$$
where
$$
\gamma_{\dal\dbe} = 
\frac{e_{\dal}\eta_{\dbe}}{2\sqrt{\eta_{\dal}\eta_{\dbe}}}.
$$

These equations are called\, {\it Darboux-Egoroff} equations.
\endproclaim

\smallskip

{\bf A sketch of the proof.} The only non-trivial commutator of basis vector
fields $(e_{\dal},\p_{\al})$ is 
$[e_{\dal},e_{\dbe}]=2\delta_{\dal\dbe}\p_{\al}$ so that the
only non-vanishing components of $C_{AB}^M$ are $C_{\dal\dbe}^{\ga}= 
2\delta_{\dal\dbe}\delta_{\dal\ga}$. Then, using the above  formulae
for $\Gamma_{AB}^C$, one obtains
the following Cristoffel symbols of the Levi-Civita connection of $g$: 

\smallskip

$$        
\nabla_{e_{\dmu}} e_{\dal}  =  \delta_{\dmu\dal}\p_{\al} +
\frac{e_{\dmu}\eta_{\dal}}{2\eta_{\dal}} e_{\dal} - 
\frac{e_{\dal}\eta_{\dmu}}{2\eta_{\dmu}} e_{\dmu}, \eqno{(4.5)}
$$
$$
\align
\nabla_{e_{\dmu}} \p_{\al}  &  = 
e_{\dmu}\left(\frac{\eta_{\al}}{2\eta_{\dal}}\right) e_{\dal}
+ \frac{e_{\al}\eta_{\dmu}}{2\eta_{\dmu}}e_{\dmu}
- \delta_{\mu\al}\sum_{\dbe} e_{\dbe}
\left(\frac{e_{\dbe}\eta_{\dmu}}{2\eta_{\dbe}}\right) e_{\dbe} \\
& \ \ \  + \ \frac{e_{\dmu}\eta_{\dal}}{2\eta_{\dal}}\p_{\al}
-\delta_{\mu\al}\frac{\eta_{\al}}{\eta_{\dal}}\p_{\al}
 +  \delta_{\mu\al}\sum_{\dbe}
\frac{e_{\dbe}\eta_{\dal}}{2\eta_{\dbe}}\p_{\be}. \tag{4.6} 
\endalign
$$
 
By Proposition 4.4.1, the connection $\nabla$ is flat if
and only if
$$ 
[\nabla_{e_{\dmu}}, \nabla_{e_{\dnu}}]=0 \ \ \ \ \text{for all} \ 
\dmu\neq \dnu.
$$
A straightforward but very tedious calculation shows that the latter equations
are equivalent to the Darboux-Egoroff equations (4.3) and (4.4).

\medskip

{\bf 4.6. Flat identity.} We say that the $\circ$-identity $e$
on a semisimple pre-Frobenius supermanifold $\cM$  is {\it flat}\
if $\nabla \bar{e}=0$, where $\nabla$ is the Levi-Civita covariant differential
(4.5)-(4.6).

\smallskip
 
\proclaim{\quad 4.6.1. Proposition}  The identity $e$ is flat
if and only if the potential $\Psi$
satisfies the equations
$$
\bar{e}\eta_{\dal}=0, 
\eqno{(4.7)}
$$
or, equivalently,
$$
\sum_{\al}\eta_{\al}= \text{\rm const}. 
\eqno{(4.8)}
$$
\endproclaim

\smallskip

{\bf Proof.}  It follows from (4.6) that
$$
\align
\nabla_{e_{\dmu}} \bar{e} & =  \sum_{\be} \nabla_{e_{\dmu}} e_{\be} \\
& =
\sum_{\dbe}e_{\dmu}\left(\frac{\eta_{\be}}{2\eta_{\dbe}}\right) e_{\dbe}
+ \sum_{\be} \frac{e_{\be}\eta_{\dmu}}{2\eta_{\dmu}}e_{\dmu}
- \sum_{\dbe} e_{\dbe}
\left(\frac{e_{\dbe}\eta_{\dmu}}{2\eta_{\dbe}}\right) e_{\dbe} \\
& \ \ \  + \ \sum_{\be}\frac{e_{\dmu}\eta_{\dbe}}{2\eta_{\dbe}}\p_{\be}
- \frac{\eta_{\mu}}{\eta_{\dmu}}\p_{\mu}
 + \sum_{\be}\frac{e_{\dbe}\eta_{\dmu}}{2\eta_{\dbe}}\p_{\be} \\
&=  \frac{\sum_{\be}e_{\be}\eta_{\dmu}}{2\eta_{\dmu}}e_{\dmu},
\endalign
$$
where we used the fact that $e_{\dmu}\eta_{\dal} + e_{\dal}\eta_{\dmu}=
2\delta_{\dmu\dal} \eta_{\al}$. 

\smallskip

\proclaim{\quad 4.6.2. Proposition} 
 Let $\cM$ be a semisimple pre-Frobenius supermanifold. The $\bar{\Pi}$-identity 
$\bar{\var}$ is flat, i.e.\ $\nabla \bar{\var}=0$ for $\nabla$ being
the Levi-Civita covariant differential,
if and only if the potential $\Psi$
satisfies the equations
$$
(\theta^{\dmu} - \theta^{\dal}) e_{\dmu}\eta_{\dal}  =  0, \eqno({4.9}) 
$$
$$
\bar{\var} \eta_{\dal}  =  \eta_{\al}. \eqno{(4.10)}
$$
\endproclaim

\smallskip

{\bf Proof}\ is a straightforward calculation. 

\medskip


{\bf 4.7  Euler field.} We want to introduce the notion of homogeneity
of a Frobenius structure by assigning the scaling degrees $1$ and $1/2$
to the canonical coordinates $u^{\al}$ and $\theta^{\dal}$ respectively
(reflecting the fact that $\p_{\al} = e_{\dal}e_{\dal}$).  With this
motivation, we define a {\em scaling field}\ on a Frobenius supermanifold $\cM$
as an even vector field $E$ satisfying 
$$
[E, \p_{\al}] = -\p_{\al}, \ \ \ \ [E, e_{\dal}]= -\frac{1}{2} e_{\dal}.
$$

\smallskip

\proclaim{\quad 4.7.1. Proposition}  If $E$ is a scaling field, then
$$
E=\sum_{\al}\left[(u^{\al} + c^{\al} + \theta^{\dal}c^{\dal}) \p_{\al}   
+ \frac{1}{2}({\theta^{\dal}} + c^{\dal}) e_{\dal}\right]
$$
for some even constants $c^{\al}$ and odd constants $c^{\dal}$.
\endproclaim

\smallskip

{\bf Proof.} Putting $E=\sum_{\al} E^{\al}\p_{\al} + E^{\dal}e_{\dal}$,
one obtains 
$$
\align
\left[E, \p_{\al}\right] = -e_{\al}  & \Longleftrightarrow 
\sum_{\be}\left[(\p_{\al}E^{\be})\p_{\be} +
(\p_{\al}E^{\dbe}) e_{\dbe}\right] = \p_{\al},  \\
\left[E, e_{\dal}\right] = -\frac{1}{2} e_{\dal} &  \Longleftrightarrow 
\sum_{\be}\left[(e_{\dal}E^{\be})\p_{\be} +
(e_{\dal}E^{\dbe}) e_{\dbe}\right] - 2 E^{\dal}\p_{\al} 
= \frac{1}{2} e_{\dal} ,
\endalign
$$
implying $\p_{\al}E^{\be}= \p_{\al} E^{\dbe} = e_{\dal} E^{\be} = e_{\dal}
E^{\dbe}=0$ for all $\al\neq \be$ as well as  $\p_{\al}E^{\dal} =0$,
$e_{\dal}E^{\al} = 2E^{\dal}$ and $e_{\dal}E^{\dal}=1/2$. Hence
$E^{\dal}= \frac{1}{2}(\theta^{\dal} + c^{\dal})$ and $E^{\al} =  
u^{\al} + c^{\al} + \theta^{\dal}c^{\dal}$
for some even constants $c^{\al}$ and odd constants $c^{\dal}$. 

\smallskip

Given a scaling field $E$,  we can and will normalize the canonical
coordinates so that $E=\sum_{\al}\left[u^{\al}\p_{\al} + 
\frac{1}{2}{\theta^{\dal}}e_{\dal}\right]$. 

\medskip

\proclaim{\quad 4.7.2. Definition}  A scaling vector field $E$ on a 
semisimple pre-Frobenius supermanifold $\cM$ is called an {\em Euler field}\ 
if $\text{Lie}_E (g) = Dg$ for some constant $D$, that is,
$$
E(g(X,Y)) - g([E,X],Y) - g(X,[E,Y])  =  Dg(X,Y), 
\eqno{(4.11)}
$$
for all vector fields $X$, $Y$. 
\endproclaim

\smallskip

\proclaim{\quad 4.7.3. Proposition}  If $E$ is an Euler field on a 
semisimple pre-Frobenius supermanifold $\cM$, then the potential $\Psi$ satisfies 
$$
E \Psi = (D-1) \Psi + \ \text{\rm const},
$$
or, equivalently
$$
E\eta_{\dal}  =  (D-\frac{3}{2}) \eta_{\dal}. 
\eqno{(4.12)}
$$
\endproclaim

{\bf Proof.} Write
(4.11) for  (i) $X=\p_{\al}$, $Y=\p_{\be}$,  (ii) 
$X=e_{\dal}$, $Y=\p_{\be}$ and (iii) $X=e_{\dal}$, $Y=e_{\dbe}$ 
to find that
$$
\align
(i) & \Leftrightarrow  \  Eg(\p_{\al}, \p_{\be}) + 2 g(\p_{\al},\p_{\be}) = D
g(\p_{\al}, \p_{\be}),  \\
(ii) & \Leftrightarrow \  Eg(e_{\dal}, \p_{\be}) +\frac{3}{2}
g(e_{\dal},\p_{\be})  = D g(e_{\dal}, \p_{\be}),\\
(iii) & \Leftrightarrow \  0=0,
\endalign
$$
which imply equation~(4.12) and
$$
E\eta_{\al}  =  (D-2) \eta_{\al}. \eqno{(4.13)}
$$
However, the latter equation is not independent ---
applying $e_{\dbe}$ to both sides of (4.12), one easily obtains
$E(e_{\dbe}\eta_{\dal}) = (D-2) e_{\dbe}\eta_{\dal}$ implying 
(4.13).  

\smallskip

\proclaim{\quad 4.7.4. Corollary}  If $E$ is an Euler field on a 
semisimple pre-Frobenius
supermanifold $\cM$, then 
$$
E\gamma_{\dmu\dnu} = - \frac{1}{2}\gamma_{\dmu\dnu}.
\eqno{(4.14)}
$$
\endproclaim

\smallskip

{\bf Proof.}
$$
\align
E\gamma_{\dmu\dnu} & = 
E\left(\frac{e_{\dmu}\eta_{\dnu}}{2\sqrt{\eta_{\dmu}\eta_{\dnu}}}\right) \\
& =  \frac{(Ee_{\dmu}\eta_{\dnu})}{2\sqrt{\eta_{\dmu}\eta_{\dnu}}}
- \frac{(E\eta_{\dmu})(e_{\dmu}\eta_{\dnu})}
{4\eta_{\dmu}\sqrt{\eta_{\dmu}\eta_{\dnu}}} 
- \frac{(E\eta_{\dnu})(e_{\dmu}\eta_{\dnu})}
{4\eta_{\dnu}\sqrt{\eta_{\dmu}\eta_{\dnu}}} \\
& =  (D-2)\gamma_{\dmu\dnu}  - \frac{1}{2}(D-\frac{3}{2})\gamma_{\dmu\dnu} 
- \frac{1}{2}(D-\frac{3}{2})\gamma_{\dmu\dnu}  \\
& =  - \frac{1}{2}\gamma_{\dmu\dnu}. 
\endalign
$$

\medskip

{\bf 4.7.5. Nullness of the $\circ$-identity}. In the presence of an
Euler field, the flatness of $e$ implies that either  $D=2$ or 
$g(\bar{e},\bar{e})$ is identically zero, i.e.\ the identity $\bar{e}\in
(\cT\cM,\circ)$ is everywhere a null vector. Indeed, 
$$
g(\bar{e},\bar{e}) = g(\sum_{\al}\p_{\al}, \sum_{\be}\p_{\be}) =
\sum_{\al}\eta_{\al},
$$
while equations (4.8) and (4.13) imply
$(D-2)\sum_{\al}\eta_{\al} = 0$. 

\medskip

{\bf 4.8. Geometry on $\cT_1$.} 
Let $\cM$ be a semisimple pre-Frobenius
supermanifold. Define a new splitting 
 of the SUSY-extension (4.1) as follows
$$
\matrix
\tilde{s}: & \cT_0 & \lon & \cT\cM, \\
& \p_{\al}\bmod \cT_1 & \lon &  \tilde{e}_{\al}:= \p_{\al}
 - \sum_{\dbe}\frac{e_{\dbe}\eta_{\dal}}{2\eta_{\dbe}} e_{\dbe}.
\endmatrix
$$
\smallskip

\proclaim{\quad 4.8.1. Lemma}  The splitting $\ts$ decomposes $\cT\cM$
into a  direct sum $\cT_1\oplus
\tilde{s}(\cT_0)$ of isotropic submodules.
\endproclaim

\smallskip

{\bf Proof.} The restriction of the Egoroff metric to
$\ts(\cT_0)$ is 
$$
\align
g(\te_{\al}, \te_{\be}) &\ = \ g(\p_{\al} 
 - \sum_{\dmu}\frac{e_{\dmu}\eta_{\dal}}{2\eta_{\dmu}} e_{\dmu},
\p_{\be}  - \sum_{\dnu}\frac{e_{\dnu}\eta_{\dal}}{2\eta_{\dnu}} e_{\dnu})
\\
&\  = \ -\delta_{\al\be}\eta_{\al} + \frac{e_{\dbe}\eta_{\dal}}{2}
+ \frac{e_{\dal}\eta_{\dbe}}{2} \\
& \ =  \ 0.
\endalign
$$

Note for future reference that $g(\ts(e_{\al}), e_{\dbe})=g(s(e_{\al}),
e_{\dbe})= \delta_{\al\be} \eta_{\dal}$. 

\medskip

{\bf 4.8.2. Distinguished connections on $\cT_1$ and $\cT_0$.}
The decomposition  $\cT\cM =\cT_1\oplus
\tilde{s}(\cT_0)$ induces a  projection
$$
p_1: \cT\cM  \rar  \cT_1. 
$$
 Then, if $\nabla: \cT\cM \rar \cT_1^*\ot \cT\cM$ is the Levi-Civita
covariant differential of the Egoroff metric,
one may define the operators 
$$
\tnabla^1: \cT_1 \rar \cT_1^* \ot\cT_1,  \ \ \ \ \ \ \ \ 
\tnabla^0: \cT_0 \rar \cT_1^* \ot\cT_0 
$$
as the compositions 
$$
\align
\tnabla^1: & \ \cT_1 \overset{i}\to\lon i(\cT_1) \overset{\nabla}\to\lon 
\cT_1^* \ot \cT\cM @>{\text{Id}\ot p_1}>>\cT_1^*\ot \cT_1, \\
\tnabla^0: & \ \cT_0 \overset{\ts}\to\lon \ts(\cT_0)\overset{\nabla}\to\lon 
\cT_1^* \ot \cT\cM @>{\text{Id}\ot p}>>\cT_1^*\ot \cT_0, 
\endalign
$$
where $i$ and $p$ are defined in (4.1).

\smallskip

Remarkably, the connections $\tnabla^1$ and $\tnabla^0$ are essentially
one and the same thing: 

\smallskip

\proclaim{\quad 4.8.3. Lemma} $\tnabla^1 (X\Pi) = (\tnabla^0 X)\Pi$, \ 
$\tnabla^0 (Y\Pi) = (\tnabla^1 Y)\Pi$,   {\em for any}\ $X\in \cT_0$,
$Y\in \cT_1$.
\endproclaim

\smallskip

{\bf Proof}.  Since 
$$
p_1(e_{\al})= \sum_{\dbe}\frac{e_{\dbe}\eta_{\dal}}
{2\eta_{\dbe}} e_{\dbe},
$$
one obtains from (4.5)
$$
\tnabla_{e_{\dmu}} e_{\dal} = 
\delta_{\dmu\dal}\sum_{\dbe}\frac{e_{\dbe}\eta_{\dal}}{2\eta_{\dbe}}
e_{\dbe} +
\frac{e_{\dmu}\eta_{\dal}}{2\eta_{\dal}}e_{\dal} -
\frac{e_{\dal}\eta_{\dmu}}{2\eta_{\dmu}}e_{\dmu}. \eqno{(4.15)}
$$ 

Analogously,
$$
\align
\tnabla_{e_{\dmu}} e_{\al}& \ = \ p\left(\nabla_{e_{\dmu}}(
\p_{\al} - \sum_{\dbe}\frac{e_{\dbe}\eta_{\dal}}{2\eta_{\dbe}} e_{\dbe})\right)
\\
& \ = \ p(\nabla_{e_{\dmu}}\p_{\al}) + 
\sum_{\dbe}\frac{e_{\dbe}\eta_{\dal}}{2\eta_{\dbe}}
p(\nabla_{e_{\dmu}} e_{\dbe}) \\
&\ = \ 
\frac{e_{\dmu}\eta_{\dal}}{2\eta_{\dal}}e_{\al} +
\delta_{\dmu\dal}\sum_{\dbe}\frac{e_{\dbe}\eta_{\dal}}{2\eta_{\dbe}}e_{\be}
- \delta_{\dmu\dal}\frac{\eta_{\mu}}{\eta_{\dmu}}e_{\mu} +
\frac{e_{\dmu}\eta_{\dal}}{2\eta_{\dmu}}e_{\mu} \\
&\ = \
\frac{e_{\dmu}\eta_{\dal}}{2\eta_{\dal}}e_{\al} +
\delta_{\dmu\dal}\sum_{\dbe}\frac{e_{\dbe}\eta_{\dal}}{2\eta_{\dbe}}
e_{\be} - \frac{e_{\dal}\eta_{\dmu}}{2\eta_{\dmu}}e_{\mu}.
\endalign
$$
Then the required statement follows. 

\smallskip

From now on we use one symbol $\tnabla$ to denote
 covariant differentials $\tnabla^1$ on
$\cT_1$, $\tnabla^0$ on $\cT_0$ and $\tnabla^1\oplus \tnabla^0$
on $\cT=\cT_1 \oplus \cT_0$. 

\medskip

{\bf 4.8.4 Metrics on $\cT$.} The Egoroff metric $g$ gives rise to 

\smallskip

(i) an even metric $h$ on $\cT_0$, $h(X,Y):= g(\Pi X,\ts(Y))$,
for any $ X, Y \in \cT_0$ (note that $g(\Pi X,\ts(Y))=
g(\Pi X,s(Y))$);

\smallskip
 
(ii) an odd metric $\tilde{g}$ on $\cT$ as the pullback of $g$
relative  to the isomorphism $\ts: \cT\rar \cT\cM$. Note that that $\cT_0$ and
$\cT_1$ are isotropic and $\tilde{g}(X,Y)= g(\ts(X),Y)= g(s(X),Y)$ 
for any $X\in \cT_0$, $Y\in \cT_1$.

\smallskip

Due to the isomorphism $\odot^2 (\cT_0^*)= \odot^2(\Pi \cT_1^*) =\Lambda^2
(\cT_1^*)$, the metric $h$ on $\cT_0$ can also be viewed as an
even non-degenerate skew-symmetric form on $\cT_1$. Explicitly, 
$h$ is given by
$$
h(e_{\al}, e_{\be}) = \delta_{\al\be}\eta_{\dal},  \ \ 
h(e_{\al}, e_{\dbe}) = 0, \ \ \
h(e_{\dal}, e_{\dbe}) = \delta_{\dal\dbe}\eta_{\dal}, 
$$
while $\tilde{g}$ satisfies 
$$
\tilde{g}(e_{\al}, e_{\be}) = 0, \ \ \
\tilde{g}(e_{\al}, e_{\dbe}) = \delta_{\al\be}\eta_{\dal}, \ \ 
\tilde{g}(e_{\dal}, e_{\dbe}) = 0. 
$$
\medskip

{\bf 4.8.5. Frobenius property.} The triple  $(\cT, \circ, \tilde{g})$ 
obviously satisfies 
$$
\tilde{g}(X,Y) = \theta (X\circ Y)
$$
for any $X,Y\in \cT$, where the odd 1-form $\theta$ is defined by
$$
\theta = {\Cal P}\delta \Psi,
$$
$\delta$ being the  SUSY-differential  and $\Cal P$ 
the parity change functor [Ma3]. In particular, one has
$$
\tilde{g}(X\circ Y, Z) =  \tilde{g}(X, Y\circ Z)
$$
for any $X,Y,Z\in \cT$ (which is a defining property of the so-called {\em Frobenius
algebras}\, [D], [H]).  

\smallskip
 
\proclaim{\quad 4.8.6 Proposition} $\tnabla h = 0$, \ $\tnabla \tilde{g}=0$.
\endproclaim

\smallskip

{\bf Proof}. Let us view $h$ as, for example, a skew-form on $\cT_1$. Then
$$
\align
(\tnabla_{e_{\dmu}} h)(e_{\dal}, e_{\dga}) & \ =  e_{\dmu}h(e_{\dal},
e_{\dga}) - h(\tnabla_{e_{\dmu}} e_{\dal}, e_{\dga}) + h(e_{\dal}, 
\tnabla_{e_{\dmu}}e_{\dga}) \\
&\ = \ \delta_{\dal\dga}e_{\dmu}\eta_{\dal} 
-\frac{1}{2}\delta_{\dal\dga}e_{\dmu}\eta_{\dal} 
+ \frac{1}{2}\delta_{\dmu\dga}e_{\dal}\eta_{\dmu} 
- \frac{1}{2}\delta_{\dal\dmu}e_{\dga}\eta_{\dal} \\ 
& \ \ \ \ -\,  \frac{1}{2}\delta_{\dal\dga}e_{\dmu}\eta_{\dal} 
+ \frac{1}{2}\delta_{\dal\dmu}e_{\dga}\eta_{\dal} 
- \frac{1}{2}\delta_{\dmu\dga}e_{\dal}\eta_{\dmu}\\ 
&\ =\ 0. 
\endalign
$$
Analogously one checks other statements. 

\medskip

{\bf 4.9. Odd identity.} 
The $\bullet$-identity $\var$ is said to be {\em flat}\  if $\tnabla \var=0$.

\smallskip

{\bf 4.9.1. Proposition.} \ $\tnabla \var=0\ \Leftrightarrow\
\sum_{\dal}\eta_{\dal}= \text{const}$.

\smallskip

{\bf Proof}.
$$
\align
\tnabla_{e_{\dmu}}\left( \sum_{\al}e_{\dal}\right) & \ = \ \sum_{\dbe}
\frac{e_{\dmu}\eta_{\dbe}}{2\eta_{\dbe}}e_{\dbe} - 
\frac{\sum_{\al} e_{\dal}\eta_{\dmu}}{2 \eta_{\dmu}}e_{\dmu} +
\sum_{\dbe}\frac{e_{\dmu}\eta_{\dbe}}{2\eta_{\dbe}}e_{\dbe} \\
&\ = \ \frac{\eta_{\mu}}{2\eta_{\dmu}} -
\frac{\sum_{\al} e_{\dal}\eta_{\dmu}}{2 \eta_{\dmu}}e_{\dmu}  \\ 
&\ = \ \frac{e_{\dmu}(\sum_{\al}\eta_{\dal})}{2 \eta_{\dmu}}e_{\dmu}  \\ 
&\ =\ 0. \ 
\endalign
$$

\smallskip

\proclaim{\quad 4.9.2. Proposition}  Flatness of $\var$ implies flatness of $e$,
or, equivalently, 
$$
\sum_{\dal}\eta_{\dal}= \text{\rm const} \ \ \Rightarrow\ \ \
\sum_{\al}\eta_{\al}= \text{\rm const}.
$$
\endproclaim


{\bf Proof.}
$$
\align 
\sum_{\dal}\eta_{\dal}= \text{const} & \ \Leftrightarrow \ 
e_{\dbe}\left(\sum_{\dal} \eta_{\dal}\right)= 0 \\
&\ \Leftrightarrow \
\eta_{\be} = \sum_{\dal:\dal\neq \dbe} e_{\dal}\eta_{\dbe} \\
&\ \Rightarrow \ \sum_{\be} \eta_{\be} =  \sum_{\dal,\dbe:\dal\neq \dbe}
e_{\dal}\eta_{\dbe} = 0. 
\endalign
$$

\medskip

{\bf 4.9.3. Orthogonality of flat identities.} 
In the presence of an
Euler field, the flatness of the odd identity $\var$ implies that either 
$D=3/2$ or $g(\bar{e},\var)$ is identically zero, i.e.\ the even and odd
identities $\bar{e}, \var\in (\cT\cM,\circ)$ are everywhere $g$-orthogonal.
Indeed, 
$$
g(e,\var) = g(\sum_{\al}\p_{\al}, \sum_{\dbe}e_{\dbe}) =
\sum_{\dal}\eta_{\dal},
$$
while Proposition 4.9.1 and equation (4.12) imply
$(D-3/2)\sum_{\dal}\eta_{\dal} = 0$.

\smallskip

Since $\tilde{g}(e,\var)=g(s(e),\var)=g(\bar{e}, \var)$, one may
reformulate
the above observation as the $\tilde{g}$-orthogonality of the identities
$(e,\var)\in (\cT,\circ)$.

\medskip

{\bf 4.10. Uniqueness and flatness of $\tnabla$.} The main justification
for introducing the connection $\tnabla$ comes from the following result.

\smallskip

\proclaim{\quad 4.10.1. Proposition}\  Let $\cM$ be a semisimple
pre-Frobenius supermanifold. Then the associated connection
$\tnabla$ is flat if and only if $\Psi$
satisfies the equation (4.3) and
$$
\sum_{\dbe: \dbe\neq \dmu,\dnu} e_{\dbe} \gamma_{\dmu\dnu} = e_{\dmu}
\gamma_{\dmu\dnu} + e_{\dnu}\gamma_{\dmu\dnu} 
\eqno{(4.16)}
$$
for all $\dmu\neq \dnu$.
\endproclaim

\smallskip

{\bf Proof}\ is a straightforward but lengthy calculation.

\smallskip

\proclaim{\quad 4.10.2. Corollary}  Let $\cM$ be a semisimple
pre-Frobenius supermanifold. If the associated connection 
$\tnabla$ is flat, then $\cM$ is semisimple Frobenius.
\endproclaim

\smallskip

{\bf Proof}.
By Proposition~{4.5.3}, it will suffice to show that 
equation~(4.16) implies equation (4.4).
 This is established by the following calculation:
$$
\align
\sum_{\al} \p_{\al} \gamma_{\dmu\dnu} & \ = \
\sum_{\dal,\dbe}e_{\dal}e_{\dbe}\gamma_{\dmu\dnu} \\
&\ = \ 2\sum_{\dal} e_{\dal}e_{\dmu}\gamma_{\dmu\dnu} + 
2\sum_{\dal} e_{\dal}e_{\dnu}\gamma_{\dmu\dnu} \\
&\ = \ 4e_{\mu}\gamma_{\dmu\dnu} - 2\sum_{\dal} e_{\dmu}(e_{\dal}\gamma_{\dmu\dnu})
+ 4e_{\nu}\gamma_{\dmu\dnu} - 2\sum_{\dal} e_{\dnu}(e_{\dal}\gamma_{\dmu\dnu})\\ 
&\ = \ 4e_{\mu}\gamma_{\dmu\dnu} - 2e_{\dmu}(2e_{\dmu}\gamma_{\dmu\dnu} +
2e_{\dnu}\gamma_{\dmu\dnu}) +
4e_{\nu}\gamma_{\dmu\dnu} - 2e_{\dnu}(2e_{\dmu}\gamma_{\dmu\dnu} +
2e_{\dnu}\gamma_{\dmu\dnu}) \\
&\ = \ -\, 4(e_{\dmu}e_{\dnu} + e_{\dnu}e_{\dmu})\gamma_{\dmu\dnu}\\
&\ = \ 0. 
\endalign
$$

\smallskip

\proclaim{\quad 4.10.3. Corollary}  Let $\cM$ be a semisimple pre-Frobenius
supermanifold. If $\var$ is flat and $\Psi$ satisfies the equation
(4.3), then $\tnabla$ is flat as well (implying that $\cM$
is Frobenius).
\endproclaim

\smallskip

{\bf Proof}. Assume $e_{\dmu} (\sum_{\dbe}\eta_{\dbe})=0$, or,
equivalently, 
$$
\eta_{\mu} = \sum_{\dbe:\dbe\neq \mu} e_{\dbe}\eta_{\dmu}.
$$
Then
$$
\align
\sum_{\dbe: \dbe\neq \dmu,\dnu} e_{\dbe} \gamma_{\dmu\dnu}
&\ = \ \sum_{\dbe:\dbe\neq \dmu,\dnu} \frac{e_{\dbe}e_{\dmu}\eta_{\dnu}}
{2\sqrt{\eta_{\dmu}\eta_{\dnu}}} - 
\sum_{\dbe:\dbe\neq \dmu,\dnu} \frac{(e_{\dbe}\eta_{\dmu})(e_{\dmu}\eta_{\dnu})}
{4\eta_{\dmu}\sqrt{\eta_{\dmu}\eta_{\dnu}}} - 
\sum_{\dbe:\dbe\neq \dnu} \frac{(e_{\dbe}\eta_{\dnu})(e_{\dmu}\eta_{\dnu})}
{4\eta_{\dnu}\sqrt{\eta_{\dmu}\eta_{\dnu}}} \\
&\ =\ - \sum_{\dbe:\dbe\neq \dnu} \frac{e_{\dmu}(e_{\dbe}\eta_{\dnu})}
{2\sqrt{\eta_{\dmu}\eta_{\dnu}}} + \frac{e_{\mu}\eta_{\dnu}}
{2\sqrt{\eta_{\dmu}\eta_{\dnu}}} - \frac{\eta_{\mu}(e_{\dmu}\eta_{\dnu})}
{4\eta_{\dmu}\sqrt{\eta_{\dmu}\eta_{\dnu}}} - 
\frac{\eta_{\nu}(e_{\dmu}\eta_{\dnu})}
{4\eta_{\dnu}\sqrt{\eta_{\dmu}\eta_{\dnu}}} \\ 
&\ =\ - \, \frac{e_{\nu}\eta_{\dmu}}
{2\sqrt{\eta_{\dmu}\eta_{\dnu}}} +
\frac{e_{\mu}\eta_{\dnu}}
{2\sqrt{\eta_{\dmu}\eta_{\dnu}}} 
- \frac{\eta_{\mu}(e_{\dmu}\eta_{\dnu})}
{4\eta_{\dmu}\sqrt{\eta_{\dmu}\eta_{\dnu}}} - 
\frac{\eta_{\nu}(e_{\dmu}\eta_{\dnu})}
{4\eta_{\dnu}\sqrt{\eta_{\dmu}\eta_{\dnu}}} \\ 
&\ =\ e_{\dmu}\gamma_{\dmu\dnu} + e_{\dnu}\gamma_{\dmu\dnu}, 
\endalign
$$
so that the equation (4.16) is satisfied.

\smallskip

For later use we give the following characterization of $\tnabla$:

\smallskip

\proclaim{\quad 4.10.4. Proposition} Let $\cM$ be a semisimple pre-Frobenius 
supermanifold. A linear connection $\nabla: \cT_1 \rar \cT_1^*\ot \cT_1$ 
satisfies the conditions

\smallskip

(a) $\nabla h =0$,

\smallskip

(b) $\nabla_{e_{\dmu}} e_{\dal} + \nabla_{e_{\dal}} e_{\dmu} =
2\delta_{\dmu\dal} \nabla_{e_{\dmu}} e_{\dmu}$,

\smallskip

(c) $h(\nabla_{e_{\dmu}} e_{\dal}, e_{\dbe})=0$ for any
$\dmu\neq \dal\neq \dbe\neq \dmu$

\smallskip

if and only if it is given by (4.15).
\endproclaim

\smallskip

{\bf Proof.}  Define $\Delta_{\dmu\dal\dbe}$ as
$$
\nabla_{e_{\dmu}} e_{\dal} = \sum_{\dbe}
\frac{\Delta_{\dmu\dal\dbe}}{2\eta_{\dbe}} e_{\dbe}.
$$
Then
$$
\align
(a) \ & \Leftrightarrow \ e_{\dmu}\eta_{\dal}\delta_{\dal\dbe} =
\frac{1}{2}(\Delta_{\dmu\dal\dbe} + \Delta_{\dmu\dbe\dal}),\\
(b)\ & \Leftrightarrow \ \Delta_{\dmu\dal\dbe} +
\Delta_{\dal\dmu\dbe} = 2 \delta_{\dal\dmu} \Delta_{\dmu\dmu\dbe},\\
(c) \ & \Leftrightarrow \ \Delta_{\dmu\dal\dbe} = 0 \ \text{for all}\ 
\dmu\neq \dal\neq \dbe\neq \dmu,
\endalign
$$
implying that the only non-vanishing components of
$\Delta_{\dmu\dal\dbe}$
are 
$$
\Delta_{\dmu\dal\dal}= \Delta_{\dmu\dal\dmu} = \Delta_{\dal\dal\dmu} = 
e_{\dmu}\eta_{\dal}. 
$$
Hence 
$$
\nabla_{e_{\dmu}} e_{\dal} = \left\{ \matrix
\sum_{\dbe}\frac{e_{\dbe}\eta_{\dal}}{2\eta_{\dbe}} e_{\dbe} &
\text{for}\ \dmu=\dal, \\
\frac{e_{\dmu}\eta_{\dal}}{2\eta_{\dal}} e_{\dal} + 
\frac{e_{\dmu}\eta_{\dal}}{2\eta_{\dmu}} e_{\dmu}=  
\frac{e_{\dmu}\eta_{\dal}}{2\eta_{\dal}} e_{\dal} - 
\frac{e_{\dal}\eta_{\dmu}}{2\eta_{\dmu}} e_{\dmu} & \text{for}\
\dmu\neq \dal,
\endmatrix \right.
$$
implying $\nabla=\tnabla$.

\newpage

\centerline{\bf \S 5. Supersymmetric Schlesinger equations}

\bigskip

{\bf 5.1. Meromorphic connections with logarithmic singularities.}  
Let $\cM$ be a complex supermanifold equipped with a SUSY structure
(4.1) and let $\cF\rar \cM$ be a locally free holomorphic 
sheaf on $\cM$. Keeping in mind semisimple Frobenius structures, we also
assume that $\cM$ comes equipped with a monomorphism $\Theta: \cT_0 \rar
\Lambda^2 \cT_1$ which, as explained in 4.4, canonically extends any
covariant differential $\nabla: \cF \rar \cT_1^*\ot V$ to a linear
connection on $\cF$.

\smallskip

Let $D$ be a complex super-submanifold of $\cM$ of codimension
$1|0$.

\smallskip

Assume first that $D$ is irreducible and that the associated 
divisor line bundle $[D]$ is free. Let $f$ be a global  basis section of
$[D]$. A holomorphic covariant differential
$$
\nabla: \cF \rar \cT_1^*\ot \cF
$$
on $\cM\setminus D$ is said to be a {\em meromorphic connection
with logarithmic singularities}\ along $D$ if there is a holomorphic 
covariant differential
$$
\nabla': \cF \rar \cT_1^*\ot \cF
$$
on $\cM$ such that
$$
\nabla - \nabla' =  A\frac{\delta f}{f}
$$
for some even holomorphic section $A\in H^0(\cM, \cF\ot \cF^*)$. 

\smallskip

Note that 

\medskip

a) this  definition
does not depend on the choice of a particular trivialization $f$ of
$[D]$ and hence can be appropriately localized and generalized;

\medskip

b) the section $A$ restricted to $D$ does not depend on the choices
made and hence gives a well-defined element of 
$H^0(D, \cF\ot \cF^*)$ which is called
the {\it residue}\ of $\nabla$ at $D$;

\medskip

c) for any local trivialization $f$ of $[D]$, the connection $\nabla$
induces a holomorphic {\em residual connection}\ 
$\nabla^{D,f}$ on $\left.\cF \right|_D$ associated with the
holomorphic covariant differential  $\left.\nabla'\right|_D=\left.(\nabla -
A\frac{\delta f}{f})\right|_D$; if $\nabla$ is flat, $\nabla^{D,f}$ is
flat for any $f$.

\medskip

{\bf 5.2. Universal isomonodromic deformation.}
 Consider $\C^{n|n}$
with its natural SUSY structure $\cT_1 \subset \cT\C^{n|n}$
spanned by  the vector fields 
$e_{\dal}= \p_{\dal} + \theta^{\dal}\p_{\al}$, where 
$(u^{\al}, \theta^{\dal})$ are natural coordinates.
Let $B$ be the universal covering of $\C^{n|n}\setminus (u^{\al}- u^{\be}
- \theta^{\dal}\theta^{\dbe}=0)$, where pairwise distinct integers $\al$ and $\be$
run over $1,\ldots,n$.

Consider a supermanifold $B\times \CP^{1|1}$ with the direct product
SUSY-structure, and  
denote by $\tilde{D}_{\al}$ the inverse image in $B\times  \CP^{1|1}$ 
of the submanifold $\la - u^{\al} - \xi \theta^{\dal} = 0$ in
$\C^{n|n}\times \CP^{1|1}$, where
 $(\la,\xi)$ are natural coordinates on a big cell of $\CP^{1|1}$. 
Furthermore, define $\tilde{D}_{\infty} = B \times \infty \subset B\times
\CP^{1|1}$, where $\infty$ stands for the codimension $1|0$ submanifold of
$\CP^{1|1}$ given by $\hat{\la}=0$, where $\hat{\la}= 1/\la$.

To any given
point $(x_0^{\al}, \theta_0^{\al})\in \C^{n|n}$ one may associate 
an embedding 
$$
\matrix
i_0: & \CP^{1|1} & \lon & B\times \CP^{1|1} \\
     &  (\la,\xi) & \lon & (x_0^{\al}, \theta_0^{\al}, \la, \xi).
\endmatrix
$$

\proclaim{\quad Theorem 5.2.1} 
Let $\cF^0$ be a locally free sheaf of rank $p|q$ on 
$\CP^{1|1}$ and let $\nabla^0$ be 
a flat meromorphic
connection on it with logarithmic singularities at $\cup_{\al=1}^n D^0_{\al}\cup
\infty$, where $D^0_{\al}\subset \CP^{1|1}$ is given by 
$\la - u_0^{\al} - \xi \theta_0^{\dal} =0$.
\smallskip

 Then there exists a locally free sheaf $\cF$ of rank $p|q$
on $B\times \CP^{1|1}$ and a flat meromorphic connection $\nabla$ on it
such that

\medskip

(a) $\nabla$ has logarithmic singularities at $\tilde{D}_{\al}$,
$\al=1,\dots, n$, and $\tilde{D}_{\infty}$;

\medskip

(b) there is a canonical isomorphism $i: i_0^*(\cF, \nabla) \rar 
(\cF^0, \nabla^0)$;

\medskip

 (c) the data $(\cF, \nabla, i)$ are unique up to unique
isomorphism.

\endproclaim

\smallskip

{\bf Comment on the proof.} According to Penkov [P], a pair
$(\cE, \nabla)$ consisting of a locally free sheaf $\cE$ on a
supermanifold $\cM$ and a flat holomorphic connection $\nabla$ on
$\cE$ is uniquely determined by the associated monodromy representation 
of $\pi_1(\cM_{\text{red}})$ on $\cE_{\text{red}}$. 
This together with the observation made in 2.2.2 about the isomorphism
of the first homotopy groups of the underlying classical manifolds, 
immediately implies that there is a pair $(\cF, \nabla)$ on 
$B\times \CP^{1|1}\setminus (\cup_{\al=1}^n \tilde{D}_{\al}\cup
\tilde{D}_{\infty})$ such that the statements (b) and (c) are true outside
singularities.

\smallskip

Using a straightforward generalization
of the original Malgrange's arguments, one may extend
$(\cF, \nabla)$ to $B\times \CP^{1|1}$ in such
a way that (a)-(c) hold. 

\medskip

{\bf 5.3. Supersymmetric Schlesinger equations.} In this subsection we
will assume that $\cF^0= T\ot \f_{\CP^{1|1}}$ where $T$ is a 
vector superspace of dimension $p|q$.

\smallskip

Using the semicontinuity principle as in Section 2.3, one may show that
there is an open subset $B'\subset B$ such that $\cF$ is free on
$B'\times \CP^{1|1}$.  Moreover, one may identify $\cF$ on $B'\times
\CP^{1|1}$ with $T\ot \f_{B'\times \CP^{1|1}}$ 
compatibly with the respective trivialization of $\cF^0$. Indeed,
one may first trivialize $\cF$ along $\tilde{D}_{\infty}$ using the residual 
connection, and then  take the
constant extension of each horizontal section along $\CP^{1|1}$. Since
$$
\delta(u^{\nu}- \la - \theta^{\dnu}\xi) = (\theta^{\dnu} - \xi)d(\theta^{\dnu}-\xi),
$$
in the chosen trivialization $\left.\cF\right|_{B'\times
\CP^{1|1}}= T\ot \f_{B'\times \CP^{1|1}}$ 
the covariant differential $\nabla$
must be of the form 
$$
\nabla = \delta + \sum_{\nu=1}^n \frac{A_{\nu}(\theta^{\dnu} - \xi)}{u^{\nu}
-\la - \theta^{\dnu}\xi} d(\theta^{\dnu}-\xi) \eqno{(5.1)}
$$
for some even meromorphic sections $A_{\nu}\in H^0(B, \cF\ot \cF^*)$.

\smallskip

\proclaim{\quad Theorem 5.3.1} The connection (5.1) is flat if and only
if
$$
e_{\dmu} A_{\nu} = -\, \frac{\theta^{\dmu} -
\theta^{\dnu}}{u^{\mu}-u^{\nu} - \theta^{\dmu}\theta^{\dnu}}
 [A_{\mu},A_{\nu}]  \eqno{(5.2)}
$$
$$
e_{\dmu}A_{\mu} =  \sum_{\nu: \nu\neq \mu} 
\frac{\theta^{\dmu} - \theta^{\dnu}}{u^{\mu}-u^{\nu} -
\theta^{\dmu} \theta^{\dnu}}  [A_{\mu},A_{\nu}] . \eqno{(5.3)}
$$
or, equivalently,
$$
dA_{\mu} =  \sum_{\nu: \nu\neq \mu} 
\frac{d(u^{\mu}-u^{\nu}-\theta^{\dmu}\theta^{\dnu})}{u^{\mu}-u^{\nu} -
\theta^{\dmu} \theta^{\dnu}} [A_{\mu}, A_{\nu}]  \eqno{(5.4)}
$$
where $d$ is the usual exterior differential on $B$.
\endproclaim

\smallskip
{\bf Proof} \ is straightforward.

\smallskip

The equations (5.2) and (5.3) are called {\em supersymmetric Schlesinger's
equations}.

\medskip

{\bf 5.4.  From Frobenius supermanifolds to strict special solutions of 
 Schlesinger's equations}. Let $\cM$ be a semisimple Frobenius
supermanifold with an Euler field $E$ and flat identities
$\var$ and $\bvar$. Define an even linear operator
$\cV: \cT_1 \rar \cT_1$ as follows, 
$$
\cV(X) = p_1(\nabla_X E) - \frac{1}{2}(D-\frac{1}{2}) X \ \ \ \mbox{for any}
\ X\in \cT_1,
$$
where $\nabla$ is the Levi-Civita connection and $p_1: \cT\cM\rar \cT_1$
the projection defined in 4.8.2.

\smallskip

\proclaim{\quad Theorem 5.4.1} 
a). Let $f_{\dal}= e_{\dal}/\sqrt{\eta_{\dal}}$, then
$$
\cV(f_{\dal}) = \sum_{\dbe:\dbe\neq
\dal}\left[\theta^{\dal}\ga_{\dbe\dal} + u^{\dbe}\p_{\dbe}\ga_{\dbe\dal}
-u^{\dal}\p_{\dal}\ga_{\dal\dbe} + \sum_{\ga:\ga\neq\al,\be}
u^{\ga}e_{\dga} \ga_{\dal\dbe}\right] f_{\dbe}. 
\eqno{(5.5)}
$$

\medskip

 b). $\cV$ is symmetric relative to $h$, i.e.
$$
h(\cV(X), Y) + h(X,\cV(Y))=0 \ \ \ \mbox{for any}
\ X,Y\in \cT_1,
$$

\medskip

c). $\tnabla \cV=0$.

\medskip

d). $\cV(\var)=  \frac{3-2D}{4} \var$.

\endproclaim

\smallskip

{\bf Proof}. Items (a) and (c) follow from a straightforward calculation which we
omit, while (b) immediately follows from (4.9) and (5.5).
To check (c), note that
$$
\align
\cV(\var) &\  =\  p_1(\nabla_{\var} E) - \frac{1}{2}(D-\frac{1}{2})\var \\
&\ = \ p_1(\nabla_E \var + [\var,E]) - \frac{1}{2}(D-\frac{1}{2})\var \\ 
&\ = \ \frac{1}{2}\var - \frac{1}{2}(D-\frac{1}{2})\var \\ 
&\ = \  \frac{3-2D}{4} \var . 
\endalign
$$

\medskip

{\bf 5.4.2. The structure connection.} Let $\cM$ be a semisimple 
pre-Frobenius supermanifold with canonical coordinates $(u^{\al},
\theta^{\dal})$,  and let $\check{\cM}$ be a 
$(n+1|n+1)$-supermanifold  $\cM \times
\CP^{1|1}\setminus {\cup_{\nu}(u^{\nu}-\la -\theta^{\dnu}\xi = 0)}$
equipped with a product SUSY structure 
$$
\check{\cT}_1:=\mbox{span}(e_{\dal}, e_{\xi}) \subset \cT\check{\cM},
$$
where $e_{\xi}= \p/\p\xi + \xi\p/\p \la$.  Let $\mbox{pr}: \check{\cM}
\rar \cM$ be a natural projection  and denote
$E_0:={\sum_\al}u^{\al}e_{\dal}$ and $E_1:= {\sum_\al}\theta^{\dal}e_{\dal}$. 

\smallskip

The following theorem
introduces  a so-called {\em structure connection}\, on
$\mbox{pr}^{*}(\cT_1)\subset \check{\cT}_1 $ which associates with any 
semisimple Frobenius structure
a 1-parameter solution of Schlesinger's equations.

\smallskip

\proclaim{\quad 5.4.3. Theorem}  For any $X,Y \in \mbox{\rm
pr}^{-1}(\cT_1)\subset \check{\cT}_1$ put
$$
\check{\nabla}_X Y = \tnabla_X Y - (\cV + \kappa \text{\rm Id})
(E_0 - \la + \xi E_1)^{-1} \bullet 
(E_1 - \xi)\bullet X \bullet Y 
\eqno{(5.6)}
$$
$$
\check{\nabla}_{e_{\xi}} Y =  (\cV + \kappa \text{\rm Id})(E_0 - \la + \xi
E_1)^{-1}\bullet (E_1 - \xi)\bullet Y.
\eqno{(5.7)}
$$
where $\kappa \in \C$.  

\smallskip

If $\cM$ is a semisimple Frobenius supermanifold with an Euler field and
flat identities $\var $ and $\tvar$, then $\check{\nabla}$ is a flat
connection on  $\text{\rm pr}^{*}(\cT_1)\subset \check{\cT}_1$ for any
$\kappa$.

\endproclaim

{\bf Proof.}  If $X=e_{\dmu}$ and $Y=e_{\dal}$, the equations (5.6) and
(5.7) take the form
$$
\check{\nabla}_{e_{\dmu}}e_{\dal} = \tnabla_{e_{\dmu}}e_{\dal} +
\delta_{\dal\dmu} \frac{[\cV(e_{\dal}) + \kappa e_{\dal}](\theta^{\dal}
- \xi)}{u^{\al} - \la - \theta^{\dal}\xi} 
$$
$$
\check{\nabla}_{e_{\xi}}e_{\dal} = -\,
\frac{[\cV(e_{\dal}) + \kappa e_{\dal}](\theta^{\dal}
- \xi)}{u^{\al} - \la - \theta^{\dal}\xi} .
$$

Then, by Proposition 5.3.1, it will suffice to show that under
the conditions stated in 5.4.2 the matrix valued fields
$$
A_{\nu\, \dal}^{\ \, \dbe} = - \delta_{\dnu\al}[\cV_{\dal\dbe} +
\kappa \delta_{\dal\dbe}]
$$
satisfy, for any $\kappa\in \C$, the Schlesinger's equations
$$
e_{\dmu}A_{\nu\,\dal}^{\ \, \dbe} + \Gamma_{\dmu\dde}^{\dbe}A_{\nu\, \dal}^{\ \,
\dde} - \Gamma_{\dmu\dal}^{\dde}A_{\nu\, \dde}^{\ \, \dbe}  
= -\, \frac{\theta^{\dmu} -
\theta^{\dnu}}{u^{\mu}-u^{\nu} - \theta^{\dmu}\theta^{\dnu}}
 [A_{\mu},A_{\nu}]_{\dbe}^{\dal}, \ \ \ \mu\neq \nu  
\eqno{(5.8)}
$$
$$
e_{\dmu}A_{\dmu\,\dal}^{\ \, \dbe} + \Gamma_{\dmu\dde}^{\dbe}A_{\dmu\, \dal}^{\ \,
\dde} - \Gamma_{\dmu\dal}^{\dde}A_{\nu\, \dde}^{\ \, \dbe}  
=  \sum_{\nu: \nu\neq \mu} 
\frac{\theta^{\dmu} - \theta^{\dnu}}{u^{\mu}-u^{\nu} -
\theta^{\dmu} \theta^{\dnu}}  [A_{\mu},A_{\mu}]_{\dal}^{\dbe}, \
\eqno{(5.9)}
$$
where 
$$
\Gamma_{\dmu\dal}^{\dbe} =
\delta_{\dal\dbe}\frac{e_{\dmu}\eta_{\dal}}{2\eta_{\dal}} 
- \delta_{\dbe\dmu}\frac{e_{\dal}\eta_{\dmu}}{2\eta_{\dmu}} 
+\delta_{\dmu\dal}\frac{e_{\dbe}\eta_{\dal}}{2\eta_{\dbe}}
$$ 
are the coefficients of the connection $\tnabla$, 
and $\cV_{\dal\dbe}$ are the coefficients of
the operator $\cV$ in the basis $e_{\dal}$.

 It is not hard to show that equations (4.9) and (5.5)
imply, for $\mu\neq \nu$,
$$
\frac{\theta^{\dmu} - \theta^{\dnu}}{u^{\mu}-u^{\nu} -
\theta^{\dmu} \theta^{\dnu}} \cV_{\dmu\dnu} = 
\frac{e_{\dmu}\eta_{\dnu}}{2\eta_{\dnu}} 
$$
which in turn imply that the terms of order $\kappa^0$ in (5.8)
and (5.9) are all equivalent to the equation $\tnabla \cV=0$ which
follows from 5.4.1(c).

Since there are no terms quadratic in $\kappa$ in (5.8)
and (5.9), it remains to show that the terms of order $\kappa^1$
cancel. Indeed, the terms linear in $\kappa$ in the l.h.s.\ of
(5.8) are
$$
-\delta_{\nu\al}\Gamma_{\dmu\dal}^{\dbe}
+\delta_{\nu\be}\Gamma_{\dmu\dal}^{\dbe}  =  
(\delta_{\nu\al}\delta_{\mu\be} - \delta_{\mu\al}\delta_{\nu\be})
\frac{e_{\dal}\eta_{\dbe}}{2\eta_{\dbe}},
$$
while the terms linear in $\kappa$ in the r.h.s.\ of (5.8) are
$$
-\frac{\theta^{\dmu} - \theta^{\dnu}}{u^{\mu}-u^{\nu} -
\theta^{\dmu} \theta^{\dnu}}\left(-\delta_{\nu\al}\delta_{\mu\be}\cV_{\dnu\dmu} +
\delta_{\mu\al}\delta_{\nu\be}\cV_{\dmu\dnu}\right)=
(\delta_{\nu\al}\delta_{\mu\be} - \delta_{\mu\al}\delta_{\nu\be})
\frac{e_{\dal}\eta_{\dbe}}{2\eta_{\dbe}}.
$$
Analogously one checks that the terms linear in $\kappa$ in (5.9)
also vanish identically. 

\medskip


{\bf 5.5. Strict special solutions of Schlesinger's equations.}
Consider the supermanifold $\cM=\C^{n|n}$ with canonical coordinates
$(u^{\al}, \theta^{\dal})$, an  $(p|q)$-dimensional vector space $T$ and
a set of even holomorphic matrix functions $A_{\nu}: \cM \rar \mbox{End}(T)$,
$\nu=1, \ldots,n$ such that the Schlesinger equations (5.4)
are satisfied. In particular, summing (5.4) over $\nu$ and
denoting $\cW= \sum_{\nu} A_{\nu}$, we find $d\cW=0$, i.e.\ $\cW\in
\mbox{End}(T)$. 

\smallskip

Theorem~5.4.3 motivates the following definition (cf.\ Section 2.4):

\smallskip

\proclaim{\quad 5.5.1. Definition}  A solution to Schlesinger's equation as above
is called\ {\rm strict special} if

\medskip

a). $\mbox{\rm rank}\ T = 0|n$;

\medskip

b). $T$ is endowed with a complex non-degenerate skew-symmetric form $h\in
\Lambda^2(T)$;

\medskip

c). $\cW= - \cV - \kappa \mbox{\rm Id}$, where $\kappa\in \C$ is
an arbitrary parameter and $\cV\in \mbox{\rm End}(T)$ is an even
operator symmetric with respect to $h$;

\medskip

d). for any $\nu$,
$$
 A_{\nu} = - (\cV + \kappa \mbox{\rm Id})P_{\nu},
\eqno{(5.10)}
$$
where $P_{\nu}: \cM\rar \mbox{\rm End}(T)$ is a set of even holomorphic
matrix functions whose values at any point of $\cM$ constitute a
complete system of  orthogonal projectors of rank 1 with respect to $h$:
$$
P_{\nu}P_{\mu}=\delta_{\mu\nu} P_{\nu}, \ \ \ \sum_{\nu}P_{\nu}= \mbox{\rm
Id}_T, \ \ \ h(\mbox{\rm Im}P_{\nu}  , \mbox{\rm Im}P_{\mu}) = 0
\eqno{(5.11)}
$$
if $\mu\neq \nu$.

\medskip

e). There is a vector $\var \in T$ such that
$$
\cV(\var) = \frac{3-2D}{4}\var
\eqno{(5.12)}
$$
for some $D\in \C$ 
and $e_{\dnu}:= P_{\nu}(\var)$ are nowhere vanishing on $\cM$.
\endproclaim

\smallskip

Theorem 5.4.3 says that to any semisimple Frobenius supermanifold
there is canonically associated a strict special solution of Schlesinger's
equations. 

\medskip


{\bf 5.6. From strict special solutions to Frobenius supermanifolds.} 
Let $(\cM, T, h, A_{\nu}, \var)$ be a strict special
solution of Schlesinger's equations.

\smallskip

\proclaim{\quad 5.6.1. Theorem} These data come from the unique structure
of semisimple split Frobenius supermanifold on $\cM$, with Euler field
and flat identities $\var$ and $\bvar$.
\endproclaim

{\bf Proof}.  Put $e_{\dnu}= P_{\nu}(\var) \in T\ot \f_{\cM}$ and identify 
$\f_{\cM}\ot T$ with $\cT_1\subset \cT\cM$ by setting 
$e_{\dnu}= \p/\p \theta^{\dnu} + \theta^{\dnu} \p/\p u^{\dnu}$. This transfers
$h$ from $T$ to $\cT_1$. Denote $\eta_{\dal}=h(e_{\dal}, e_{\dal})$. 
Define the multiplication in  
$\cT\cM$ by the formulae 4.3.2. 

To prove Theorem 5.6.1 it will suffice to show that
1) $\eta_{\dal}=e_{\dal}\Psi$ for an odd function $\Psi$ (potentiality);
2) $\Psi$ satisfies the Darboux-Egoroff equations (4.3) and
(4.4) (flatness of the Egoroff metric); 3) the equation  (4.12) is
satisfied (Euler property);
4) the equations (4.9) and (4.10) are satisfied
(flatness of $\bvar$). That the condition 
$\sum_{\dal} \eta_{\dal}$=const (flatness of $\var$) is satisfied 
follows immediately from the definition of $\eta_{\dal}$ and the 
 fact that $\sum_{\dal} \eta_{\dal}= 
h(\var,\var)$.

\smallskip

{\bf Step 1 (potentiality)}. Since 
$$
\sum_{\dal}e_{\dal} = \sum_{\al} P_{\al}(\var) = \Id(\var) = \var
$$
and $h(\mbox{Im}P_{\al}  , \mbox{Im}P_{\be}) = 0$ if $\al\neq \be$, we  have 
$h_{\dal}=h(\var,e_{\dal})$ and hence, in view of 5.5.1(c,d),
$$
h(\var, A_{\nu}(\var)) = - h(\var, (\cV+\kappa \Id)e_{\dnu}) = 
h(\cV(\var), e_{\dnu}) - \kappa h(\var,e_{\dnu}) =
\frac{3-2D-4\kappa}{4}\eta_{\dnu}. \eqno{(5.13)}
$$

Let $\nabla$ be the unique flat connection in $\cT_1\simeq T\times \cM$
which makes constants sections of $T\times \cM$ horizontal. Obviously, 
$\tnabla$ preserves $h$ and satisfies $\nabla \var =0$. Then derivating
(5.13) we find, for every $\mu\neq \nu$,
$$
\frac{3-2D-4\kappa}{4}e_{\dmu}\eta_{\dnu} =  h(\nabla_{e_{\dmu}}\var,
A_{\nu} \var) - h(\var,\nabla_{e_{\dmu}}(A_{\nu}\var)) 
=  h(\var, ({e_{\dmu}}A_{\nu})\var).  \eqno{(5.14)}
$$
Since, in view of (5.2),
$$
{e_{\dmu}} A_{\nu} = -\, \frac{\theta^{\dmu} -
\theta^{\dnu}}{u^{\mu}-u^{\nu} - \theta^{\dmu}\theta^{\dnu}}
 [A_{\mu},A_{\nu}]  = -\,{e_{\dnu}} A_{\dmu} \ \ \ \forall\ \mu\neq \nu,
$$
we find
$$
e_{\dmu}\eta_{\dnu} + e_{\dnu}\eta_{\dmu} = 0 \ \ \ \ \ \ \forall\ \mu\neq \nu,
$$
or, equivalently,
$$
e_{\dmu}\eta_{\dnu} + e_{\dnu}\eta_{\dmu} = 2\delta_{\dmu\dnu}\eta_{\mu}
 \ \ \ \ \ \ \forall\ \mu,\nu,
$$
where $\eta_{\mu}:= e_{\dmu}\eta_{\dmu}$. Analogous calculations,
involving Schlesinger's equations (5.2) and (5.3), show
that
$$
e_{\dmu}\eta_{\nu} - \p_{\nu}\eta_{\dmu}=0, \ \ \ \ \p_{\mu}\eta_{\nu}-
\p_{\nu}\eta_{\mu}=0.
$$
Finally, defining the 1-form $\om=\sum_{\al}[d\theta^{\dal}(\eta_{\dal}
- \theta^{\dal}\eta_{\al}) + du^{\al}\eta_{\al}]$,
 it is straightforward to check that the latter three equations 
are equivalent to $d\om=0$. Hence $\om = d\Psi$
for some odd function $\Psi$, i.e.\  $\eta_{\dal}= e_{\dal}\Psi$.

\medskip

{\bf Step 2 (flatness of the Egoroff metric).} Let $\Psi$ be as above
and $g$ the associated Egoroff metric. Let us prove that $g$ is flat.

\smallskip

 By Corollary~4.10.2,
it will suffice to show that the flat connection $\nabla$ coincides
with the connection (4.15), i.e.\ that $\nabla$ satisfies all
three conditions of Proposition~4.10.4. 

\smallskip

The condition 4.10.4(a) is
obviously satisfied. 

\smallskip

Since $\kappa$ is arbitrary, we
may assume without loss of generality that $\cW$ is invertible and
rewrite Schlesinger's equation (5.2) in the form
$$
{e_{\dmu}} P_{\nu} = -\, \frac{\theta^{\dmu} -
\theta^{\dnu}}{u^{\mu}-u^{\nu} - \theta^{\dmu}\theta^{\dnu}}
\cW^{-1} [A_{\mu},A_{\nu}]. 
$$
Then for every $\dmu\neq \dnu$ we have
$$
\align
\nabla_{e_{\dmu}}e_{\dnu} + \nabla_{e_{\dnu}}e_{\dmu} &\ = \ 
\nabla_{e_{\dmu}}(P_{\nu}\var) + \nabla_{e_{\dnu}}(P_{\mu}\var) \\
&\ =\ (e_{\dmu}P_{\nu})\var + (e_{\dnu}P_{\mu})\var \\
&\ = \ -\, \frac{\theta^{\dmu} -
\theta^{\dnu}}{u^{\mu}-u^{\nu} - \theta^{\dmu}\theta^{\dnu}}
\cW^{-1} [A_{\mu},A_{\nu}] \var
- \frac{\theta^{\dnu} -
\theta^{\dmu}}{u^{\nu}-u^{\mu} - \theta^{\dnu}\theta^{\dmu}}
\cW^{-1} [A_{\nu},A_{\mu}] \var \\
&\ = \ 0, 
\endalign
$$
or, equivalently,
$$
\nabla_{e_{\dmu}} e_{\dnu} + \nabla_{e_{\dnu}} e_{\dmu} =
2\delta_{\dmu\dnu} \nabla_{e_{\dmu}} e_{\dmu} \ \ \ \ \ \ \ \forall\ \mu,\nu.
$$

Thus, it remains to check the  condition 4.10.4(c) for all
$\dmu\neq\dnu\neq\dal\neq\dmu$ :
$$
\align
h(\nabla_{e_{\dmu}}e_{\dnu},e_{\dal}) &\ = \  
h((e_{\dmu}P_{\nu})\var, P_{\al}\var)   \\
&\ =\ -\, \frac{\theta^{\dmu} -
\theta^{\dnu}}{u^{\mu}-u^{\nu} - \theta^{\dmu}\theta^{\dnu}}
h(\cW^{-1} [A_{\mu},A_{\nu}]\var,  P_{\al}\var)   \\
&\ = \
-\, \frac{\theta^{\dmu} -
\theta^{\dnu}}{u^{\mu}-u^{\nu} - \theta^{\dmu}\theta^{\dnu}}
h(P_{\mu}(\ldots) + P_{\nu}(\ldots), P_{\al}\var)   \\
&\ = \ 0.
\endalign
$$
This establishes the flatness of the Egoroff metric.

\smallskip

{\bf Step 3 (Euler property).} Schlesinger's equations (5.2)-(5.3) imply
$$
\align
E A_{\nu}&\ = \ \sum_{\mu:\mu\neq \nu}\left(u^{\mu}\p_{\mu}A_{\nu} +
\frac{1}{2}\theta^{\dmu}e_{\dmu}A_{\nu}\right) + u^{\nu}\p_{\nu}A_{\nu} +
\frac{1}{2}\theta^{\dnu}e_{\dnu}A_{\dnu} \\
&\ = \ \sum_{\mu:\mu\neq \nu}\left(-
\frac{u^{\mu}[A_{\mu},A_{\nu}]}{u^{\mu}-u^{\nu} - \theta^{\dmu}\theta^{\dnu}}
- \frac{\theta^{\dmu}(\theta^{\dmu}-
 \theta^{\dnu})[A_{\mu},A_{\nu}]}{2(u^{\mu}-u^{\nu} -
\theta^{\dmu}\theta^{\dnu})}\right. \\
&\ \ \ \ \left.
+ \frac{u^{\nu}[A_{\nu},A_{\nu}]}{u^{\nu}-u^{\mu} - \theta^{\dnu}\theta^{\dmu}}
+ \frac{\theta^{\dnu}(\theta^{\dnu} -\theta^{\dmu})[A_{\nu},A_{\mu}]}
{2(u^{\nu}-u^{\mu} - \theta^{\dnu}\theta^{\dmu})}\right) \\
&\ = \ -\sum_{\mu:\mu\neq \nu}[A_{\mu},A_{\nu}].
\endalign
$$

Then, using (5.13), we find 
$$
\align 
\frac{3-2D-4\kappa}{4}E\eta_{\dnu} &\ = \ h(\var, (EA_{\nu})\var) \\
&\ = \ - h(\var, \sum_{\mu:\mu\neq \nu}[A_{\mu},A_{\nu}]\var) \\
&\ = \ - h(\var, [\cV + \kappa\Id, (\cV+\kappa\Id)P_{\nu}]\var)\\
&\ = \ -\,
\left(\frac{3-2D}{2}\right)\left(\frac{3-2D-4\kappa}{4}\right)\eta_{\dnu}.
\endalign
$$
Hence $E\eta_{\dnu}= (D-\frac{3}{2})\eta_{\dnu}$. 

\smallskip

{\bf Step 4 (flatness of $\bvar$).} One finds from (5.14):
$$
\align
\frac{3-2D-4\kappa}{4}(\theta^{\dmu}-\theta^{\dnu})e_{\dmu}\eta_{\dnu}
&\ = \ (\theta^{\dmu}-\theta^{\dnu}) h(\var, ({e_{\dmu}}A_{\nu})\var)\\
&\ = \ -(\theta^{\dmu}-\theta^{\dnu}) h(\var, \frac{\theta^{\dmu} -
\theta^{\dnu}}{u^{\mu}-u^{\nu} - \theta^{\dmu}\theta^{\dnu}}
 [A_{\mu},A_{\nu}] \var)\\
&\ = \ 0.
\endalign
$$
Hence the equation (4.9) is satisfied. Analogously one checks
that (4.10) is valid as well. 

\smallskip

This completes the
proof of Theorem 5.6.1.

\bigskip

\newpage

\centerline{\bf Bibliography}

\medskip

[BM] K.~Behrend, Yu.~Manin. {\it Stacks of stable maps and
Gromov--Witten invariants.} Duke Math. J., 85:1 (1996), 1--60.

\smallskip

[D] B.~Dubrovin. {\it Geometry of 2D topological field theories.}
In: Springer LNM, 1620 (1996), 120--348.

\smallskip

[H] N.~Hitchin. {\it Frobenius manifolds (notes by D. Calderbank.)}
Preprint, 1996.

\smallskip

[KM] M.~Kontsevich, Yu.~Manin. {\it Gromov-Witten classes, quantum
cohomology, and enumerative geometry.} Comm. Math. Phys.,
164:3 (1994), 525--562.

\smallskip

[Mal1] B.~Malgrange. {\it D\'eformations de syst\`emes diff\'erentielles
et microdiff\'erenti\-elles.} In: S\'eminaire de l'ENS 1979--1982,
Progress in Math. 37, Birkh\"auser, Boston (1983), 353--379.

\smallskip

[Mal2] B.~Malgrange. {\it La clasification des connections irr\'egulieres
\`a une variable.} ibid, 381--399.

\smallskip

[Mal3] B.~Malgrange. {\it Sur les d\'eformations isomonodromiques.
I. Singularit\'es r\'eguli\`eres.}
ibid, 401--426.

\smallskip

[Mal4] B.~Malgrange. {\it Sur les d\'eformations isomonodromiques.
II. Singularit\'es irr\'eguli\`eres.} ibid, 427--438.

\smallskip

[Ma1] Yu.~Manin. {\it Frobenius manifolds, quantum cohomology,
and moduli spaces (Chapters I, II, III).} Preprint MPI 96--113, 111 pp.

\smallskip

[Ma2] Yu.~Manin. {\it Sixth Painlev\'e equation, universal elliptic curve,
and mirror of $\bold{P}^2$.} Preprint MPI 96--114 and alg--geom/9605010.

\smallskip

[Ma3] Yu.~Manin. {\it Gauge field theory and complex geometry.} Second
Edition, Springer, 1997.

\smallskip

[Ma4] Yu.~Manin. {\it Topics in Noncommutative geometry.}
Princeton University Press, 1991.

\smallskip

[P] I.~Penkov. {\it $\Cal D$-modules on supermanifolds}. Inv.\ Math.
71 (1983), 501-512.

\smallskip

[S]  C.~Sabbah. {\it Frobenius manifolds: isomonodromic
deformations and infinitesimal period mappings.}
Preprint, 1996.

\smallskip

[Sch]  L.~Schlesinger. {\it \"Uber eine Klasse von 
Differentialsystemen beliebiger Ordnung mit
festen kritischer Punkten.} J. f\"ur die reine
und angew. Math., 141 (1912), 96--145.

\enddocument